\documentclass[twoside,11pt]{article}

\usepackage{jmlr2e}

\usepackage{amsmath}
\usepackage{mathrsfs}
\usepackage{amssymb}
\usepackage{amssymb}
\usepackage{graphicx}
\usepackage{subfig}

\usepackage{enumitem}
\usepackage{epstopdf}
\usepackage{color}
\usepackage{multirow}
\usepackage{verbatim}
\usepackage{algorithm}
\usepackage{algorithmic} 
\usepackage{bm}
\usepackage{booktabs}

\ifx\assumption\undefined
\newtheorem{assumption}[theorem]{Assumption}
\fi
\newcommand{\argmin}{\mathop{\mathrm{argmin}}}

\newcommand{\be}{\begin{equation}}
\newcommand{\ee}{\end{equation}}
\newcommand{\ba}{\begin{array}}
\newcommand{\ea}{\end{array}}

\newcommand{\M}{\mathcal{M}}
\newcommand{\St}{\mathrm{St}}
\newcommand{\x}{X}
\newcommand{\st}{\textrm{s.t.}}
\newcommand{\Retr}{\textrm{Retr}}
\newcommand{\T}{\mathrm{T}}
\newcommand{\Id}{\mathrm{Id}}
\renewcommand{\v}{V}
\newcommand{\bzeta}{\mathbf{\zeta}}
\newcommand{\Pc}{\mathcal{P}}
\newcommand{\Gr}{\mathrm{Gr}}
\renewcommand{\O}{\mathcal{O}}

\newcommand{\grad}{\mathrm{{grad}}}

\jmlrheading{1}{2000}{1-48}{4/00}{10/00}{meila00a}{Wang, Ma and Xue}

\ShortHeadings{Riemannian Stochastic Proximal Gradient Methods}{Wang, Ma and Xue}
\firstpageno{1}

\begin{document}

\title{Riemannian Stochastic Proximal Gradient Methods for Nonsmooth Optimization over the Stiefel Manifold}

\author{\name Bokun Wang \email bokunw.wang@gmail.com \\
       \addr Department of Computer Science \\
				The University of Iowa \\
				Iowa City, IA 52242
       \AND
       \name Shiqian Ma \email sqma@ucdavis.edu \\
       \addr Department of Mathematics \\ 
       University of California \\
       One Shields Avenue \\
       Davis, CA 95616 
       \AND
       \name Lingzhou Xue \email lzxue@psu.edu \\
       \addr Department of Statistics \\ Pennsylvania State University \\
       University Park, PA 16802}

\editor{}

\maketitle


\begin{abstract}
Riemannian optimization has drawn a lot of attention due to its wide applications in practice. Riemannian stochastic first-order algorithms have been studied in the literature to solve large-scale machine learning problems over Riemannian manifolds. However, most of the existing Riemannian stochastic algorithms require the objective function to be differentiable, and they do not apply to the case where the objective function is nonsmooth. In this paper, we present two Riemannian stochastic proximal gradient methods for minimizing nonsmooth function over the Stiefel manifold. The two methods, named R-ProxSGD and R-ProxSPB, are generalizations of proximal SGD and proximal SpiderBoost in Euclidean setting to the Riemannian setting. Analysis on the incremental first-order oracle (IFO) complexity of the proposed algorithms is provided. Specifically, the R-ProxSPB algorithm finds an $\epsilon$-stationary point with $\O(\epsilon^{-3})$ IFOs in the online case, and $\O(n+\sqrt{n}\epsilon^{-2})$ IFOs in the finite-sum case with $n$ being the number of summands in the objective. Experimental results on online sparse PCA and robust low-rank matrix completion show that our proposed methods significantly outperform the existing methods that use Riemannian subgradient information. 
\end{abstract}

\begin{keywords} 
Riemannian Optimization, Stochastic Gradient Descent, SPIDER, Manifold Proximal Gradient Method, Online Sparse PCA
\end{keywords}

\section{Introduction}

We consider the following composite optimization problem over the Stiefel manifold 
$\M:=\St(d,r)=\{\x\in\mathbb{R}^{d\times r}\mid \x^\top \x = I_r\}$:
\begin{equation}\label{prob_1}
    \min_{X\in\mathcal{M}} F(X):= f(X) + h(X),
\end{equation}
where $f(X)$ takes one of the following two forms: 
\begin{itemize}
\item Online case:
    \begin{equation}\label{online-f} f(X) := \mathbb{E}_\pi[f(X; \pi)], \end{equation}
    where $\mathbb{E}_\pi$ is the expectation with respect to the random variable $\pi$.
\item Finite-sum case:
    \begin{equation}\label{finite-sum-f} f(X) :=  \frac{1}{n}\sum_{i=1}^n f_i(X), \end{equation}
    where $n$ denotes the number of data and is assumed to be extremely large.
\end{itemize}
Throughout this paper, we assume that $f(\cdot;\pi)$, $f_i(\cdot)$ and thus $f(\cdot)$ are all smooth, $h$ is convex and possibly nonsmooth. Here the smoothness and convexity are interpreted when the function in question is considered as a function in the ambient Euclidean space. Note that since \eqref{online-f} involves an expectation, and \eqref{finite-sum-f} involves extremely large $n$, we assume that the full gradient information of $f$ is not available and only stochastic estimators to the gradient of $f$ can be obtained.

Problem \eqref{prob_1} with $f$ being \eqref{online-f} and \eqref{finite-sum-f} appears frequently in machine learning applications. In the online case \eqref{online-f}, $f(X;\pi)$ denotes the loss function corresponding to data $\pi$; and in the finite-sum case \eqref{finite-sum-f}, $f_i(X)$ denotes the loss function corresponding to the $i$-th sample data. Function $h$ is usually a regularizer that can promote certain desired structure of the solution. For example, letting $h(X)=\|X\|_1:=\sum_{ij}|X_{ij}|$ serves the purpose of promoting the sparsity of solution $X$.

One important application of \eqref{prob_1} in the online case is the online sparse PCA, which can be cast as
\begin{equation}\label{online-spca}
    \min_{\x} \quad \mathbb{E}_{Z\in\mathcal{D}}[\|Z - \x\x^\top Z\|_2^2] + \mu\|\x\|_1, \st, \x\in \mathcal{M},
\end{equation}
where $\mu>0$ is a weighting parameter, $\mathcal{D}$ denotes the distribution of the random online data $Z$, and the $\ell_1$ norm is used to promote the sparsity of the eigenvectors. In this case, $r$ is the desired number of principal components.
{
For PCA, each principal component is a linear combination of all variables, and it is usually difficult to interpret the derived principal components, especially when the dimension is high. Simple thresholding is an ad hoc way to estimate sparse loadings for better interpretability, but it may result in misleading results in various respects \citep{cadima1995loading}. By solving a manifold optimization problem, sparse PCA estimates sparse loadings to achieve a good balance between dimension reduction and interpretability. Sparse PCA has been widely used in many research fields such as medical imaging, ecology, and neuroscience. In the landmark-based shape analysis of the CC brain structure, \cite{sjostrand2007sparse} found that sparse PCA is useful to derive localized and interpretable patterns of variability while PCA did not provide much interpretational value. \cite{gravuer2008strong} applied sparse PCA to perform the dimension reduction before fitting the aggregated boosted trees model, and the sparsity helps the interpretability of their model. Recently, \cite{baden2016functional} used sparse PCA to study the functional diversity of mouse retinal ganglion cells through a clustering framework and found that SPCA leads to better cluster quality than PCA.} Although PCA and sparse PCA have been studied extensively in the literature, studies for online sparse PCA, i.e., sparse PCA with streaming data, seem to be very limited \citep{yang2015streaming,online-spca-wang-lu-2016}. In this paper, we propose efficient stochastic Riemannian algorithms for solving this important application.

\subsection{Related Works}
Riemannian optimization has been an active research area in the last decade, due to its wide applications in machine learning, signal processing, statistics and so on. The monograph by \citet{absil2009optimization} studied optimization algorithms on matrix manifolds in depth. 
{Recently, Riemannian optimization with nonsmooth objective has attracted a lot of attention due to its applications in sparse PCA \citep{Jolliffe2003}, compressed modes in physics \citep{Ozolins2013}, unsupervised feature selection \citep{Yang2011,tang2012unsupervised}, sparse blind deconvolution \citep{Zhang-cvpr-2017}, to name just a few. Many deterministic algorithms for solving Riemannian optimization with nonsmooth objective have been studied recently, including Riemannian subgradient method \citep{li2019nonsmooth}, manifold proximal gradient method (ManPG) \citep{chen2018proximal}, Riemannian proximal gradient method \citep{Huang-Wei-2019}, manifold proximal point algorithm \citep{chen2020manifold}, manifold proximal linear algorithm \citep{Wang-ManPL-2021} and so on.}
When the loss function $f$ takes the expectation or finite-sum form as in \eqref{online-f} and \eqref{finite-sum-f}, stochastic algorithms are usually in demand because we have only access to noisy stochastic gradients of $f$ instead of the full gradient. When the nonsmooth regularizer $h$ vanishes, that is, when \eqref{prob_1} reduces to a smooth problem with $f$ given by \eqref{online-f} or \eqref{finite-sum-f}, there exist stochastic algorithms for solving it. {In particular, R-SGD \citep{bonnabel2013stochastic}, R-SVRG \citep{zhang2016first}, R-SRG \citep{kasai2018riemannian} and R-SPIDER \citep{zhou2018faster,zhang2018r} can all be used to solve it}. Among these algorithms, R-SVRG, R-SRG and R-SPIDER all utilize the variance reduction techniques \citep{johnson2013accelerating,SAGA-2014} to improve the convergence rate of R-SGD. On the other hand, when the nonsmooth regularizer $h$ presents but the manifold constraint vanishes in \eqref{prob_1}, i.e., when $\mathcal{M}$ is the Euclidean space, there exist stochastic proximal gradient algorithms for solving these unconstrained problems in Euclidean space. Popular methods include ProxSGD \citep{rosasco2014convergence}, ProxSVRG \citep{xiao2014proximal}, ProxSARAH \citep{pham2019proxsarah} and ProxSpiderBoost \citep{wang2018spiderboost}. However, to the best of our knowledge, when both nonsmooth regularizer $h$ and manifold constraint $X\in\mathcal{M}$ present as in \eqref{prob_1}, there is no stochastic algorithm that can solve them.  In this paper, we close this gap by proposing two stochastic algorithms, namely R-ProxSGD and R-ProxSPB, for solving \eqref{prob_1} with $f$ being \eqref{online-f} or \eqref{finite-sum-f}, i.e., Riemannian optimization with nonsmooth objectives. Our algorithms are inspired by the ManPG algorithm that is recently proposed by \citet{chen2018proximal} for solving the nonsmooth Riemannian optimization problem \eqref{prob_1}. ManPG assumes that the full gradient of $f$ can be obtained, and thus it is a deterministic algorithm, while our R-ProxSGD and R-ProxSPB are the first stochastic algorithms for solving \eqref{prob_1} without using subgradient information. Recently, \citet{li2019nonsmooth} showed that when the objective function is weakly convex, Riemannian stochastic subgradient Method (R-Subgrad) has $\O(\epsilon^{-4})$ iteration complexity for obtaining an $\epsilon$-stationary point. 

	\begin{table}[h]
		%
		\center
		\resizebox{\columnwidth}{!}{
		\begin{tabular}{ccc}
			\toprule
			Objective & Euclidean & Riemannian \\
			\midrule
			\multirow{5}*{Smooth}  & SGD \citep{nemirovski2009robust}  & R-SGD \citep{bonnabel2013stochastic} \\
			~ & SVRG \citep{johnson2013accelerating} & R-SVRG \citep{zhang2016first} \\
			~ & SARAH \citep{nguyen2017sarah} &  R-SRG \citep{kasai2018riemannian}\\
			~ & SPIDER \citep{fang2018spider} & R-SPIDER \citep{zhou2018faster,zhang2018r}\\ ~ & SpiderBoost \citep{wang2018spiderboost} & \textbf{R-SpiderBoost} (ours)
			 \\
			 \midrule 
			  & ProxSGD \citep{rosasco2014convergence} & \textbf{R-ProxSGD} (ours) \\
			Non- & ProxSVRG \citep{xiao2014proximal} & N/A \\
			smooth & ProxSARAH \citep{pham2019proxsarah} & N/A \\
			~ & ProxSpiderBoost \citep{wang2018spiderboost} & \textbf{R-ProxSPB} (ours) \\
			\bottomrule
		\end{tabular}}
		\vspace{-.2cm}
		\caption{Summary of existing methods and our methods in Euclidean and Riemannian settings.}
        \label{tabloo}
	\end{table}
	
\begin{table}[htbp]
\setlength\tabcolsep{3.5pt}
\center
 \begin{tabular}{lccl}
  \toprule
  Algorithms & Step size &  Finite-sum & Online\\
  \midrule
 ManPG \citep{chen2018proximal} & constant & $\O(n\epsilon^{-2})$ & N/A\\
 R-ProxSGD & {constant}   & N/A & {$\O(\epsilon^{-4})$}\\
 R-ProxSPB & constant & $\O(n+\sqrt{n}\epsilon^{-2})$ & $\O(\epsilon^{-3})$\\
  \bottomrule
 \end{tabular}
  \caption{Comparison of IFO complexity for nonsmooth Riemannian optimization methods over the Stiefel manifold. 
  }
  \label{theoretical_res}
\end{table}

\subsection{Our Contributions}
The contributions of this paper lie in several folds.
\begin{enumerate}[label=(\roman*)]
\item First, we propose two stochastic algorithms for solving \eqref{prob_1}. These two algorithms, named R-ProxSGD and R-ProxSPB, are Riemannian generalizations of their counterparts in the Euclidean setting: ManPG \citep{chen2018proximal} and ProxSpiderBoost \citep{wang2018spiderboost}. On the other hand, they can also be viewed as generalizations of their smooth counterparts, R-SGD and R-SpiderBoost, to the nonsmooth case. {However, we emphasize here that although the design of these algorithms are straightforward, proving their convergence is more involved, due to the presence of stochastic gradients information.} 
    In Table \ref{tabloo} we give a summary of existing methods and our proposed methods in different cases: the objective is smooth or nonsmooth and the constraint is Riemannian manifold or Euclidean space. Note that when the nonsmooth function $h$ vanishes, our R-ProxSPB reduces to a Riemannian SpiderBoost algorithm (R-SpiderBoost) that solves Riemannian optimization with smooth objective. It seems that R-SpiderBoost is also new in the literature. 
\item Second, we prove the convergence of the proposed two algorithms and analyze their incremental first-order oracle (IFO) complexity results. Specifically, we analyze the IFO complexity of R-ProxSGD for the online setting problem, i.e., \eqref{prob_1} with $f$ being \eqref{online-f}; and R-ProxSPB for both the online setting problem and the finite-sum setting problem, i.e., \eqref{prob_1} with $f$ being \eqref{finite-sum-f}. In Table \ref{theoretical_res} we summarize the IFO complexity results of our proposed algorithms and the existing ManPG algorithm, as they are the only algorithms that can solve the nonsmooth Riemannian optimization problem \eqref{prob_1} with known IFO complexity results.
\item Third, we conduct numerical experiments for solving online sparse PCA \eqref{online-spca} and robust low-rank matrix completion problems to demonstrate the advantages of the proposed methods. 
\end{enumerate}

\begin{remark}
{We provide some further remark about the proposed algorithms R-ProxSGD and R-ProxSPB. Our algorithms incorporated several concepts, including Riemannian algorithm, proximal algorithm, stochastic algorithm, and variance reduction. We point out that they are all well motivated and justified. Note that the problem \eqref{prob_1} has three items that need to be taken care of: the manifold constraint, the smooth function $f$ and the nonsmooth function $h$. First, to deal with the manifold constraint, a {\it Riemannian} algorithm needs to be adopted. Second, since we do not have access to the full gradient information of the smooth function $f$, we need to design a stochastic algorithm that utilizes the noisy gradient information only. Third, to handle the nonsmooth function $h$, we need to design a proximal algorithm. Last, the variance reduction technique is adopted to reduce the variance of the stochastic gradients, and thus to accelerate the convergence of the algorithm.}
\end{remark}

{\bf Organization.} The rest of the paper is organized as follows. Section \ref{prelim} introduces the necessary notation and assumptions. Our new algorithms and their convergence and complexity results are presented in Section \ref{approach}. The experimental results are reported in Section \ref{experiments}. Finally, we make some concluding remarks in Section \ref{conclu}. The detailed proofs of the theorems and lemmas are provided in the appendix.

\section{Preliminaries}\label{prelim}

In this work, we consider the Riemannian submanifold $(\mathcal{M}, \mathfrak{g})$ where $\mathcal{M}$ is the Stiefel manifold and $\mathfrak{g}$ is the Riemannian metric on $\mathcal{M}$ that is induced from the Euclidean inner product. {That is, for any $x\in \M$, $\xi \in \T_x\M$ and $\zeta\in \T_x\M$, we have $\langle \xi,\zeta\rangle_x = \langle \xi, \zeta \rangle$, where $\T_x\M$ denotes the tangent space of $\M$ at $x$.} For smooth function $f$, we use $\mathrm{grad} f(X)$ to denote the full Riemannian gradient of $f$ at $X$, and $\nabla f({X})$ represents the full Euclidean gradient of $f$ at $X$. With an abuse of notation, when there is no ambiguity, we use $f_i$ to denote the component function in the online case \eqref{online-f}, i.e., $f_i(X):=f(X;\pi_i)$, though it is still used as a component function in the finite-sum case \eqref{finite-sum-f}.
For a mini-batch set $\mathcal{S}$, $\nabla f_{\mathcal{S}} (\x):= \frac{1}{|\mathcal{S}|} \sum_{i\in\mathcal{S}} \nabla f_i(\x)$ denotes the stochastic Euclidean gradient estimated on $\mathcal{S}$. We use $\mathcal{F}_t$ to denote all randomness occurred up to (include) the $t$-th iteration of any algorithm. When there is no ambiguity, we use $\|\mathbf{a}\|$ to denote the Frobenius norm when $\mathbf{a}$ is a matrix and the Euclidean norm when $\mathbf{a}$ is a vector.

{A classical geometric concept in the study of manifolds is the exponential mapping, which defines a geodesic curve on the manifold. However, the exponential mapping is difficult to compute in general. The concept of a retraction \citep{absil2009optimization}, which is a first-order approximation of the exponential mapping and can be more amenable to computation, is given as follows.}
\begin{definition}\label{def_retraction}\cite[Definition 4.1.1]{absil2009optimization}
A retraction on a differentiable manifold $\mathcal{M}$ is a smooth mapping $\Retr$ from the tangent bundle $\T\mathcal{M}$ onto $\mathcal{M}$ satisfying the following two conditions (here $\Retr_\x$ denotes the restriction of $\Retr$ onto $\T_\x \mathcal{M}$):
\begin{enumerate}
\item $\Retr_\x(0)=\x, \forall \x\in\mathcal{M}$, where $0$ denotes the zero element of $\T_X\mathcal{M}$.
\item For any $\x\in\M$, it holds that
    \[\lim_{\T_\x\M\ni\xi\rightarrow 0}\frac{\|\Retr_\x(\xi)-(\x+\xi)\|}{\|\xi\|} = 0.\]
\end{enumerate}
\end{definition}
\begin{remark}
Here and thereafter, when we talk about the summation $\x+\xi$, we always treat $\x$ and $\xi$ as elements in the ambient Euclidean space so that their sum is well defined. The second condition in Definition~\ref{def_retraction} ensures that $\Retr_{\x}(\xi)  =X +\xi +\mathcal{O}(\|\xi\|^2)$ and $D\Retr_{\x}(0)=\Id$, where $D\Retr_{\x}$ is the differential of $\Retr_\x$ and $\Id$ denotes the identity mapping. For more details about retraction, we refer the reader to \cite{absil2009optimization,boumal2018global} and the references therein.
\end{remark}
The retraction onto the Euclidean space is simply the identity mapping; i.e., $\Retr_X(\xi)=X+\xi$. For the Stiefel manifold $\St(d,r)$, common retractions
include the exponential mapping \citep{EdelmanAriasSmith1999}
\[\Retr_X^{\mathrm{exp}}(\xi) = [X,Q] \exp\left(\begin{bmatrix}
-X^\top\xi & -R^\top \\
R&0
\end{bmatrix}\right)\begin{bmatrix}
I_r\\0
\end{bmatrix}, \]
where $QR=-(I_d-XX^\top)\xi$ is the unique QR factorization; the polar decomposition
\[\Retr_{X}^{\mathrm{polar}}(\xi)=(X+\xi)(I_r +\xi^\top\xi)^{-1/2};\]
the QR decomposition
\[\Retr_{X}^{\mathrm{QR}}(\xi)=\mathrm{qf}(X+\xi),\]
where $\mathrm{qf}(A)$ is the $Q$ factor of the QR factorization of $A$;
the Cayley transformation \citep{Wen-Yin-2013}
\[\Retr_{X}^{\mathrm{cayley}}(\xi)=\left(I_d-\frac{1}{2}W(\xi)\right)^{-1}\left(I_d+\frac{1}{2}W(\xi)\right)X,\]
where $W(\xi)=(I_d-\frac{1}{2}XX^\top)\xi X^\top-X\xi^\top(I_d-\frac{1}{2}XX^\top)$.

In this paper, we adopt the assumption that the retraction that we use is invertible, the same as what is assumed in existing works \citep{kasai2018riemannian,zhou2018faster}. {We use $\Gamma_\x^{Y}$ to denote the vector transport from $\x$ to $Y$ satisfying $\Retr_\x{(\xi)} = Y$. Vector transport $\Gamma: \T\M \bigoplus \T\M \rightarrow \T\M$, $(\xi, \zeta) \mapsto \Gamma_\x^Y(\zeta)$ is associated with the retraction $\Retr$, where $\xi, \zeta \in \T_\x\M$.}

The following assumptions regarding the retraction and vector transport are necessary to our analysis.

\begin{assumption}\label{ass_retr_all}
\begin{enumerate}[label=(\roman*)]
    \item (see \cite{kasai2018riemannian}). All of the iterates $\{\x_t\}_{t=1}^{T+1}$ are in a totally retractive neighborhood $\mathcal{U}\subset\M$ of an optimum $X^*$: $\{\mathrm{Retr}_{X_t}(\xi_t)\}\in\mathcal{U}$ with $X_{t+1} = \mathrm{Retr}_{X_t}(\zeta_t)$, $\zeta_t \in \T_{X_t}\M$.
    \item (see \cite{kasai2018riemannian}). Suppose that $\mathrm{Exp}_\x: \T_\x\M \rightarrow \M$ denotes the exponential mapping and $\mathrm{Exp}_\x^{-1}: \M \rightarrow \T_\x\M$ is its inverse mapping. There exist $c_R, c_E>0$ such that $\|\mathrm{Exp}_{X}^{-1}(Y) - \mathrm{Retr}_{X}^{-1}(Y)\|\leq c_R\|\mathrm{Retr}_{X}^{-1}(Y)\|$, $ \forall X, Y\in\mathcal{U}$ and $\|\Retr_X^{-1}(Y)\|\leq c_E\|\xi\|$ if $\Retr_X(\xi)=Y$.
    \item (see \cite{boumal2018global}). For all $X\in\mathcal{M}$ and $\xi\in \T_X\mathcal{M}$, there exist constants $M_1>0$ and $M_2>0$ such that the following two inequalities hold:
        \begin{eqnarray}\label{retr1}
            & \|\mathrm{Retr}_X(\xi) - X\| \leq M_1\|\xi\|\\\label{retr2}
            & \|\mathrm{Retr}_X(\xi) - (X+\xi)\| \leq M_2\|\xi\|^2.
        \end{eqnarray}
     \item {$f(\cdot;\pi)$, $f_i(\cdot)$ and $f(\cdot)$ are all twice continuously differentiable.}
\end{enumerate}
\end{assumption}

\begin{assumption}\label{ass_vec_trans}
(see \cite{kasai2018riemannian}). {The vector transport is isometric on the manifold $\M$, i.e., $\|\Gamma_{X}^{Y}(\zeta)\|=\|\zeta\|$ for $\x, Y\in \M$, $\xi,\zeta \in \T_\x\M$ and $\Retr_\x{(\xi)} = Y$.}
\end{assumption} 

Besides, we impose some assumptions on $f(X)$ and its first-order oracle, which are also required in previous work on smooth Riemannian optimization with retraction and vector transport \citep{kasai2018riemannian, zhou2018faster}.

\begin{assumption}[Upper-bounded Hessian of $f$] \label{ass_upper_bounded_hess}
Every individual loss $f_i(X)$ is twice continuously differentiable and the individual Hessian of every $f_i(X)$ is bounded as $\|\nabla^2 f_i(X)\|\leq L_H$. $f(X)$ has upper-bounded Hessian in $\mathcal{U}\in\M$ with respect to the retraction $\mathrm{Retr}_X(\cdot)$ if there exists $L_R>0$ such that $\frac{d^2f(\mathrm{Retr}_\x(t\xi))}{dt^2}\leq L_R$ for all $X\in\mathcal{U}, \xi\in \T_X\M$ with $\|\xi\|=1$ and all $t$ such that $\mathrm{Retr}_\x(\tau\xi) \in\mathcal{U}$ for all $\tau\in[0,t]$.
\end{assumption}

\begin{assumption}[Bounded variance]\label{ass_bounded_variance}
Stochastic gradient oracle of every individual loss $f_i(X)$ is bounded $\|\nabla f_i(X)\|\leq G$ and its variance is also bounded $\mathbb{E}_i[\|\nabla f_i(X) - \nabla f(X)\|^2]\leq \sigma^2$.
\end{assumption}

Moreover, we make the following assumption on the regularization term $h(\x)$.

\begin{assumption}\label{ass_h}
The regularization function $h$ is convex and $L_h$-Lipschitz continuous, i.e., $\|h(\x)-h(Y)\|\leq L_h\|\x-Y\|$, $\forall \x,Y\in\M$.
\end{assumption}

We now give the definition of the stationary point of problem \eqref{prob_1}, which is standard in the literature, see \citep{yang2014optimality,chen2018proximal}. 

\begin{definition}[Stationary point]\label{stationary_1}
$X\in \M$ is a stationary point of \eqref{prob_1} if it satisfies:
\begin{equation}\label{epsilon_stationary_1_def}
    0\in \hat{\partial} F(\x) := {\grad f(\x)} + \mathrm{Proj}_{\T_X\M}{\partial} h(X),
\end{equation}
\end{definition}
where $\grad f(\x)$ is the Riemannian gradient of $f$ at $\x$, and $\hat{\partial} F(X)$ is the generalized Clarke subdifferential at $\x$ (see Definition \ref{def-clarke} in Appendix).

The computational costs of the algorithms are evaluated in terms of IFO complexity.
\begin{definition}\label{oracle_complex_def}
An IFO takes an index $i\in \{1,\ldots,n\}$ and returns $(f_i(X), \nabla f_i(X))$ for the finite-sum case \eqref{finite-sum-f}, or $(f(X;\pi_i), \nabla_\x f(X;\pi_i))$ for the online case \eqref{online-f}.
\end{definition}

\section{Riemannian Stochastic Proximal Gradient Methods}\label{approach}

In this section, we introduce our Riemannian stochastic proximal gradient algorithms and provide their non-asymptotic convergence results. Proofs of the theorems are provided in the appendix.

\subsection{The Main Framework}\label{subproblem_sec}

The main framework of our Riemnannian stochastic proximal gradient algorithms is inspired by the ManPG algorithm \citep{chen2018proximal}. The ManPG algorithm aims to solve the nonsmooth Riemannian optimization problem \eqref{prob_1} by assuming that the full gradient of $f$ can be accessed. Therefore, it is a deterministic algorithm. ManPG is a generalization of the proximal gradient method from Euclidean setting to the Riemannian setting. The proximal gradient method for solving $\min_X F(X):=f(X)+h(X)$ in the Euclidean setting generates the iterates as follows:
\begin{equation}\label{alg:proximal_euclidean}
X_{t+1}:=\argmin_Y f(X_t) +\langle\nabla f(X_t),Y-X_t\rangle + \frac{1}{2\gamma}\|Y-X_t\|^2 + h(Y).
\end{equation}
In other words, one minimizes the quadratic function $Y \mapsto f(X_t) +\langle\nabla f(X_t),Y-X_t\rangle + \frac{1}{2\gamma}\|Y-X_t\|^2 + h(Y)$ of $F$ at $X_t$ in the $t$-th iteration, where $\gamma>0$ is a parameter that can be regarded as the stepsize. {It is known that this quadratic function can bound $F$ from above when $\gamma\leq 1/L$, where $L$ is the Lipschitz constant of $\nabla f$.}
The subproblem \eqref{alg:proximal_euclidean} corresponds to the proximal mapping of $h$ and the efficiency of the proximal gradient method relies on the assumption that \eqref{alg:proximal_euclidean} is easy to solve.
For \eqref{prob_1}, in order to deal with the manifold constraint, one needs to ensure that the descent direction lies in the tangent space. This motivates the following subproblem for finding the descent direction $\xi_t$ in the $t$-th iteration:
\begin{eqnarray}\label{subproblem-ManPG}
\begin{array}{ll}
      \mathbf{\xi}_t = & \argmin_\mathbf{\xi}:=  \langle \nabla f(\x_t), \mathbf{\xi}\rangle + \frac{1}{2\gamma}\|\mathbf{\xi}\|^2 + h(X_t +\mathbf{\xi}) \\
      & \mathrm{s.t.}\quad \mathbf{\xi}\in \T_{X_t}\mathcal{M},
\end{array}
\end{eqnarray}
and then a retraction step is performed to keep the iterate feasible to the manifold constraint:
\begin{equation}\label{retraction-ManPG}
X_{t+1} := \Retr_{\x_t}(\eta_t\xi_t).
\end{equation}
It is shown that the ManPG algorithm \eqref{subproblem-ManPG}-\eqref{retraction-ManPG} finds an $\epsilon$-stationary point of \eqref{prob_1} in $O(\epsilon^{-2})$ iterations. It was shown in \citep{chen2018proximal} that ManPG performs better than some existing algorithms for solving the sparse PCA problem. The ManPG algorithm was extended successfully later to solving problems with two block variables \citep{A-ManPG-2019} such as another sparse PCA formulation \citep{Zou-spca-2006} and the sparse CCA problem \citep{hardoon2011sparse}.

Motivated by the success of the ManPG algorithm, when we only have the access to stochastic gradient of $f$, we design a stochastic version of ManPG to solve \eqref{prob_1}. In particular, each iteration of our proposed algorithm consists of two steps: (i) finding the descent direction, and (ii) performing retraction. The basic framework of our proposed algorithm is to simply replace the full gradient in ManPG by a stochastic estimator to the gradient. This leads to the following updating scheme of the proposed framework: 
\begin{eqnarray}\label{subproblem}
\begin{array}{ll}
      \mathbf{\zeta}_t = & \argmin_\mathbf{\zeta} \phi_t(\mathbf{\zeta}) \!:=\!  \langle \v_t, \mathbf{\zeta}\rangle + \frac{1}{2\gamma}\|\mathbf{\zeta}\|^2 + h(X_t +\mathbf{\zeta})\\
      & \mathrm{s.t.}\quad \mathbf{\zeta}\in \T_{X_t}\mathcal{M},
\end{array}
\end{eqnarray}
and
\begin{equation}\label{retraction}
X_{t+1} := \Retr_{\x_t}(\eta_t\bzeta_t),
\end{equation}
where $\gamma>0$ and $\eta_t>0$ are step sizes, and $\v_t$ denotes a stochastic estimation of the Euclidean gradient $\nabla f(\x_t)$. Specific choices of $\v_t$ will be discussed in Sections \ref{sec:R-ProxSGD} and \ref{sec:R-SpiderBoost}. 
Note that for the Stiefel manifold $\mathcal{M}$, the tangent space is given by $\T_\x\M = \{\zeta\mid \zeta^\top\x+\x^\top\zeta=0\}$. Therefore, the constraint in \eqref{subproblem} is a linear equality constraint. Since we assume that $h$ is a convex function, it follows that the subproblem \eqref{subproblem} is a convex problem. This convex problem can be efficiently solved using the semi-smooth Newton method \citep{xiao2018regularized}. We refer the readers to \cite{xiao2018regularized} and \cite{chen2018proximal} for more details on how to solve \eqref{subproblem} efficiently.

To prepare for the analysis of IFO complexity, we need to define the $\epsilon$-stationary solution and the $\epsilon$-stochastic stationary point.
\begin{definition}[$\epsilon$-stationary point and $\epsilon$-stochastic stationary point]\label{epsilon_stat_point}
Define
\begin{equation}\label{def-G}G(\x, \nabla f(\x), \gamma) = (\x - \mathrm{Retr}_{\x}(\xi))/\gamma,\end{equation}
where
\begin{equation}\label{def-xi-after-G}\xi := \argmin_{\xi \in \T_\x\mathcal{M}} \{\langle \nabla f(X), \xi\rangle + \frac{1}{2\gamma}\|\xi\|^2 + h(\x +\xi)\}.
\end{equation}
$\x$ is called an $\epsilon$-stationary point of \eqref{prob_1} if $\|G(\x, \nabla f(\x), \gamma)\|\leq \epsilon$.
When the sequence $\{\x_t\}$ is generated by a stochastic algorithm (stochastic process), we call $\x_t$ an $\epsilon$-stochastic stationary point if $\mathbb{E}[\|G(\x_t, \nabla f(\x_t), \gamma)\|]\leq\epsilon$, where the expectation $\mathbb{E}$ is taken for all randomness before $\x_t$ is generated.
\end{definition}

\begin{remark}
Note that $\xi$ defined in \eqref{def-xi-after-G} is the solution to \eqref{subproblem} with full gradient $\v_t = \nabla f(\x_t)$. In the Euclidean space, $\Retr_{\x_t}(\gamma \xi_t)$ reduces to $\x_t + \gamma \xi_t$ and $\xi_t = \mathrm{prox}_{\gamma h}(\x_t - \gamma\nabla f(\x))-\x$. Thus, $G(\x_t, \nabla f(\x_t), \gamma)$ defined in \eqref{def-G} is analogous to the proximal gradient in the Euclidean space.
\end{remark}

\subsection{R-ProxSGD: Riemannian Stochastic Proximal Gradient Descent Algorithm}\label{sec:R-ProxSGD}

In this section, we design the basic Riemannian proximal stochastic gradient descent method (R-ProxSGD) by choosing $\v_t$ as the mini-batch stochastically sampled gradients. More specifically, in the $t$-th iteration of R-ProxSGD, we randomly sample a mini-batch set $\mathcal{S}_t$, and define $\v_t =  \frac{1}{|\mathcal{S}_t|}\sum_{{i_t}\in\mathcal{S}_t} \nabla f_{i_t}(X_{t})$, which is an unbiased gradient estimator with bounded variance. That is, $\mathbb{E}[\v_t] = \nabla f(X_t)$ and $\mathbb{E}[\|\v_t-\nabla f(X_t)\|^2] \leq \frac{\sigma^2}{|\mathcal{S}_t|}$.
Our R-ProxSGD is described in Algorithm \ref{prox-r-sgd-algo-chart}. 

\begin{algorithm}[tb]
		\caption{R-ProxSGD}\label{prox-r-sgd-algo-chart}
		\begin{algorithmic}[1]
			\STATE \textbf{Input:} initial point $X_0\in\mathcal{M}$, parameters $\eta\in (0,1)$, {\begin{equation}\label{def-gamma}\gamma = \frac{2\eta}{2\tilde{L}\eta^2+\eta+1}, \quad \mbox{where} \quad \tilde{L} = L_R/2+L_hM_2.\end{equation}}
			\FOR{$t=0,1,..., T-1$}
			\STATE Compute the stochastic gradient by randomly sampling a mini-batch set $\mathcal{S}_t$ and calculating the unbiased stochastic gradient estimator:
			\begin{equation}\nonumber
			    \v_t = \nabla f_{\mathcal{S}_t}(\x_t) :=  \frac{1}{|\mathcal{S}_t|}\sum_{{i_t}\in\mathcal{S}_t} \nabla f_{i_t}(X_{t})
			\end{equation}
			\STATE Proximal step: obtain $\mathbf{\zeta}_t$ by solving the subproblem \eqref{subproblem}.
			\STATE Retraction step: $X_{t+1} = \mathrm{Retr}_{X_t} (\eta_t\mathbf{\zeta}_t)$, with {$\eta_t:=\eta$}.
			\ENDFOR
			\STATE \textbf{Output:} $X_{\nu}$, where $\nu$ is uniformly sampled from $\{1,...,T\}$. 
		\end{algorithmic}
	\end{algorithm}

We have the following iteration and IFO complexity results for R-ProxSGD for solving the online case problem \eqref{prob_1} with $f$ being \eqref{online-f}. The proof is given in the appendix.

\begin{theorem}\label{theorem_prox_r_sgd}
In R-ProxSGD, we set the batch size $|\mathcal{S}_t| := s = \O(\epsilon^{-2})$ for all $t$, and $\gamma$ is chosen as in \eqref{def-gamma}. Under this parameter setting, the number of iterations needed by R-ProxSGD for obtaining an $\epsilon$-stochastic stationary point of the online case problem \eqref{prob_1} with $f$ being \eqref{online-f}, is {$T=\O(\epsilon^{-2})$}. Moreover, the IFO complexity of the R-ProxSGD algorithm for obtaining an $\epsilon$-stochastic stationary point in the online setting \eqref{prob_1} with $f$ being \eqref{online-f} is {$\O(\epsilon^{-4})$}.
\end{theorem}

\begin{remark}
In Theorem \ref{theorem_prox_r_sgd}, since we require the batch size to be $\O(\epsilon^{-2})$, the results only hold for the online case problem, and do not hold for the finite-sum case problem.
\end{remark}

\subsection{R-ProxSPB: Riemannian Proximal SpiderBoost Algorithm}\label{sec:R-SpiderBoost}

Note that the convergence and complexity results of R-ProxSGD do not apply to the finite-sum case problem. In this section, we propose a Riemannian proximal SpiderBoost algorithm (R-ProxSPB) that can solve both the online case problem and the finite-sum case problem. More importantly, we can show that R-ProxSPB has an improved IFO complexity comparing with R-ProxSGD for the online case problem. For smooth problems in the Euclidean setting, there exist many works that use the variance reduction technique to improve the convergence speed of SGD, such as SVRG \citep{johnson2013accelerating}, SAGA \citep{SAGA-2014}, SARAH \citep{nguyen2017sarah}, SPIDER \citep{fang2018spider} and SpiderBoost \citep{wang2018spiderboost}. In particular, the SpiderBoost algorithm proposed by \citet{wang2018spiderboost} achieves the same complexity bound as SPIDER, but in practice SpiderBoost can converge faster because it allows a constant step size, while SPIDER requires an $\epsilon$-dependent step size that can be too conservative in practice.
Some of these algorithms have been extended to the Riemannian optimization with smooth objective functions, such as R-SVRG \citep{zhang2016first}, R-SRG \citep{kasai2018riemannian} and R-SPIDER \citep{zhang2018r,zhou2018faster}. It was found that R-SRG and R-SPIDER equipped with the biased R-SARAH estimator consistently outperform the R-SVRG algorithm. Inspired by the SpiderBoost algorithm, we propose a Riemannian proximal SpiderBoost algorithm, named R-ProxSPB, which is a generalization of SpiderBoost to nonsmooth Riemannian optimization. When the nonsmooth function $h$ vanishes, our R-ProxSPB algorithm reduces to a Riemannian SpiderBoost algorithm (R-SpiderBoost) for Riemannian optimization with smooth objective function, which seems to be new in the literature as well.

Our R-ProxSPB algorithm is described in Algorithm \ref{prox-r-spiderboost-algo-chart}. R-ProxSPB specifies a constant integer $q$. When the iteration number $t$ is a multiple of $q$, mini-batch $\mathcal{S}_t^1$ is sampled and unbiased stochastic gradient estimator is used; while for other iterations, mini-batch $\mathcal{S}_t^2$ is sampled and R-SARAH estimator \eqref{R-SARAH} is used. Comparing with R-ProxSGD (Algorithm \ref{prox-r-sgd-algo-chart}), a significant difference of R-ProxSPB is that it allows a constant step size $\eta$ instead of a diminishing step size.
That the constant step size is allowed is due to the biased stochastic gradient estimator R-SARAH, which leads to variance reduction of the stochastic gradients, and thus improves the convergence rate. This has been justified in several variance reduced stochastic algorithms such as SVRG, SAGA, SPIDER and SpiderBoost and so on. A constant step size usually leads to a faster algorithm both theoretically and practically. In fact, we can prove the following convergence rate and IFO complexity results of R-ProxSPB, which indeed improve the results of R-ProxSGD.

\begin{algorithm}[tb]
		\caption{R-ProxSPB}\label{prox-r-spiderboost-algo-chart}
		\begin{algorithmic}[1]
			\STATE \textbf{Input:} initial point $X_0\in\mathcal{M}$, parameters $\eta>0$, $\gamma>0$, integers $q$, $T$.
			\FOR{$t=0,..., T-1$}
			\IF{$\mathrm{mod}(t,q)=0$}
			\STATE Randomly sample a mini-batch $\mathcal{S}_t^1$ and calculate $\v_t = \nabla f_{\mathcal{S}_t^1}(X)$ satisfying:
            \[
			    \mathbb{E}[\|\v_t - \nabla f(X_t)\|^2]\leq \frac{\sigma^2}{|\mathcal{S}_t^1|}
			\]
			\ELSE
			\STATE Randomly sample a mini-batch $\mathcal{S}_t^2$ and calculate $\v_t$ by the R-SARAH estimator:
            \begin{equation}\label{R-SARAH}
                \v_t = \nabla f_{\mathcal{S}_t^2}(X_t) - \Gamma_{X_{t-1}}^{X_t}\big( \nabla f_{\mathcal{S}_t^2}(X_{t-1}) - \v_{t-1}\big)
            \end{equation}

	\ENDIF
	\STATE Proximal step: obtain $\mathbf{\zeta}_t$ by solving the subproblem \eqref{subproblem}.

    \STATE Retraction step: $X_{t+1} = \mathrm{Retr}_{X_t} (\eta\mathbf{\zeta}_t)$.
			\ENDFOR
			\STATE \textbf{Output:} $X_{\nu}$, $\nu$ is uniformly sampled from $\{1,...,T\}$.
		\end{algorithmic}
	\end{algorithm}

\begin{theorem}\label{theorem_spiderboost}
    In R-ProxSPB (Algorithm \ref{prox-r-spiderboost-algo-chart}), we set $\eta = \min\big(\frac{1}{2(L_R/2+L_hM_2)}, \frac{1}{\sqrt{2c_E \Theta^2}}\big)$, $\gamma=\frac{2}{5}$, and $|\mathcal{S}_t^2|=q$ for all $t$, where $\Theta$ is a constant that will be specified in the proof. Under this parameter setting, we have the following convergence rate and IFO complexity results of R-ProxSPB.
    \begin{itemize}
        \item[(i).] For the finite-sum case problem, i.e., problem \eqref{prob_1} with $f$ being \eqref{finite-sum-f}, we set $q=\sqrt{n}$,  $|\mathcal{S}_t^1| = n$, for all $t$. R-ProxSPB returns an $\epsilon$-stochastic stationary point of \eqref{prob_1} after $T = \O(\epsilon^{-2})$ iterations. Moreover, the IFO complexity of R-ProxSPB for obtaining an $\epsilon$-stochastic stationary point of \eqref{prob_1} is $\O(\sqrt{n}\epsilon^{-2}+n)$.
        \item[(ii).] For the online case problem, i.e., problem \eqref{prob_1} with $f$ being \eqref{online-f}, we set $q=\O(\epsilon^{-1})$, $|\mathcal{S}_t^1| = \O(\epsilon^{-2})$, for all $t$. R-ProxSPB returns an $\epsilon$-stochastic stationary point of \eqref{prob_1} after $T = \O(\epsilon^{-2})$ iterations. Moreover, the IFO complexity of R-ProxSPB for obtaining an $\epsilon$-stochastic stationary point of \eqref{prob_1} is $\O(\epsilon^{-3})$.
    \end{itemize}
\end{theorem}

\begin{remark}
Here we summarize some comparisons of the two proposed algorithms. For the online case problem, R-ProxSPB has a better IFO complexity than R-ProxSGD. R-ProxSPB allows constant step size $\eta$, but R-ProxSGD needs a diminishing step size $\eta_t$. The convergence results of R-ProxSPB in Theorem \ref{theorem_spiderboost} covers the finite-sum case problem, which is still lacking for the R-ProxSGD algorithm. We also need to point out that, though R-ProxSPB is faster than R-ProxSGD in theory, it involves more tuning parameters and the R-SARAH estimator might be difficult to compute for certain manifolds. Therefore, for certain applications, R-ProxSGD could be more favorable in practice.
\end{remark}

\section{Numerical Experiments}\label{experiments}


We compare our proposed algorithms R-ProxSGD and R-ProxSPB with several baselines on the online sparse PCA problem \eqref{online-spca}. The experiments are performed on two real datasets: \texttt{coil100} \citep{nene1996columbia} and \texttt{mnist} \citep{lecun1998mnist}. The \texttt{coil100} dataset contains $n=7,200$ RGB images of 100 objects taken from different angles with {$d=1024$}. The \texttt{mnist} dataset has $n=80,000$ grayscale digit images of size {$d=28\times28=784$}.

\subsection{Online Sparse PCA Problem}

\subsubsection{Comparison with Riemannian stochastic subgradient method}

First, we compare our proposed algorithms R-ProxSGD and R-ProxSPB with the Riemannian stochastic subgradient method (R-Subgrad). R-Subgrad for solving \eqref{online-spca} iterates as follows:
\[\begin{array}{ll}
\xi_t    & := -\mathrm{Proj}_{\x_t}(-2Z_{i_t}Z_{i_t}^\top\x_t + \mu~\mathrm{sign}(\x_t)), \\
\x_{t+1} & := \Retr_{\x_t}(\eta_t\xi_t),
\end{array}
\]
where $Z_{i_t}$ is a randomly sampled data. Here the projection operation is defined as: $\mathrm{Proj}_\x(Y) = Y - \x~\mathrm{sym}(\x^\top Y)$ and $\mathrm{sym}(\x) = \frac{1}{2}(\x+\x^\top)$.

For R-Subgrad and R-ProxSGD, we use the diminishing step size $\eta_t = \frac{\eta_0}{\sqrt{t+1}}$. For R-ProxSPB, we use the constant step size $\eta$ as suggested in our theory. Because some of the problem-dependent constants cannot be directly estimated from the datasets, we perform grid search to tune $\eta_0$ and $\eta$ for all algorithms from $\{5\times 10^{-5}, 10^{-4}, 5\times 10^{-4}, ..., 1\}$. The best $\eta_0$ and $\eta$ on different datasets are reported in the appendix. For R-ProxSGD, we set $|\mathcal{S}_t|=100$. For R-ProxSPB, we set $|\mathcal{S}_t^1|=n$ and $|\mathcal{S}_t^2|=q=100$. 

\begin{figure*}[htbp]
\minipage{0.33\textwidth}
  \includegraphics[width=\linewidth]{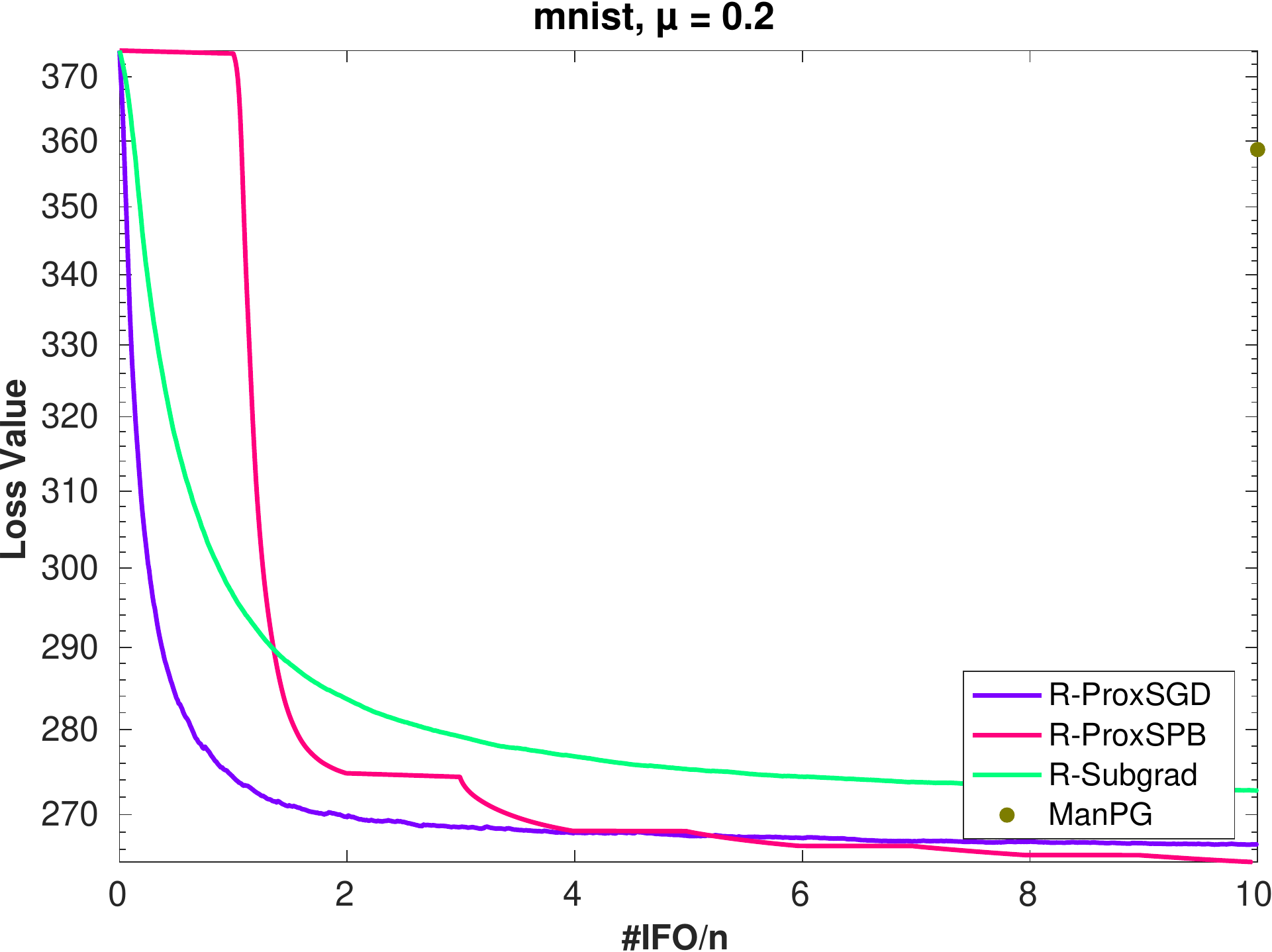}
\endminipage\hfill
\minipage{0.33\textwidth}
  \includegraphics[width=\linewidth]{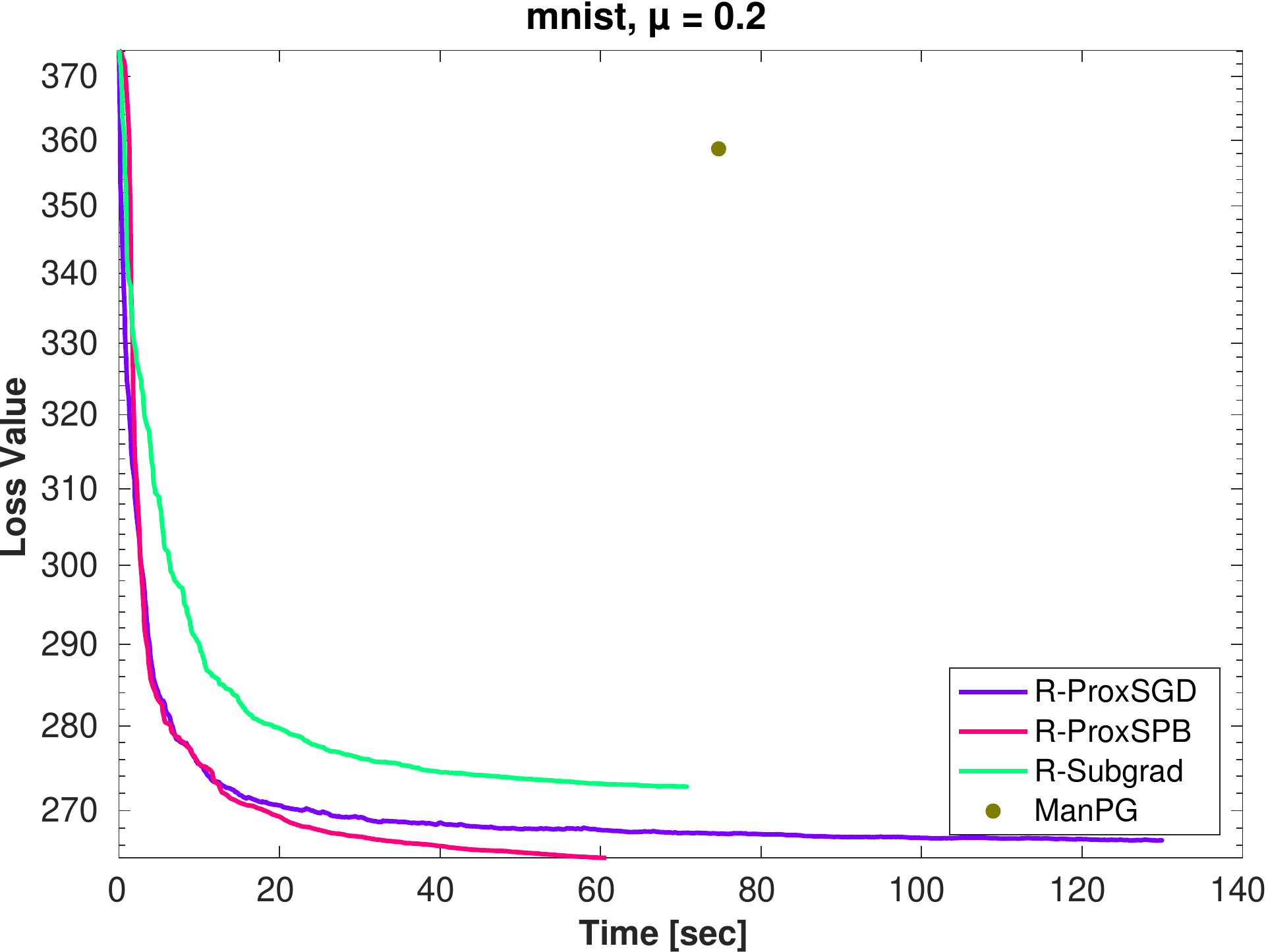}
\endminipage\hfill
\minipage{0.33\textwidth}%
  \includegraphics[width=\linewidth]{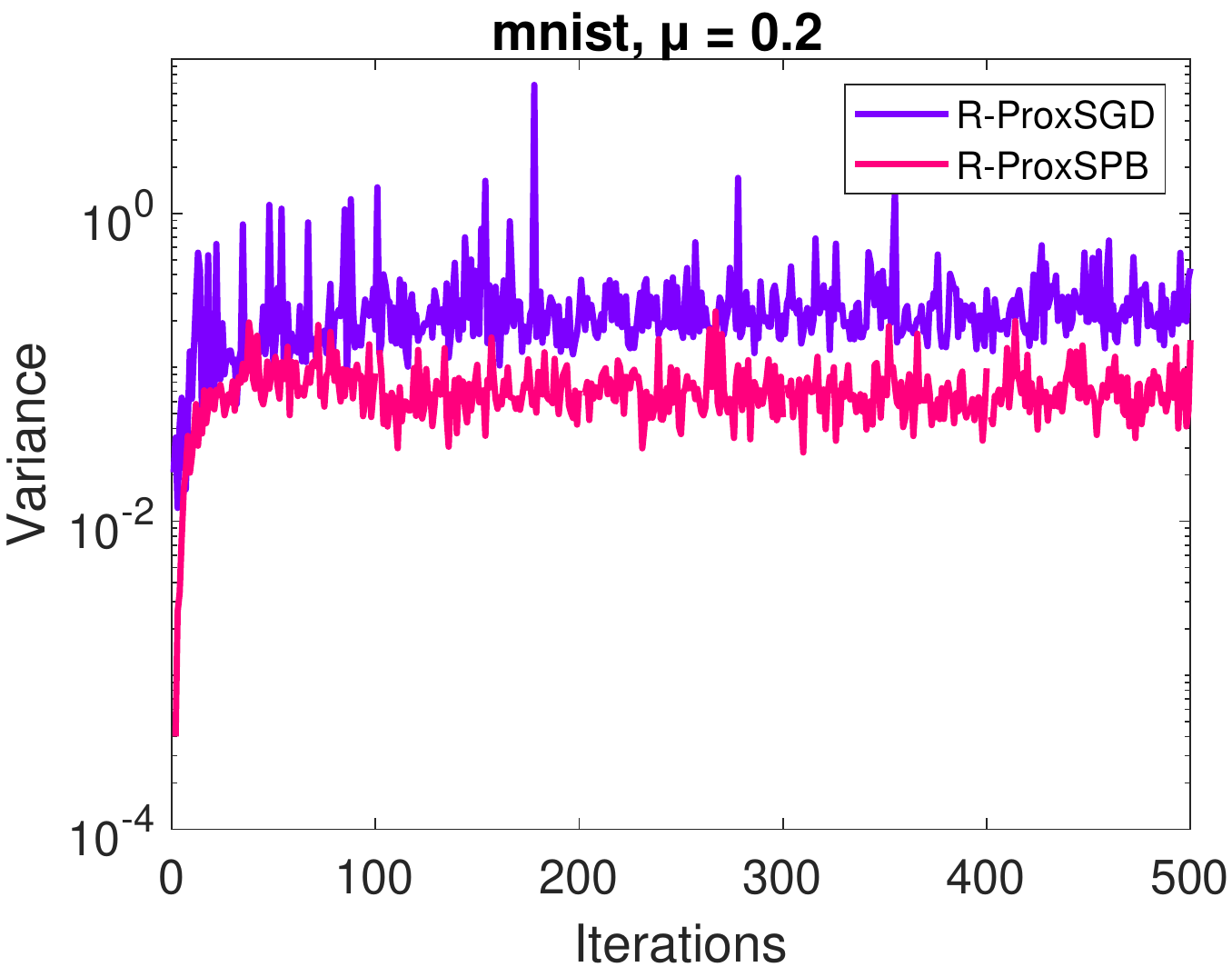}
\endminipage\hfill
\minipage{0.33\textwidth}
  \includegraphics[width=\linewidth]{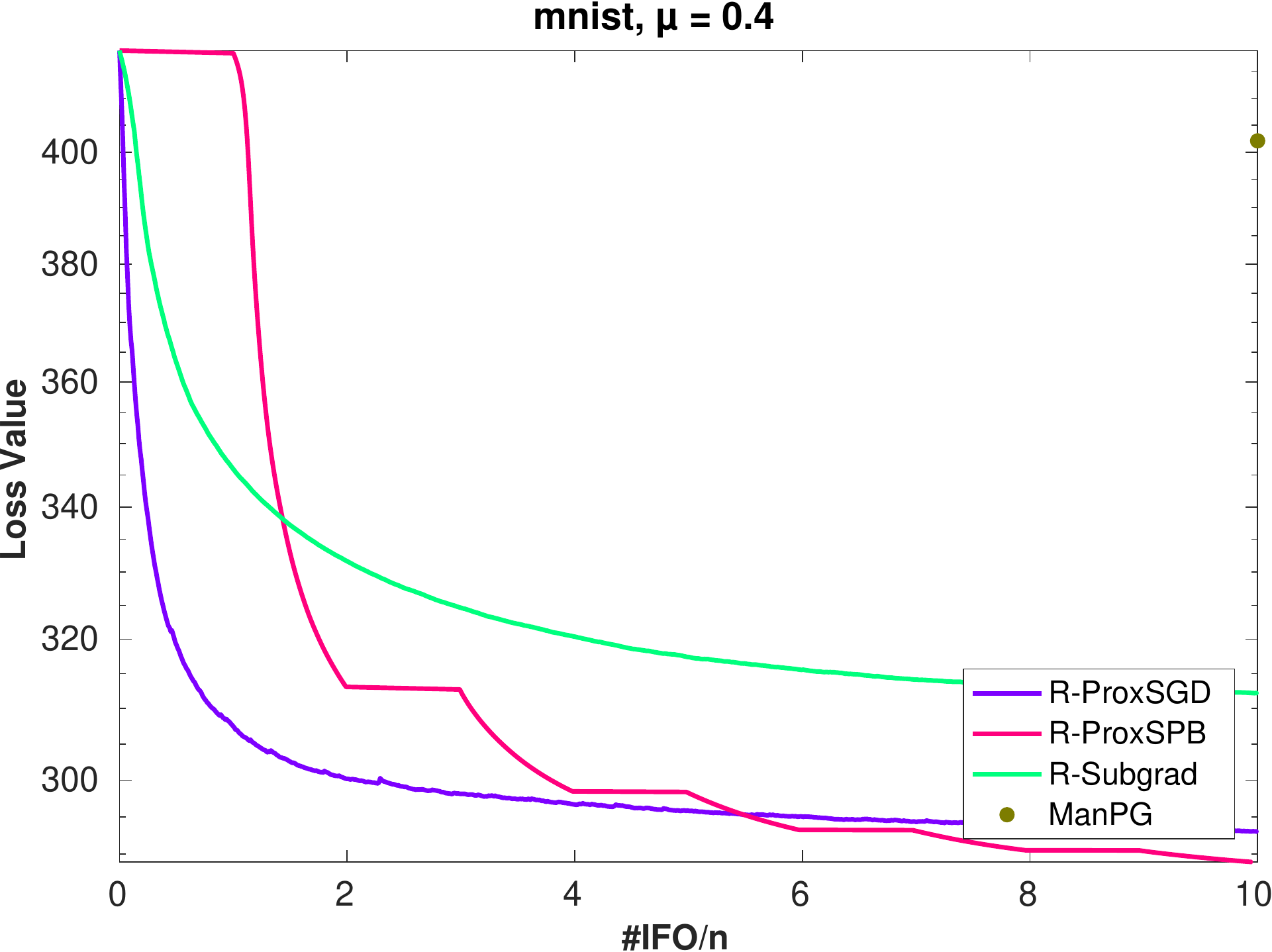}
\endminipage\hfill
\minipage{0.33\textwidth}
  \includegraphics[width=\linewidth]{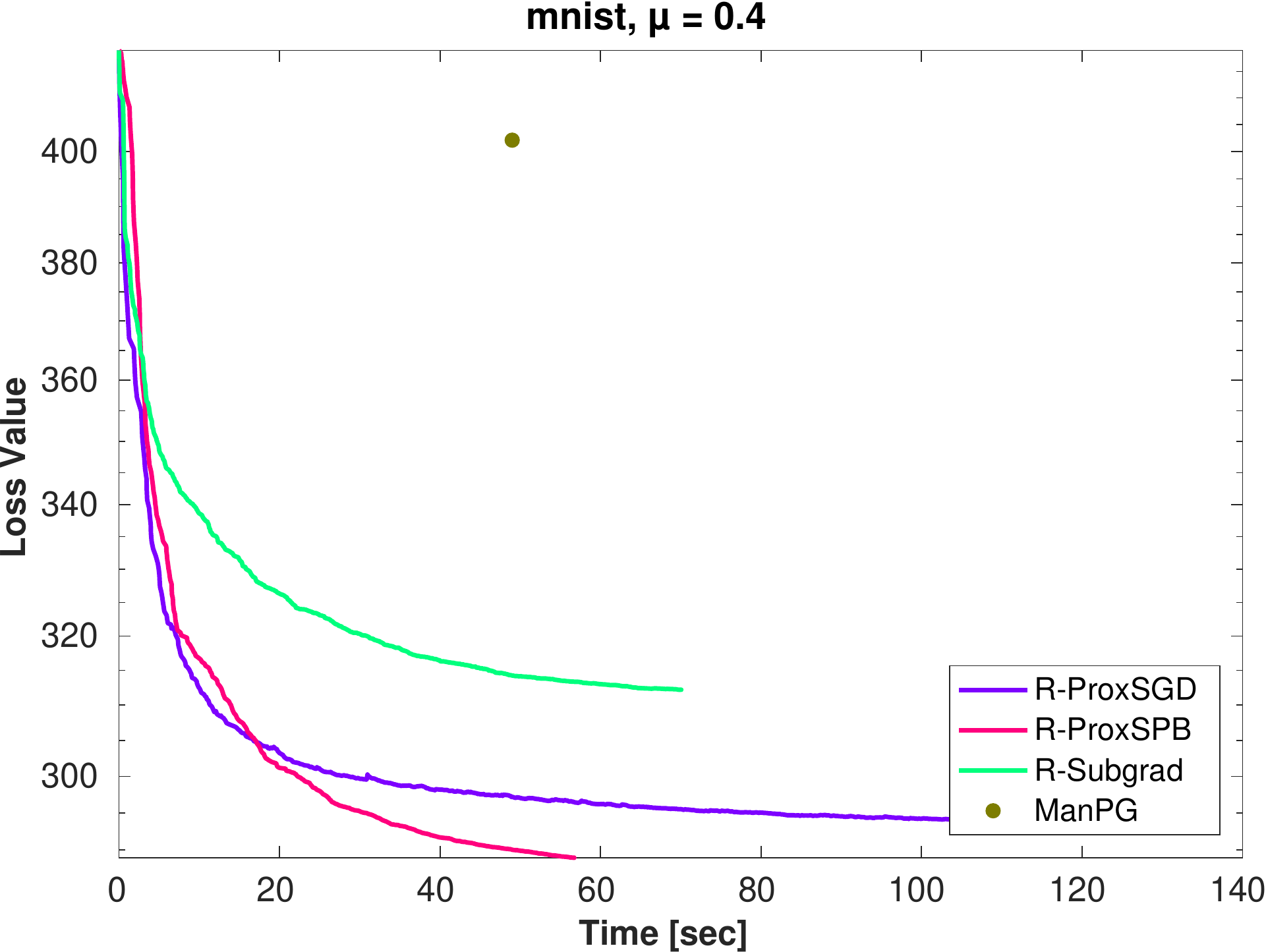}
\endminipage\hfill
\minipage{0.33\textwidth}%
  \includegraphics[width=\linewidth]{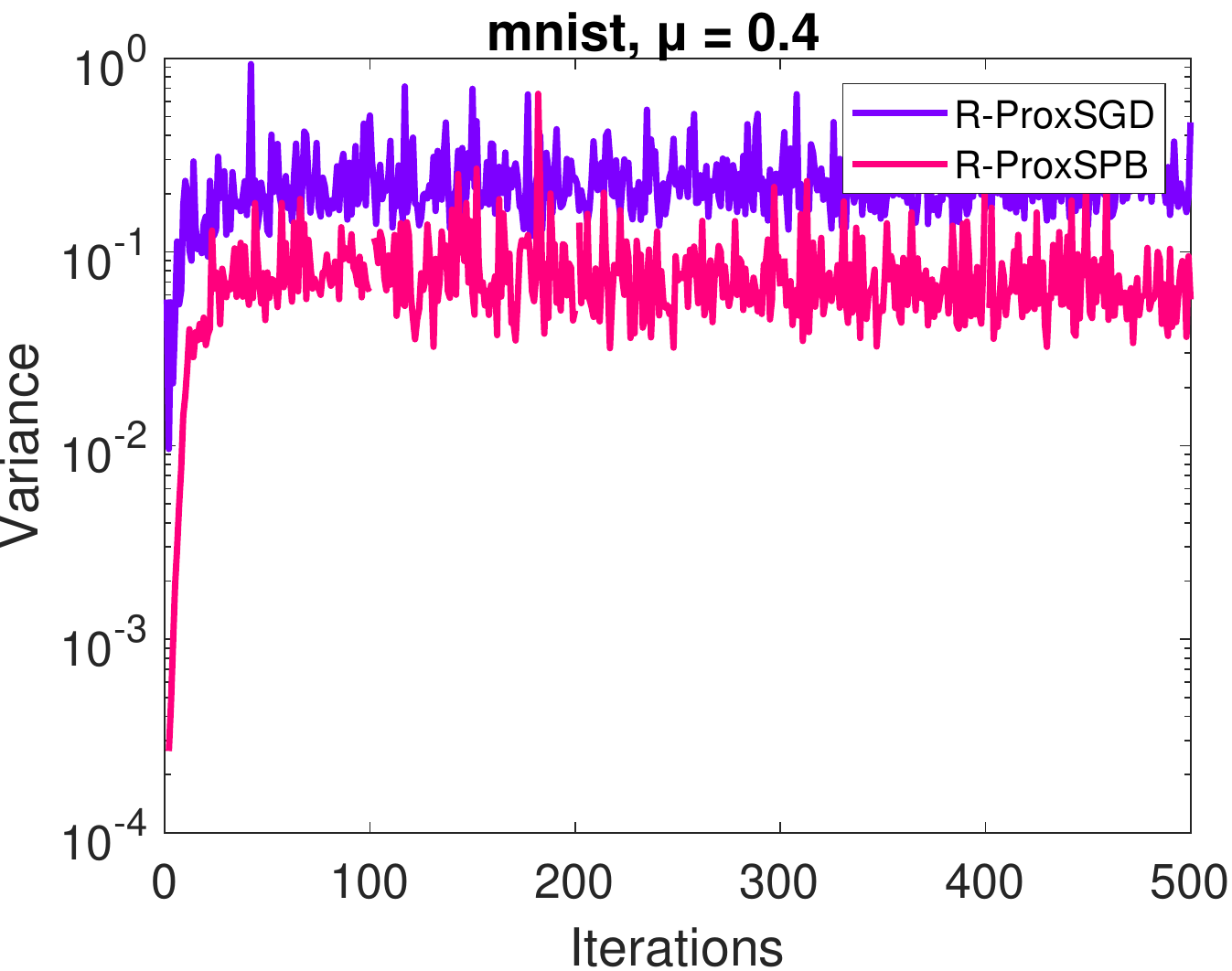}
\endminipage\hfill

\caption{Experimental results on the \texttt{mnist} dataset with $\mu = 0.2$ and  $0.4$.}
\label{fig1}
\end{figure*}

\begin{figure*}[htbp]
\minipage{0.33\textwidth}
  \includegraphics[width=\linewidth]{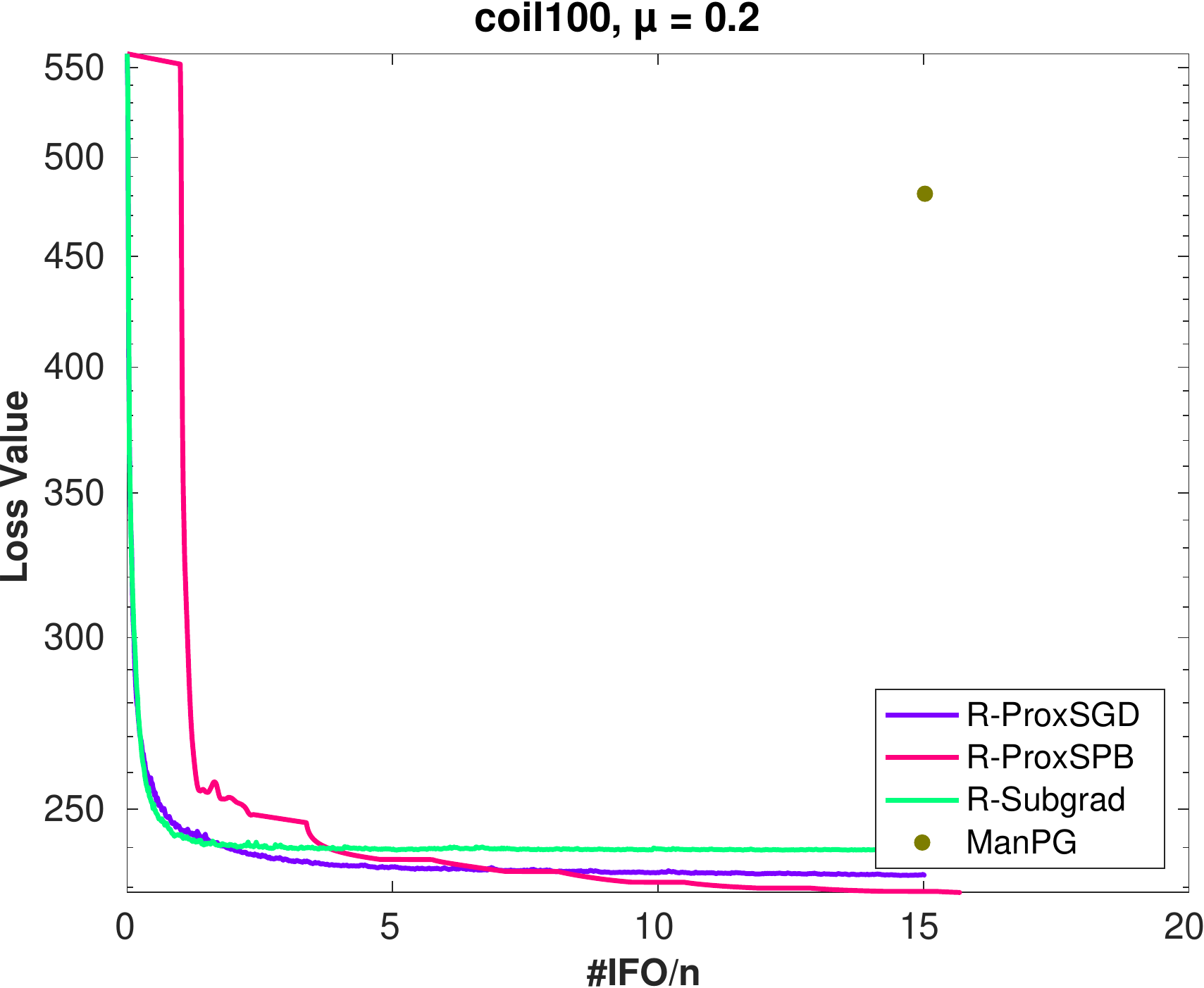}
\endminipage\hfill
\minipage{0.33\textwidth}
  \includegraphics[width=\linewidth]{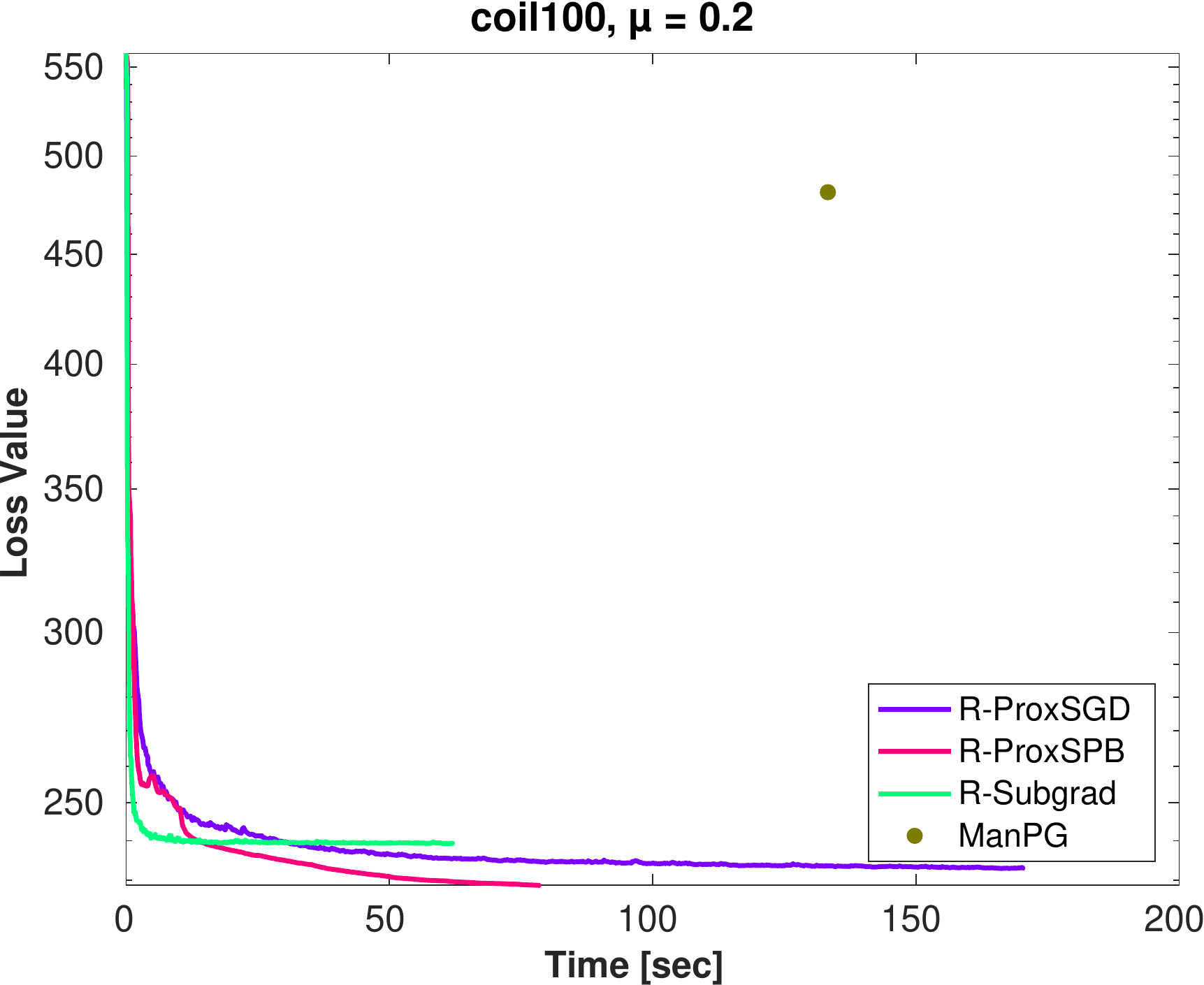}
\endminipage\hfill
\minipage{0.33\textwidth}%
  \includegraphics[width=\linewidth]{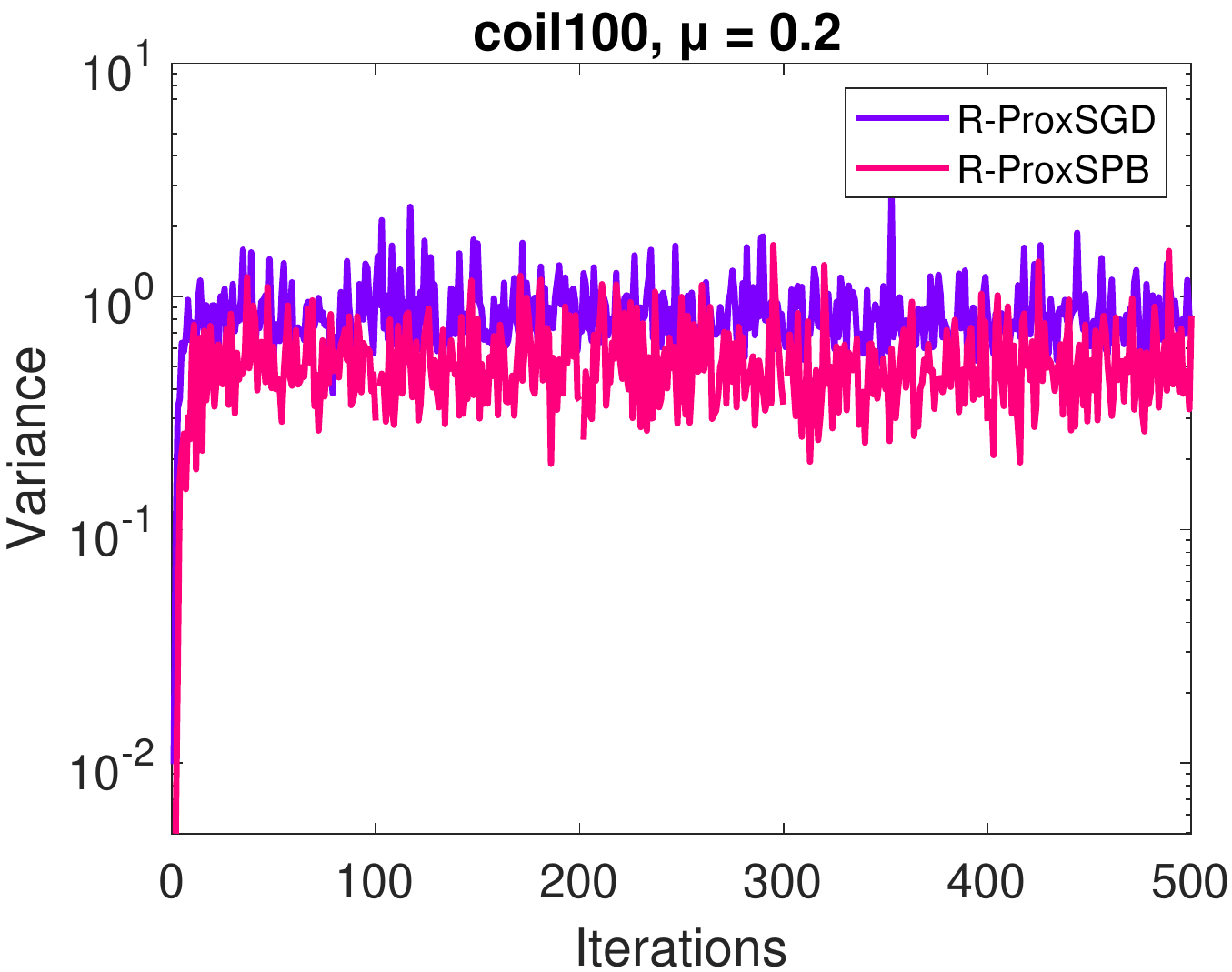}
\endminipage\hfill
\minipage{0.33\textwidth}
  \includegraphics[width=\linewidth]{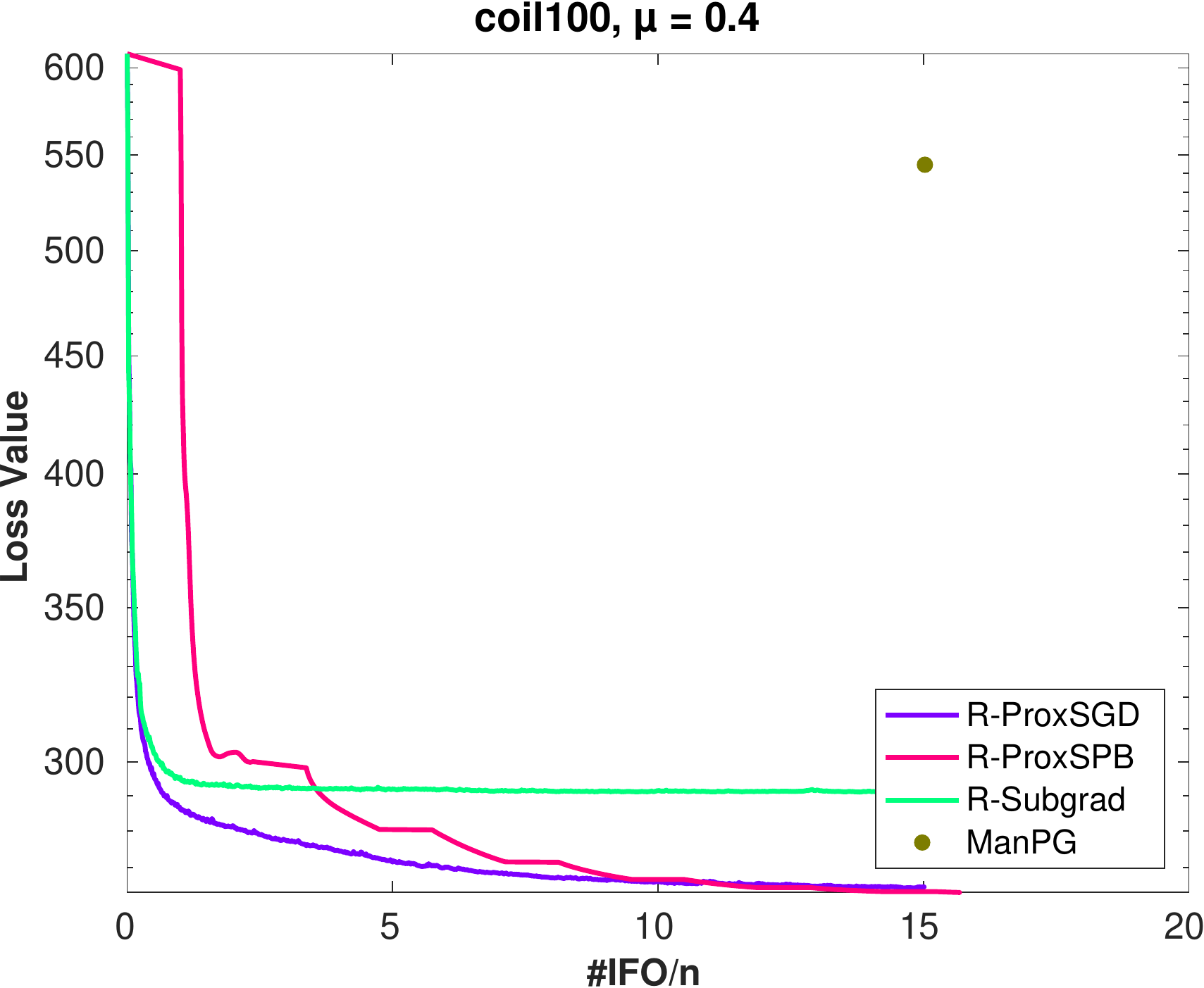}
\endminipage\hfill
\minipage{0.33\textwidth}
  \includegraphics[width=\linewidth]{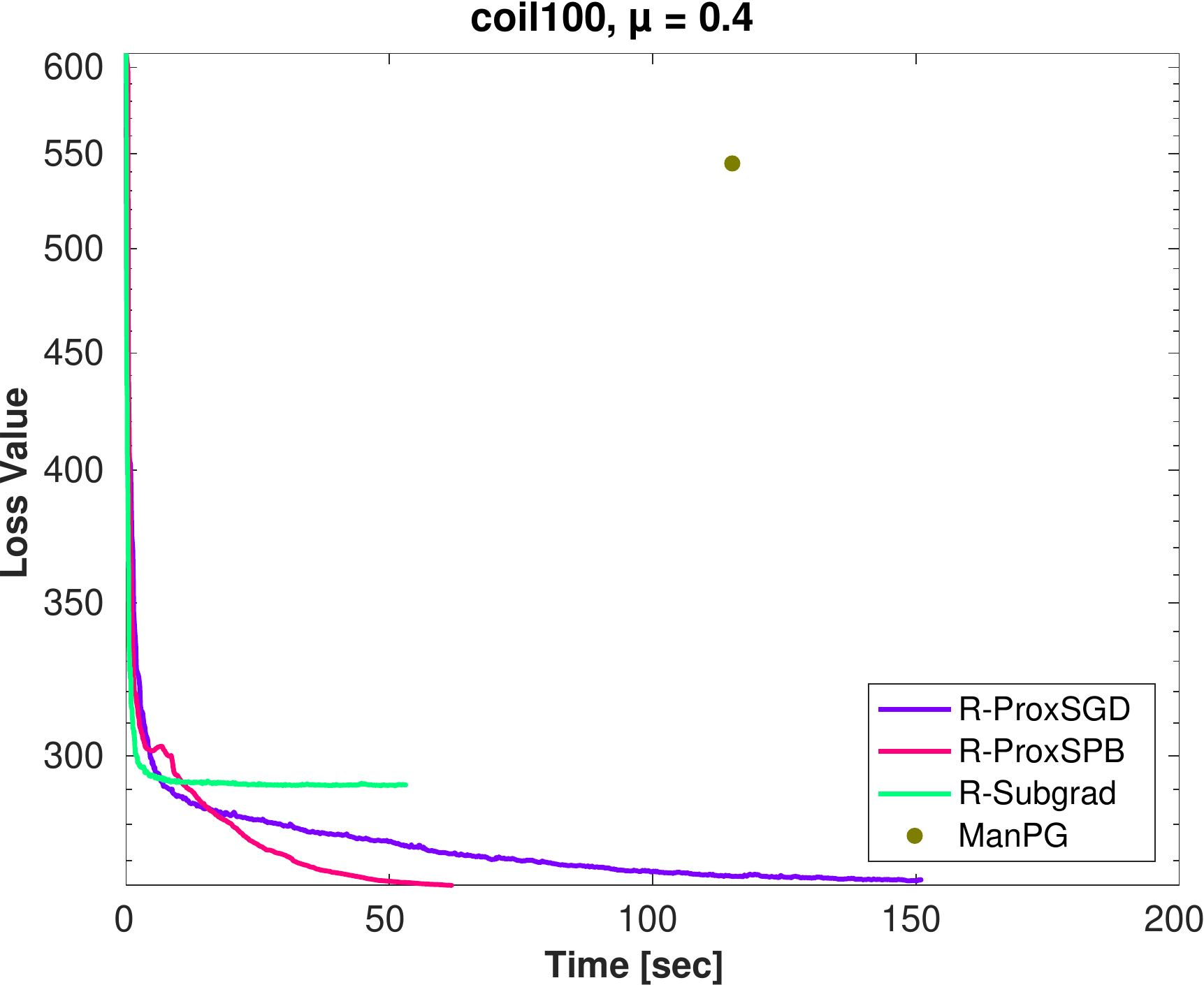}
\endminipage\hfill
\minipage{0.33\textwidth}%
  \includegraphics[width=\linewidth]{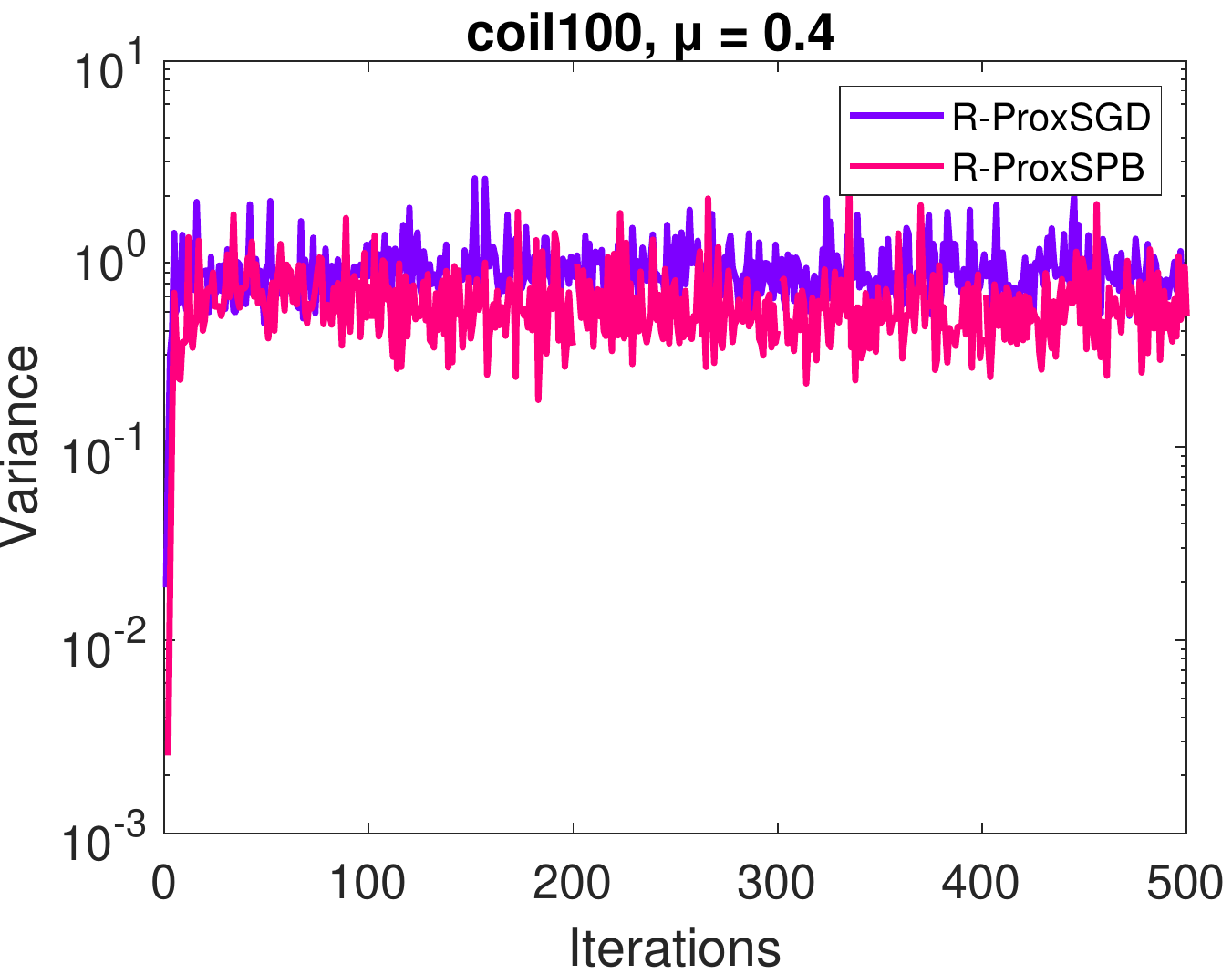}
\endminipage\hfill
\caption{Experimental results on the \texttt{coil100} dataset with $\mu = 0.2$ and $0.4$. }
\label{fig2}
\end{figure*}

All algorithms are implemented in Matlab and we use the \texttt{Manopt} \citep{manopt} package to compute vector transport, retraction and Riemannian gradient. {Since all of R-ProxSGD, R-ProxSPB, and R-Subgrad aim to solve the same problem \eqref{online-spca}, we evaluate the performance of those algorithms based on the objective function value $\mathbb{E}_{Z\in\mathcal{D}}[\|Z - XX^\top Z\|_2^2] + \mu\|X\|_1$ (``loss value'' in Figures \ref{fig1} and \ref{fig2}).} The experimental results are shown in Figures \ref{fig1} and \ref{fig2}. In particular, Figures \ref{fig1} and \ref{fig2} give results for $r=10$. More specifically, in Figure \ref{fig1} we report the results on the \texttt{mnist} dataset, and in Figure \ref{fig2} we report the results on the \texttt{coil100} dataset, both with two choices of $\mu$: $\mu=0.2$ and $\mu=0.4$. Note that $\mu$ is the parameter in \eqref{online-spca} controlling the sparsity of the solution. {In the first column of Figures \ref{fig1} and \ref{fig2}, we report the loss value in \eqref{online-spca} versus the number of IFO divided by $n$. In the second column of Figures \ref{fig1} and \ref{fig2}, we report the loss value versus the CPU time (in seconds).} In the third column of Figures \ref{fig1} and \ref{fig2}, we report the variance of gradient estimation versus the number of iterations, which is adopted in \cite{defazio2018ineffectiveness}. 

All the results in Figures \ref{fig1} and \ref{fig2} indicate that both R-ProxSGD and R-ProxSPB consistently outperform R-Subgrad in terms of CPU time and the number of IFO calls. Moreover, these figures show that R-Subgrad is not able to reduce the loss value to a desired accuracy, comparing with R-ProxSGD and R-ProxSPB. Furthermore, these results also show that R-ProxSPB usually performs better than R-ProxSGD, which is consistent with our theoretical results on the complexity bounds. Figures \ref{fig1} and \ref{fig2} also imply that R-ProxSPB is effective to reduce the variance of the stochastic gradient on both datasets. We perform grid search to tune $\eta_0$ (used in R-Subgrad and R-ProxSGD) and $\eta$ (used in R-ProxSPB) from $\{5\times 10^{-5}, 10^{-4}, 5\times 10^{-4}, ..., 1\}$. The best $\eta_0$ and $\eta$ on different data sets are reported in Table~\ref{best_eta}.

\begin{figure*}[htbp]
	\minipage{0.33\textwidth}
	\includegraphics[width=\linewidth]{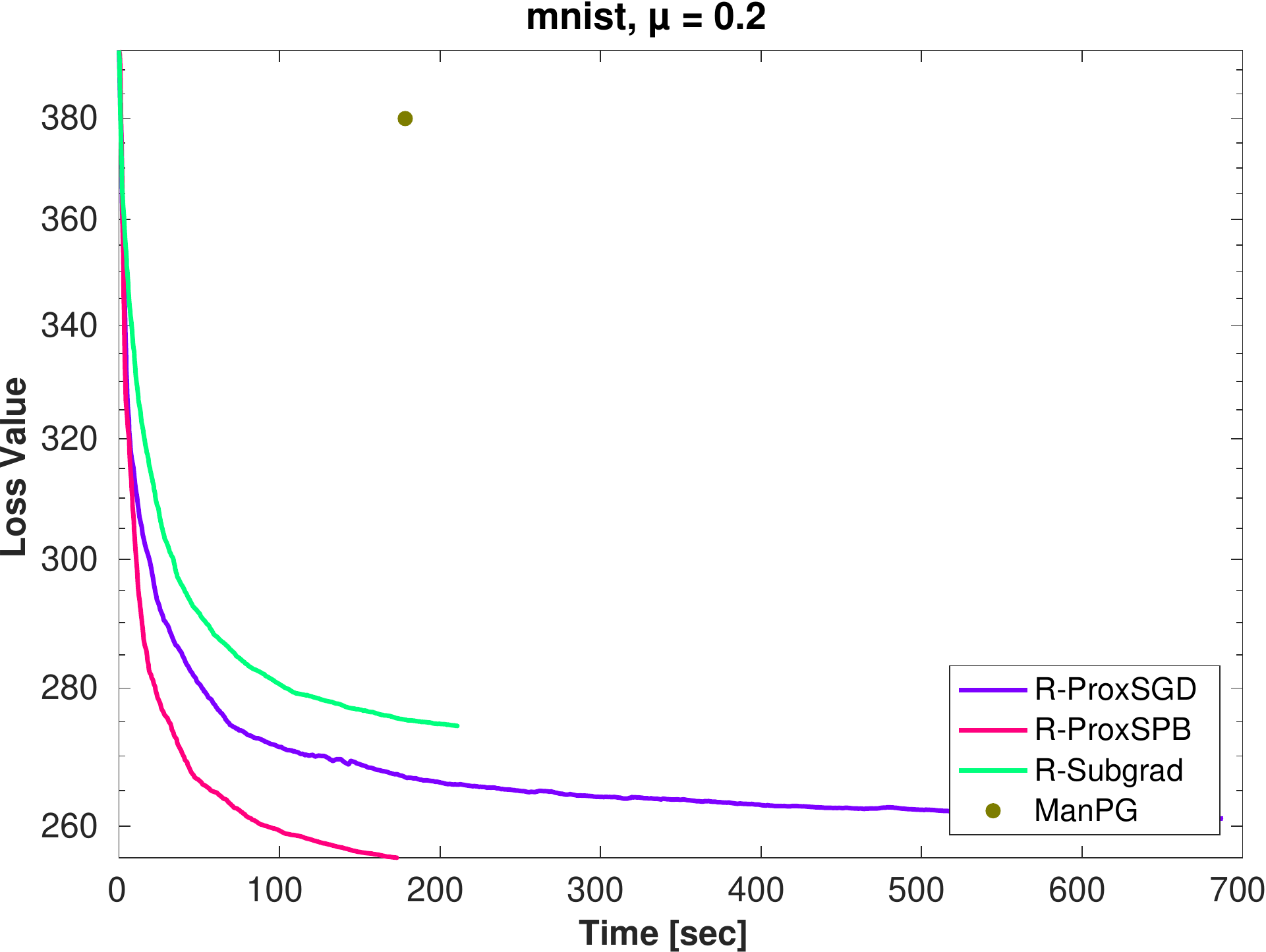}
	\endminipage\hfill
	\minipage{0.33\textwidth}
	\includegraphics[width=\linewidth]{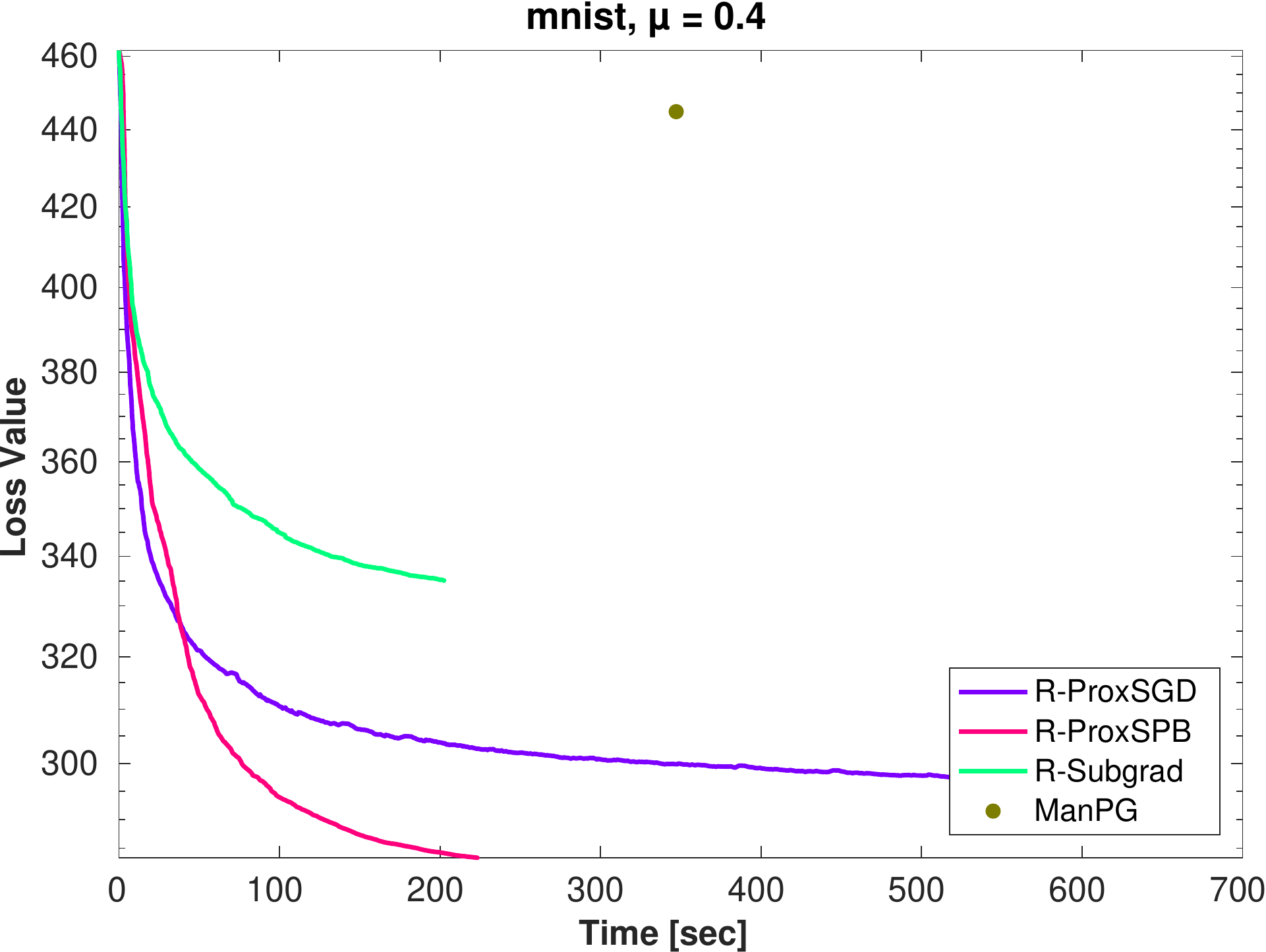}
	\endminipage\hfill
	\minipage{0.33\textwidth}%
	\includegraphics[width=\linewidth]{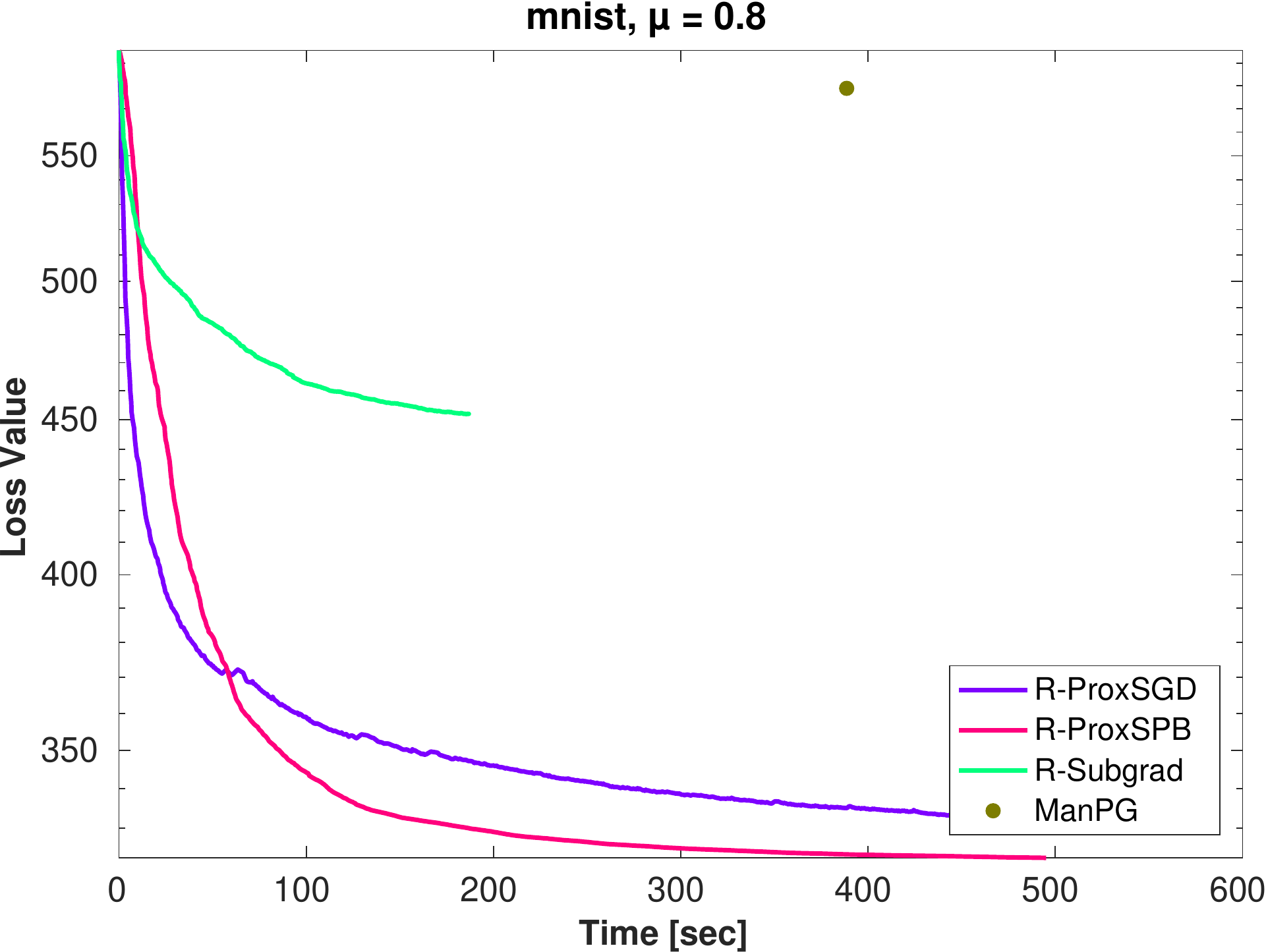}
	\endminipage\hfill
	\minipage{0.33\textwidth}
	\includegraphics[width=\linewidth]{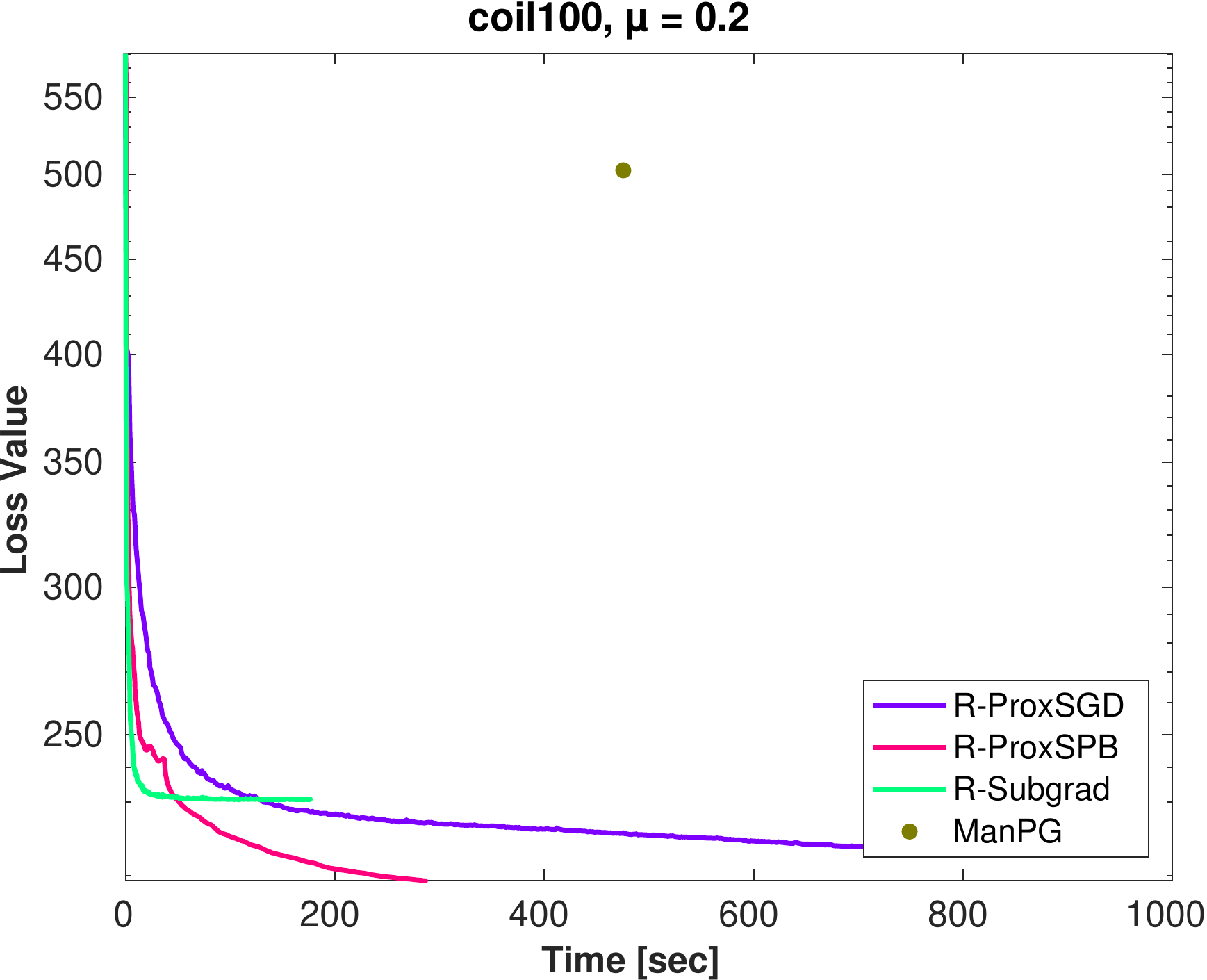}
	\endminipage\hfill
	\minipage{0.33\textwidth}
	\includegraphics[width=\linewidth]{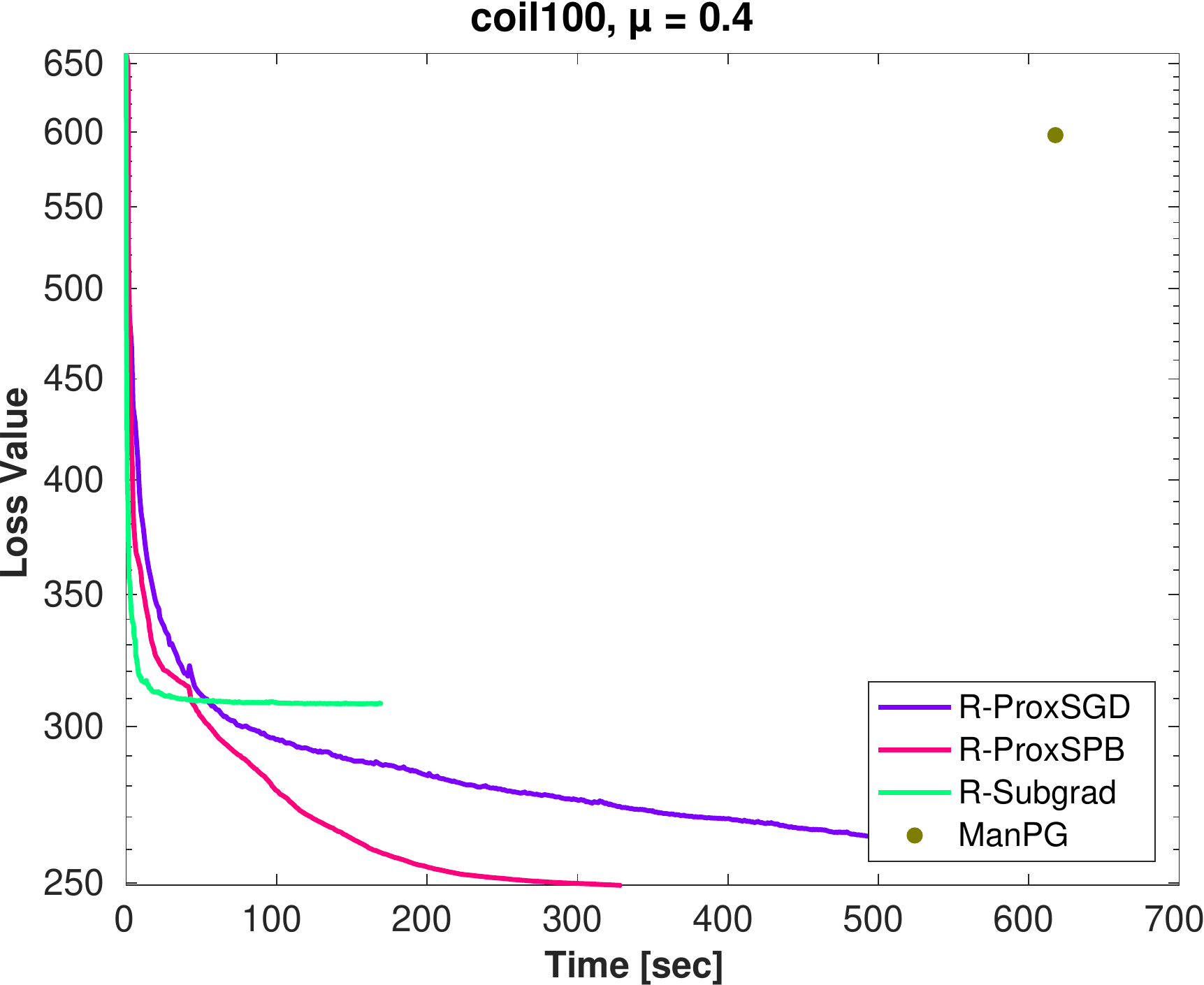}
	\endminipage\hfill
	\minipage{0.33\textwidth}%
	\includegraphics[width=\linewidth]{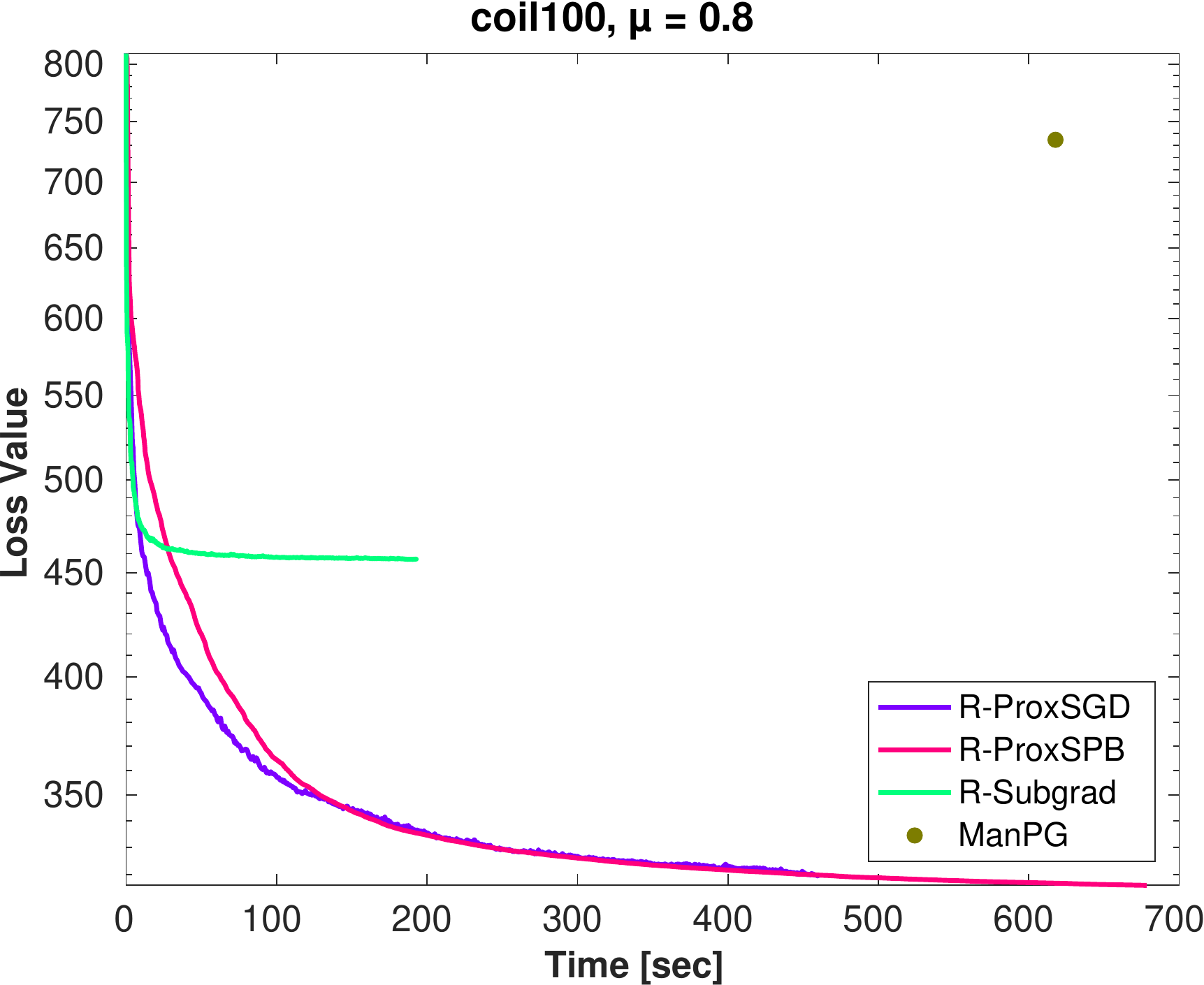}
	\endminipage\hfill
	\caption{Loss value versus runtime on two datasets with $r = 15$ and $\mu=0.2,0.4,0.8$.}
	\label{fig3}
\end{figure*}

Figure \ref{fig3} gives more results on the case $r = 15$ and $\mu = 0.2$, $0.4$, $0.8$, and here we only present the loss value versus the CPU time. These results further justify the advantages of our proposed R-ProxSGD and R-ProxSPB algorithms.

\begin{table}[H]
	\vspace{-.1cm}
	\center
		\resizebox{\columnwidth}{!}{
	\begin{tabular}{|lccc|lccc|}
		\hline
		\multicolumn{4}{|c|}{\texttt{mnist} data set}&	\multicolumn{4}{c|}{\texttt{coil} data set}\\
		\hline
		$\mu$ & R-Subgrad & R-ProxSGD & R-ProxSPB & $\mu$ & R-Subgrad & R-ProxSGD & R-ProxSPB\\
		\hline
		0.2 & 0.01 & 0.005 & 0.005 & 0.2 & 0.005 & 0.01 & 0.005\\
		0.4 & 0.01 & 0.01 & 0.005 & 	0.4 & 0.01 & 0.01 & 0.005 \\
		\hline
	\end{tabular}}
	\caption{Chosen $\eta_0$ (for R-Subgrad and R-ProxSGD) and $\eta$ (for R-ProxSPB) for the reported results on \texttt{mnist} and \texttt{coil} data sets.}
	\label{best_eta}\vspace{-.2cm}
\end{table}

{
\subsubsection{Comparison with ManPG and SPAMS}

To justify the necessity of introducing the stochasticity, we compare our R-ProxSGD and R-ProxSPB with the deterministic algorithm ManPG~\citep{chen2018proximal}, which also solves the same problem in \eqref{online-spca} but assumes that the full gradient information for the smooth part is available. As shown in Figures \ref{fig1} and \ref{fig2}, ManPG leads to larger loss values than our R-ProxSGD and R-ProxSPB given the same budget of gradient oracles or running time. We also point out that if the problem is online, then ManPG is not applicable.

Moreover, we also compare our R-ProxSPB algorithm with the SPAMS algorithm~\citep{mairal2010online} for the online sparse PCA problem. We run R-ProxSPB for 1000 iterations with batch size 100. For fair comparison, we run SPAMS using the same batch size and the same number of gradient oracles. Since the problem formulation of online sparse PCA in SPAMS is different from \eqref{online-spca}, we cannot directly compare the objective function value. Instead, we consider the explained variance and sparsity metrics as suggested in \cite{yang2015streaming}. The explained variance is defined as $\frac{\text{tr}(X^\top AA^\top X)}{\text{tr}(XX^\top)}$, where $A\in\mathbb{R}^{d\times n}$ is the data matrix and $X\in\mathbb{R}^{d\times r}$ is the model parameter. The sparsity is defined as the number of elements in $X$ whose absolute value is larger than the threshold 0.001.  As shown in Table~\ref{added_exp}, our R-ProxSPB leads to better sparsity while achieving comparable explained variance.

\begin{table}[htbp]
	\centering
			\resizebox{\columnwidth}{!}{
	\begin{tabular}{lcc|lcc} 
		\toprule 
		\multicolumn{3}{c|}{\texttt{coil100} Data Set} &   \multicolumn{3}{c}{\texttt{mnist} Data Set}\\
		\midrule
		Algorithms &  Explained Variance & Sparsity &Algorithms &  Explained Variance & Sparsity\\ 
		\midrule 
		SPAMS  & 0.0132 & 240 & SPAMS & 0.0179 & 185\\ 
		R-ProxSPB & 0.0120 & 22 &  R-ProxSPB  & 0.0190 & 18\\ 
		\bottomrule 
	\end{tabular} }
	\caption{Comparison of R-ProxSPB and SPAMS~\citep{mairal2010online}. } 
		\label{added_exp}
\end{table}

}

\subsection{Robust Low-Rank Matrix Completion}

Robust low-rank matrix completion is closely related to the robust PCA problem. The robust PCA aims to decompose a given matrix $M\in\mathbb{R}^{m\times n}$ into the superposition of a low-rank matrix $L$ and a sparse matrix $S$. Robust low-rank matrix completion is the same as robust PCA, except that only a subset of the entries of $M$ is observed. The convex formulations of them are studied extensively in the literature and we refer the reader to the recent survey \citep{Ma-Aybat-PIEEE}. A typical convex formulation of robust low-rank matrix completion is given as follows:
\begin{equation}\label{rmc-convex}
\min_{L,S} \ \|L\|_* + \gamma\|S\|_1, \ \st, \ \Pc_\Omega(L+S) = \Pc_\Omega (M),
\end{equation}
where $\|L\|_*$ denotes the nuclear norm of $L$ and it sums the singular values of $L$, $\Omega$ is a subset of the index set $\{(i,j)\mid 1\leq i\leq m, 1\leq j\leq n\}$, and the projection operator $\Pc_\Omega$ is defined as: $[\Pc_\Omega(Z)]_{ij}=Z_{ij}$, if $(i,j)\in\Omega$, and $[\Pc_\Omega(Z)]_{ij}=0$ otherwise. Due to the presence of the nuclear norm in \eqref{rmc-convex}, algorithms for solving \eqref{rmc-convex} usually require computing the SVD of an $m\times n$ matrix in every iteration, which can be time consuming when $m$ and $n$ are large. Recently, some nonconvex formulations of robust low-rank matrix completion were proposed because they allow more efficient and scalable algorithms. In \cite{Huang-ManPGMC-2020}, the authors proposed the following nonconvex formulation of robust low-rank matrix completion:
\begin{eqnarray}\label{rmc-huang}
\min_{\mathbb{U}\in\Gr(m,r),V\in\mathbb{R}^{r\times n},S\in\mathbb{R}^{m\times n}} \frac{1}{2}\|\mathcal{P}_\Omega (UV - M + S)\|^2_F + \frac{\lambda}{2}\|\mathcal{P}_{\bar{\Omega}} (UV)\|^2_F + \gamma \|\mathcal{P}_\Omega (S)\|_1,
\end{eqnarray}
where $\Gr(m,r)$ denotes the Grassmann manifold, which is the set of $r$-dimensional vector subspaces of $\mathbb{R}^m$, and we use $U\in\mathbb{R}^{m\times r}$ to denote a basis of the subspace $\mathbb{U}\in\Gr(m,r)$. 
In \eqref{rmc-huang}, the low-rank matrix $L$ is replaced by $UV$ with $U\in\mathbb{R}^{m\times r}$, $V\in\mathbb{R}^{r\times n}$, and $r$ is the estimation of the rank of $L$; the term $\frac{\lambda}{2}\|\mathcal{P}_{\bar{\Omega}} (UV)\|^2_F$ is added as a regularizer and $\lambda>0$ is sufficiently small indicating that we have a small confidence of the components of $UV$ on $\bar{\Omega}$ being zeros; the constraint $\mathbb{U}\in\Gr(m,r)$ is added to remove the scaling ambiguity of $U$ and $V$. The nonconvex formulation \eqref{rmc-huang} was motivated by some recent works on Riemannian optimization \citep{boumal2011rtrmc,Cambier-Absil-2016}. Note that, for fixed $U$ and $S$, the variable $V$ in \eqref{rmc-huang} can be uniquely determined. By denoting
\be\label{def-f-bar}
\bar{f}(U,V,S) = \frac{1}{2}\|\mathcal{P}_{\Omega}(UV - M + S)\|_F^2 + \frac{\lambda}{2}\|\mathcal{P}_{\bar{\Omega}}(UV)\|_F^2,
\ee
and
\be\label{def-V}
V_{U,S}:=\argmin_V \ \bar{f}(U,V,S), \mbox{ and } f(U,S) = \bar{f}(U,V_{U,S},S),
\ee
we can rewrite \eqref{rmc-huang} as
\begin{eqnarray}\label{rmc-huang-1}
\min_{\mathbb{U}\in\Gr(m,r), S\in\mathbb{R}^{m\times n}} f(U,S) + \gamma \|\mathcal{P}_\Omega (S)\|_1,
\end{eqnarray}
which is a Riemannian optimization problem with nonsmooth objective. Note that although the manifold is the Grassmann manifold instead of the Stiefel manifold, our algorithms discussed in Section \ref{approach} can be directly applied to \eqref{rmc-huang-1}. To see this, first note that as suggested in \citep{boumal2011rtrmc}, without loss of generality, we can restrict matrix $U$ as an orthonormal basis of $\mathbb{U}$. Therefore, we have
\[\|\mathcal{P}_{\bar{\Omega}}(UV)\|_F^2=\|UV\|_F^2-\|\mathcal{P}_\Omega\|_F^2=\|V\|_F^2-\|\mathcal{P}_\Omega\|_F^2,\]
and thus we can rewrite $\bar{f}(U,V,S)$ and $f(U,S)$ as
\begin{equation}\label{def-bar-f-new}
\bar{f}(U,V,S) =  \frac{1}{2}\|\mathcal{P}_\Omega(UV - M + S)\|_F^2 + \frac{\lambda}{2} \|V\|^2_F - \frac{\lambda}{2}\|\mathcal{P}_\Omega(UV)\|_F^2.
\end{equation}
\begin{equation}\label{def-f-U-S-new}
f(U, S) =  \frac{1}{2}\|\mathcal{P}_\Omega(UV_{U,S} - M + S)\|_F^2 + \frac{\lambda}{2} \|V_{U,S}\|^2_F - \frac{\lambda}{2}\|\mathcal{P}_\Omega(UV_{U,S})\|_F^2.
\end{equation}
From \eqref{def-V} we know that $\nabla_V \bar{f}(U,V_{U,S},S)=0$. Therefore,
\[\nabla_U f(U,S) = \nabla_U \bar{f}(U,V_{U,S},S) =\nabla_1 \hat{f}(U,V_{U,S},S),\]
where
\begin{equation}\label{def-hat-f}
\hat{f}(U,V_{U,S},S):=\frac{1}{2}\|\mathcal{P}_\Omega(UV_{U,S} - M + S)\|_F^2 - \frac{\lambda}{2}\|\mathcal{P}_\Omega(UV_{U,S})\|_F^2 = \sum_{(i,j)\in\Omega}\hat{f}_{ij}(U,V_{U,S},S),
\end{equation}
and
\[\hat{f}_{ij}(U,V_{U,S},S) = \frac{1}{2}(UV_{U,S} - M + S)_{ij}^2-\frac{\lambda}{2}(UV_{U,S})_{ij}^2.\]
That is, $\hat{f}$ in \eqref{def-hat-f} has a natural finite-sum structure, and a stochastic gradient approximation to $\nabla_U f(U,S)$ is given by $\nabla_1\hat{f}_{ij}(U,V_{U,S},S)$ with randomly sampled index pair $(i,j)\in\Omega$. It is easy to verify that  
\[\nabla_1\hat{f}_{ij}(U,V_{U,S},S)=(u_i^\top v_j-M_{ij}+S_{ij}-\lambda u_i^\top v_j) \bar{V}_j^\top,\]
where $u_i^\top$ denotes the $i$-th row of $U$, and $v_j$ denotes the $j$-th column of $V_{U,S}$, and
\[\bar{V}_j = \begin{bmatrix} 0 & 0 & \cdots & v_j & \cdots & 0 \end{bmatrix}.\]
That is, $\bar{V}_j \in\mathbb{R}^{r\times m}$ is a matrix whose $j$-th column is $v_j$ and all other columns are zeros. Clearly, when computing $\nabla_U f_{ij}(U,S)$, we only need to access $u_i^\top$ and $v_j$ and we do not need to access the whole matrix $U$ and $V_{U,S}$ and compute the matrix multiplication $UV_{U,S}$, and this is very useful when $m$ and $n$ are large.

\begin{figure*}[t]%
	\centering
	\subfloat[Original Image]{\includegraphics[width=.249\linewidth]{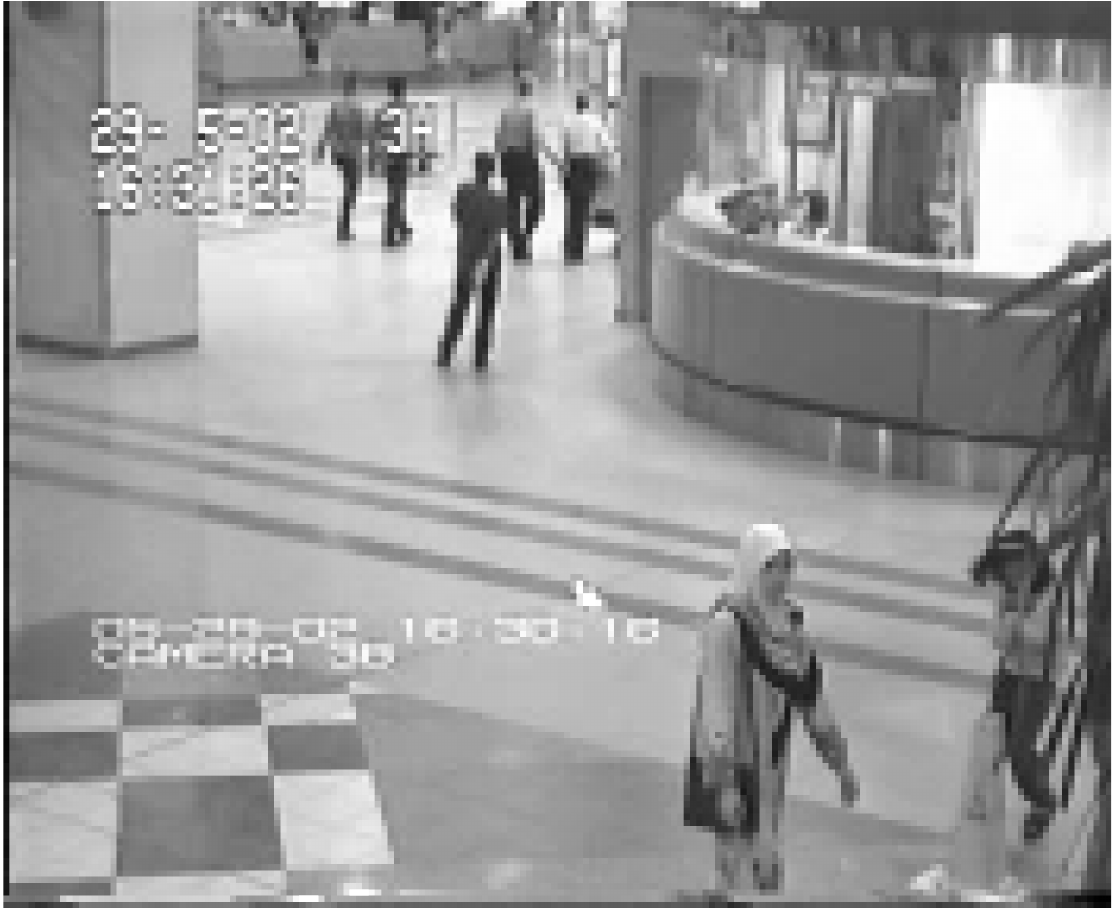}}
	\subfloat[R-ProxSPB\\time=2.58]{\includegraphics[width=.249\linewidth]{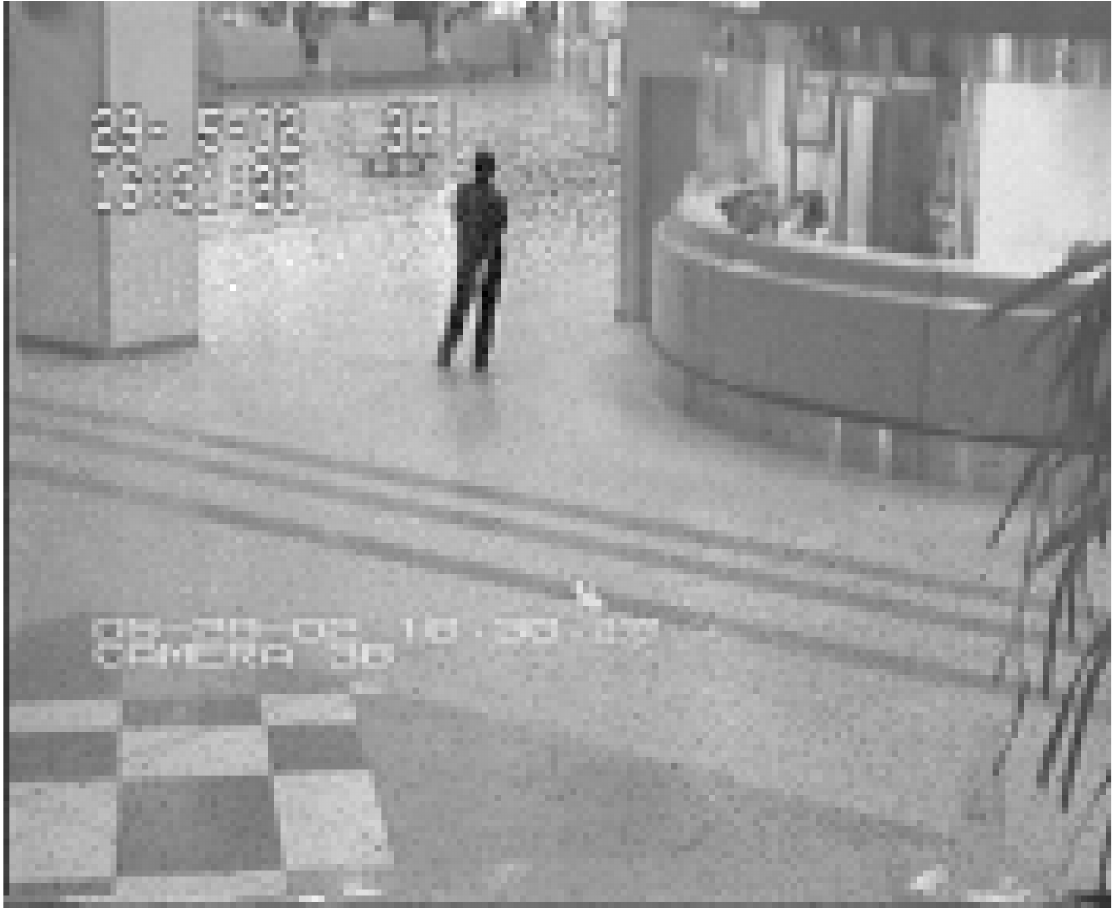}}
	\subfloat[R-ProxSGD\\ time=5.76]{\includegraphics[width=.249\linewidth]{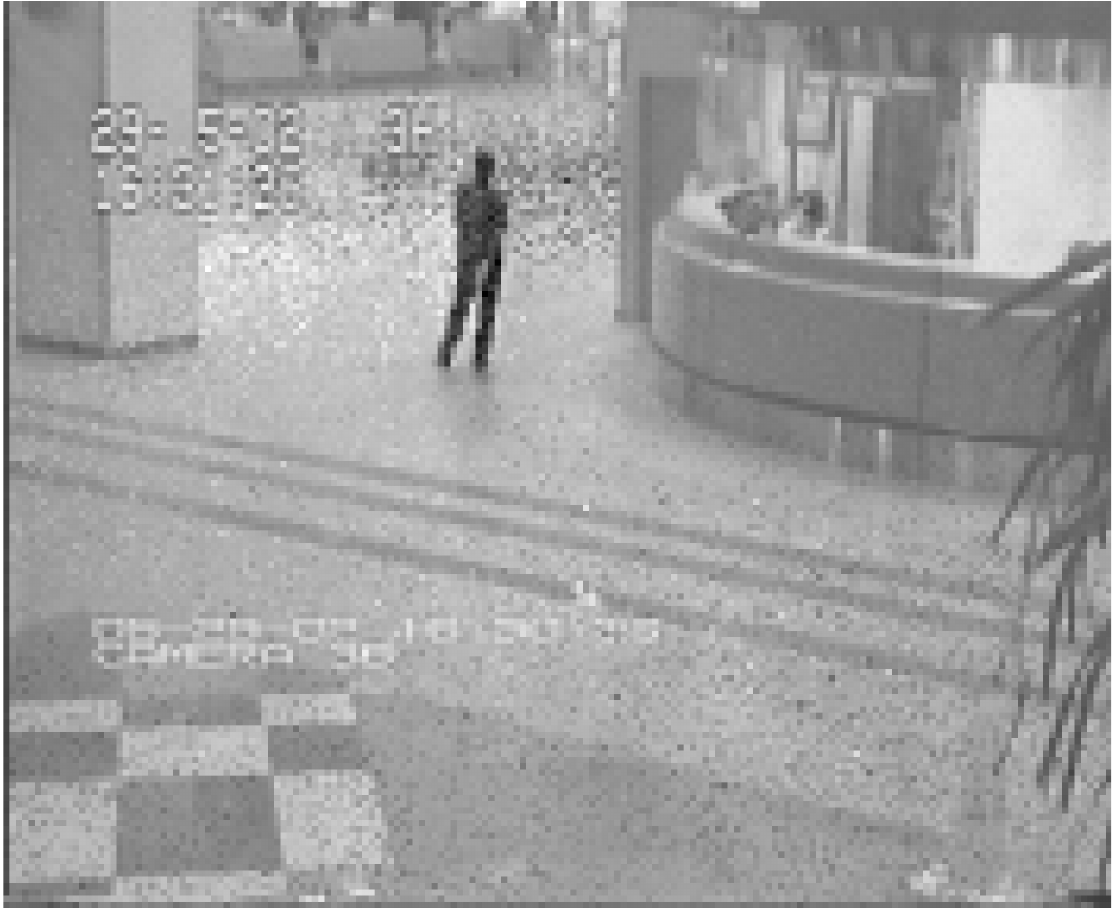}} 
	\subfloat[R-Subgrad\\time=2.26]{\includegraphics[width=.249\linewidth]{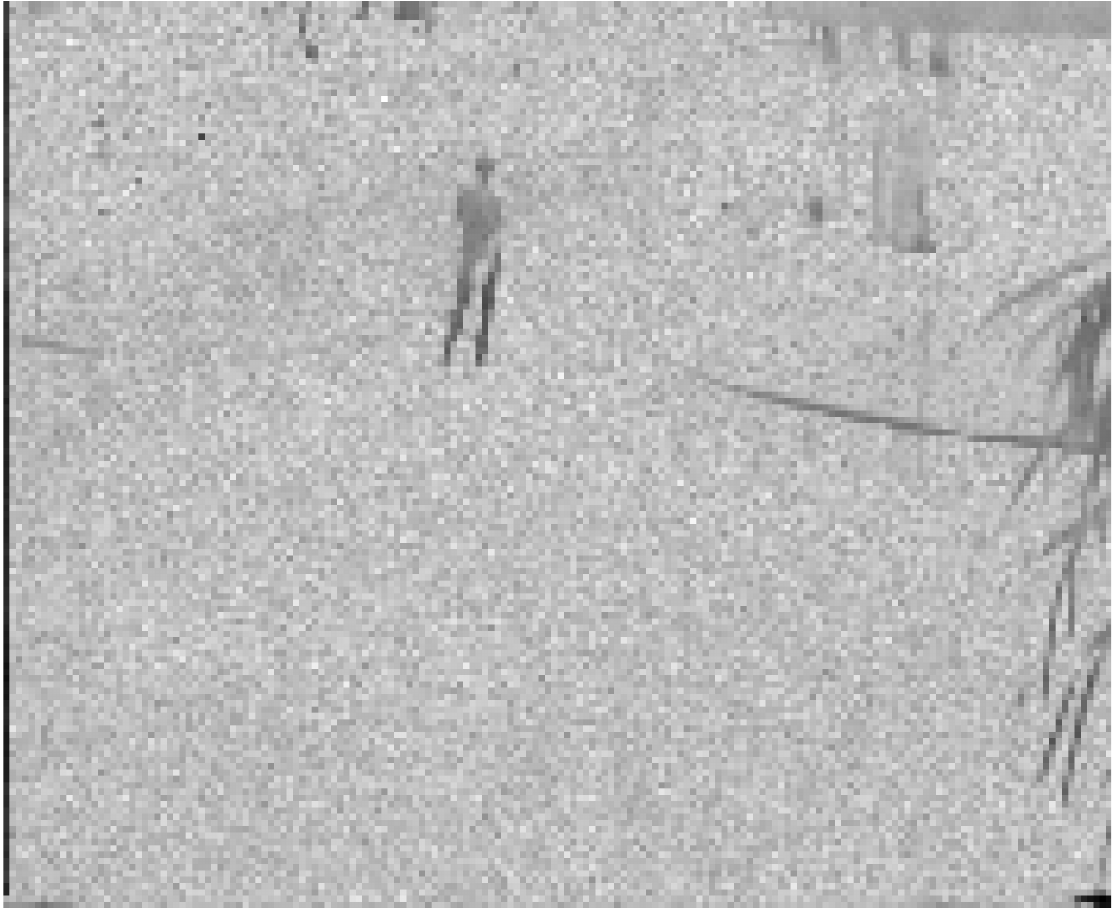}} \\
	\subfloat[Original Image]{\includegraphics[width=.249\linewidth]{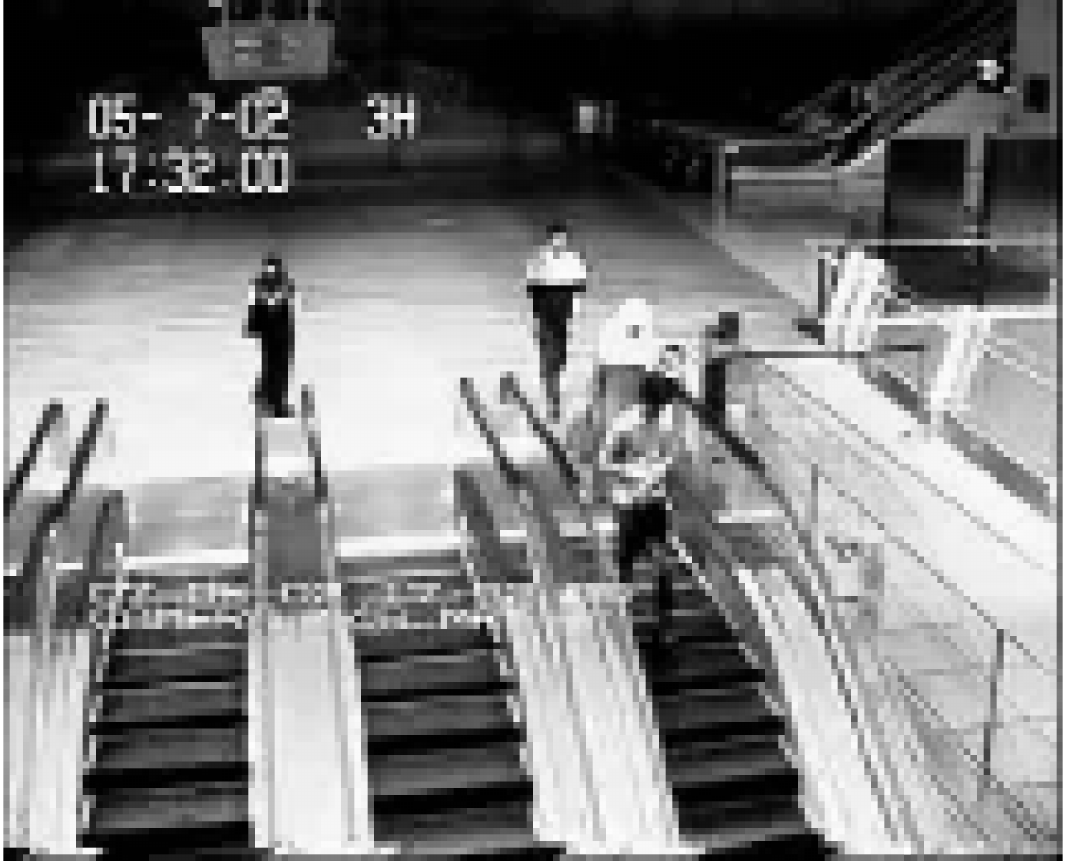}}
	\subfloat[R-ProxSPB\\ time=3.11]{\includegraphics[width=.249\linewidth]{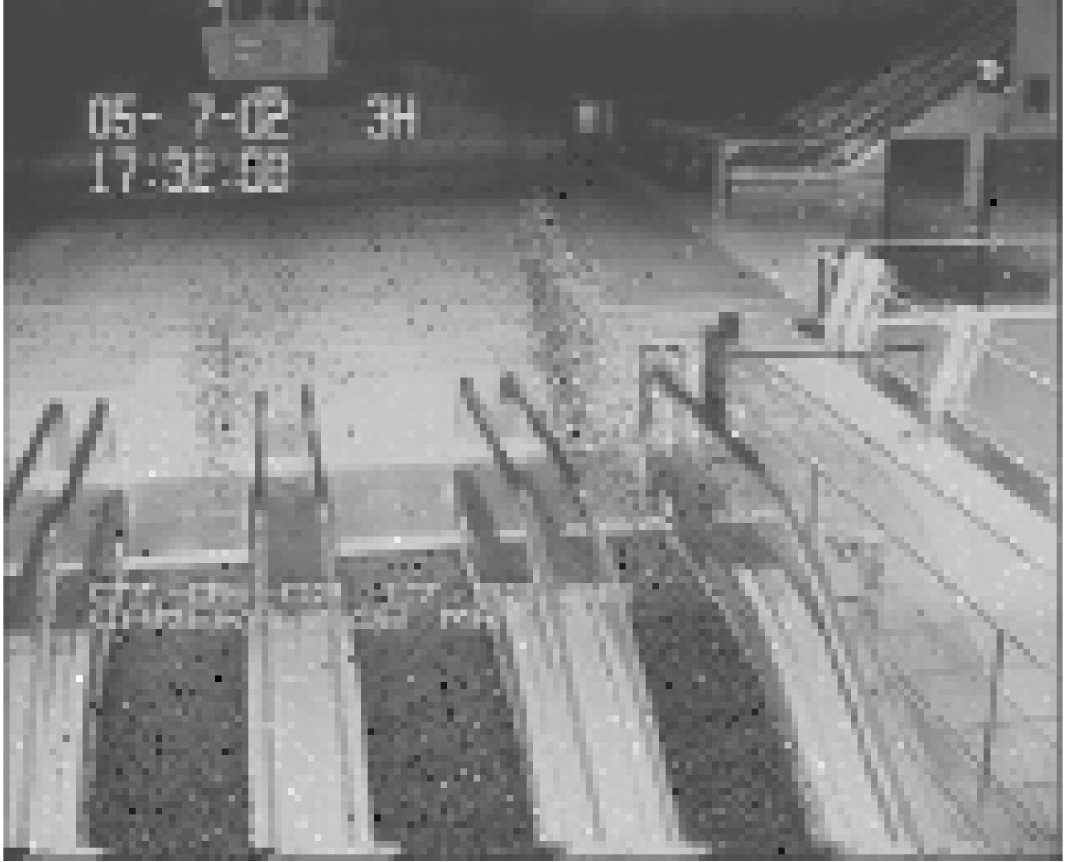}}
	\subfloat[R-ProxSGD\\  time=6.99]{\includegraphics[width=.249\linewidth]{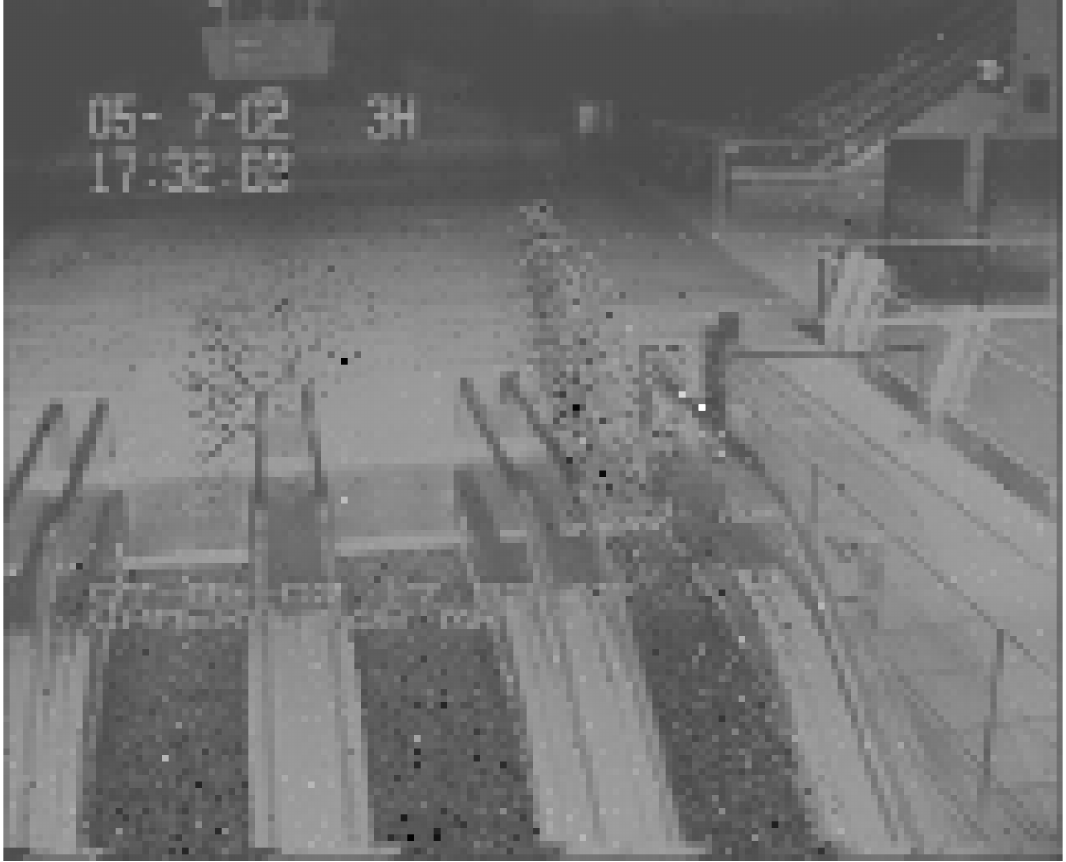}} 
	\subfloat[R-Subgrad\\ time=2.45]{\includegraphics[width=.249\linewidth]{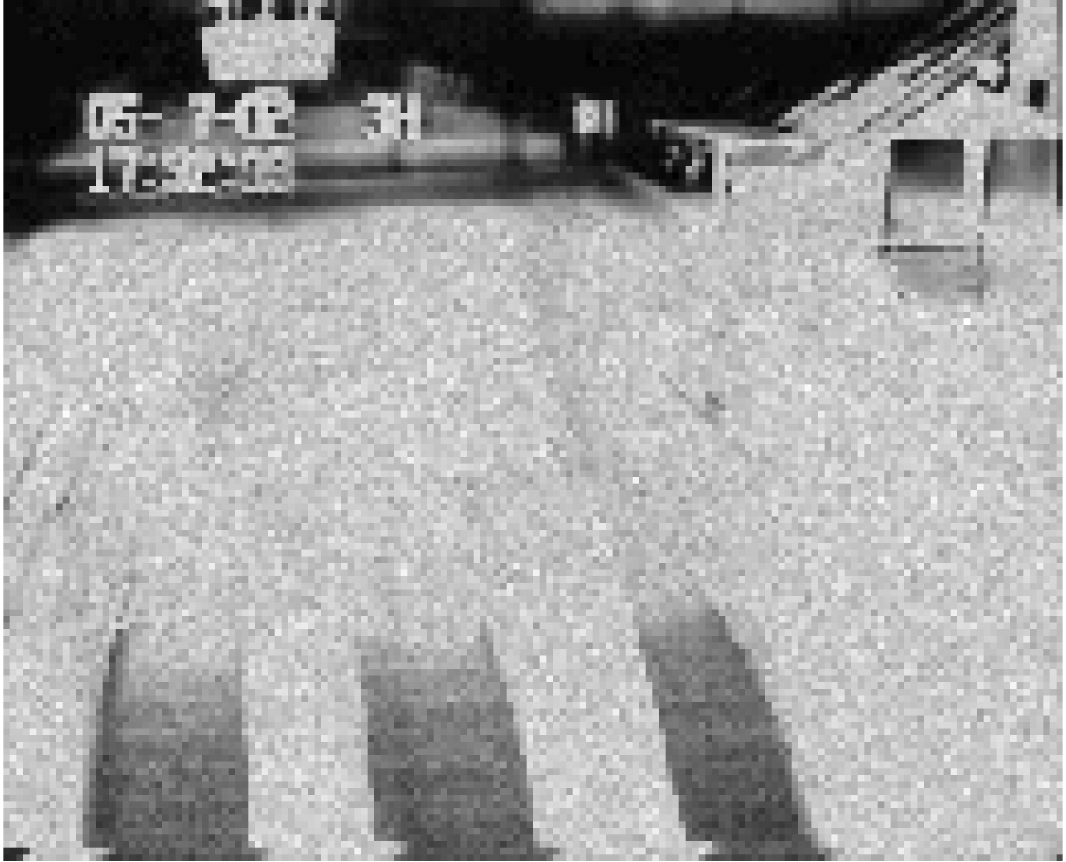}} \\
	\caption{First row: background estimation from partial observations on the ``Hall of a business building'' data set; Second row:  background estimation from partial observations on the ``Airport elevator'' data set.}
\label{rmc_fig}
\end{figure*}

{
We applied our R-ProxSGD and R-ProxSPB algorithms to solve the robust low-rank matrix completion problem \eqref{rmc-huang-1} on some real data for video background estimation \citep{li2004statistical} and we again compared their performance with R-Subgrad. We consider two surveillance video datasets: ``Hall of a business building'' and ``Airport elevator''. The data matrix $X^*$ is obtained by vectorizing each grayscale frame of the video. 
We then randomly sample 50\% of the indices to obtain $\Omega$, and then sample the entries of $X^*$ from $\Omega$ to get $M$. A sparse matrix $S^*$ was then added to $M$. In R-ProxSGD and R-Subgrad, we randomly sample 10\% of the known entries as a batch in each iteration. In R-ProxSPB, we set $|\mathcal{S}_t^1| = |\Omega|$ and $q=5$. The initial step sizes $\eta_0$ are tuned from $\{10^{-j}/|\Omega|,~i = 0, 1,\ldots, 4\}$.

In Figure \ref{rmc_fig}, we present the experimental results on the problem with those two real datasets. For fair comparison, we report the results of all algorithms using the same budget of stochastic gradients, which is $4|\Omega|$. The results in Figure \ref{rmc_fig} clearly show the advantage of our R-ProxSPB and R-ProxSGD algorithms over R-Subgrad algorithm. 
}

\section{Conclusion}\label{conclu}

In this paper, we considered the nonsmooth Riemannian optimization problems with nonsmooth regularizer in the objective. We designed Riemannian stochastic algorithms that do not need subgradiet information for solving this class of problems. Specifically, we proposed two Riemannian stochastic proximal gradient algorithms: R-ProxSGD and R-ProxSPB to solve this problem. The two proposed algorithms are generalizations of their counterparts in Euclidean space to Riemannian manifold setting. 
We analyzed the iteration complexity and IFO complexity of the proposed algorithms for obtaining an $\epsilon$-stationary point. Numerical results on solving online sparse PCA and robust low-rank matrix completion are conducted which demonstrate that our proposed algorithms outperform significantly the Riemannian stochastic subgradient method. Future work includes extending the current results to more general Riemannian manifolds.

\section*{Acknowledgement}
The authors would like to thank Shixiang Chen for fruitful discussions. The authors are very grateful for the associate editor and the two reviewers for very constructive comments and suggestions that led to significant improvement of the presentation of this paper.

\bibliography{reference}

\begin{thebibliography}{50}
\providecommand{\natexlab}[1]{#1}
\providecommand{\url}[1]{\texttt{#1}}
\expandafter\ifx\csname urlstyle\endcsname\relax
  \providecommand{\doi}[1]{doi: #1}\else
  \providecommand{\doi}{doi: \begingroup \urlstyle{rm}\Url}\fi

\bibitem[Absil et~al.(2009)Absil, Mahony, and Sepulchre]{absil2009optimization}
P-A Absil, Robert Mahony, and Rodolphe Sepulchre.
\newblock \emph{Optimization algorithms on matrix manifolds}.
\newblock Princeton University Press, 2009.

\bibitem[Baden et~al.(2016)Baden, Berens, Franke, Ros{\'o}n, Bethge, and
  Euler]{baden2016functional}
Tom Baden, Philipp Berens, Katrin Franke, Miroslav~Rom{\'a}n Ros{\'o}n,
  Matthias Bethge, and Thomas Euler.
\newblock The functional diversity of retinal ganglion cells in the mouse.
\newblock \emph{Nature}, 529\penalty0 (7586):\penalty0 345--350, 2016.

\bibitem[Bonnabel(2013)]{bonnabel2013stochastic}
Silvere Bonnabel.
\newblock Stochastic gradient descent on {R}iemannian manifolds.
\newblock \emph{IEEE Transactions on Automatic Control}, 58\penalty0
  (9):\penalty0 2217--2229, 2013.

\bibitem[Boumal et~al.(2014)Boumal, Mishra, Absil, and Sepulchre]{manopt}
N.~Boumal, B.~Mishra, P.-A. Absil, and R.~Sepulchre.
\newblock {M}anopt, a {M}atlab toolbox for optimization on manifolds.
\newblock \emph{Journal of Machine Learning Research}, 15:\penalty0 1455--1459,
  2014.
\newblock URL \url{http://www.manopt.org}.

\bibitem[Boumal and Absil(2011)]{boumal2011rtrmc}
Nicolas Boumal and Pierre-Antoine Absil.
\newblock {RTRMC}: A {R}iemannian trust-region method for low-rank matrix
  completion.
\newblock In \emph{Advances in neural information processing systems}, pages
  406--414, 2011.

\bibitem[Boumal et~al.(2019)Boumal, Absil, and Cartis]{boumal2018global}
Nicolas Boumal, Pierre-Antoine Absil, and Coralia Cartis.
\newblock Global rates of convergence for nonconvex optimization on manifolds.
\newblock \emph{IMA Journal of Numerical Analysis}, 39\penalty0 (1):\penalty0
  1--33, 2019.

\bibitem[Cadima and Jolliffe(1995)]{cadima1995loading}
Jorge Cadima and Ian~T Jolliffe.
\newblock Loading and correlations in the interpretation of principal
  compenents.
\newblock \emph{Journal of applied Statistics}, 22\penalty0 (2):\penalty0
  203--214, 1995.

\bibitem[Cambier and Absil(2016)]{Cambier-Absil-2016}
L.~Cambier and P.-A. Absil.
\newblock Robust low-rank matrix completion by {R}iemannian optimization.
\newblock \emph{SIAM J. Sci. Comput.}, 38\penalty0 (5):\penalty0 S440--S460,
  2016.

\bibitem[Chen et~al.(2020{\natexlab{a}})Chen, Deng, Ma, and
  So]{chen2020manifold}
Shixiang Chen, Zengde Deng, Shiqian Ma, and Anthony Man-Cho So.
\newblock Manifold proximal point algorithms for dual principal component
  pursuit and orthogonal dictionary learning.
\newblock \emph{arXiv preprint https://arxiv.org/abs/2005.02356},
  2020{\natexlab{a}}.

\bibitem[Chen et~al.(2020{\natexlab{b}})Chen, Ma, So, and
  Zhang]{chen2018proximal}
Shixiang Chen, Shiqian Ma, Anthony Man-Cho So, and Tong Zhang.
\newblock Proximal gradient method for nonsmooth optimization over the
  {S}tiefel manifold.
\newblock \emph{SIAM J. Optimization}, 30\penalty0 (1):\penalty0 210--239, Jan
  2020{\natexlab{b}}.

\bibitem[Chen et~al.(2020{\natexlab{c}})Chen, Ma, Xue, and Zou]{A-ManPG-2019}
Shixiang Chen, Shiqian Ma, Lingzhou Xue, and Hui Zou.
\newblock An alternating manifold proximal gradient method for sparse principal
  component analysis and sparse canonical correlation analysis.
\newblock \emph{INFORMS Journal on Optimization}, 2\penalty0 (3):\penalty0
  192--208, 2020{\natexlab{c}}.

\bibitem[Defazio et~al.(2014)Defazio, Bach, and Lacoste-Julien]{SAGA-2014}
A.~Defazio, F.~Bach, and S.~Lacoste-Julien.
\newblock {SAGA}: A fast incremental gradient method with support for
  non-strongly convex composite objectives.
\newblock In \emph{NIPS}, 2014.

\bibitem[Defazio and Bottou(2019)]{defazio2018ineffectiveness}
Aaron Defazio and L{\'e}on Bottou.
\newblock On the ineffectiveness of variance reduced optimization for deep
  learning.
\newblock \emph{NeurIPS}, 2019.

\bibitem[Edelman et~al.(1999)Edelman, Arias, and Smith]{EdelmanAriasSmith1999}
A.~Edelman, T.~A. Arias, and S.~T. Smith.
\newblock The geometry of algorithms with orthogonality constraints.
\newblock \emph{SIAM J. Matrix Anal. Appl.}, 20\penalty0 (2):\penalty0 303--353
  (electronic), 1999.
\newblock ISSN 0895-4798.

\bibitem[Fang et~al.(2018)Fang, Li, Lin, and Zhang]{fang2018spider}
Cong Fang, Chris~Junchi Li, Zhouchen Lin, and Tong Zhang.
\newblock {SPIDER}: Near-optimal non-convex optimization via stochastic
  path-integrated differential estimator.
\newblock In \emph{Advances in Neural Information Processing Systems}, pages
  689--699, 2018.

\bibitem[Gravuer et~al.(2008)Gravuer, Sullivan, Williams, and
  Duncan]{gravuer2008strong}
Kelly Gravuer, Jon~J Sullivan, Peter~A Williams, and Richard~P Duncan.
\newblock Strong human association with plant invasion success for trifolium
  introductions to new zealand.
\newblock \emph{Proceedings of the National Academy of Sciences}, 105\penalty0
  (17):\penalty0 6344--6349, 2008.

\bibitem[Hardoon and Shawe-Taylor(2011)]{hardoon2011sparse}
David~R Hardoon and John Shawe-Taylor.
\newblock Sparse canonical correlation analysis.
\newblock \emph{Machine Learning}, 83\penalty0 (3):\penalty0 331--353, 2011.

\bibitem[Hosseini and Pouryayevali(2011)]{hosseini2011generalized}
S~Hosseini and MR~Pouryayevali.
\newblock Generalized gradients and characterization of epi-{L}ipschitz sets in
  {R}iemannian manifolds.
\newblock \emph{Nonlinear Analysis: Theory, Methods \& Applications},
  74\penalty0 (12):\penalty0 3884--3895, 2011.

\bibitem[Huang et~al.(2020)Huang, Ma, and Lai]{Huang-ManPGMC-2020}
M.~Huang, S.~Ma, and L.~Lai.
\newblock Robust low-rank matrix completion via an alternating manifold
  proximal gradient continuation method.
\newblock \emph{https://arxiv.org/abs/2008.07740}, 2020.

\bibitem[Huang and Wei(2019)]{Huang-Wei-2019}
Wen Huang and Ke~Wei.
\newblock Riemannian proximal gradient methods.
\newblock \emph{arXiv preprint arXiv:1909.06065}, 2019.

\bibitem[Johnson and Zhang(2013)]{johnson2013accelerating}
Rie Johnson and Tong Zhang.
\newblock Accelerating stochastic gradient descent using predictive variance
  reduction.
\newblock In \emph{Advances in neural information processing systems}, pages
  315--323, 2013.

\bibitem[Jolliffe et~al.(2003)Jolliffe, N.Trendafilov, and Uddin]{Jolliffe2003}
I.~Jolliffe, N.Trendafilov, and M.~Uddin.
\newblock A modified principal component technique based on the lasso.
\newblock \emph{Journal of computational and Graphical Statistics}, 12\penalty0
  (3):\penalty0 531--547, 2003.

\bibitem[Kasai et~al.(2018)Kasai, Sato, and Mishra]{kasai2018riemannian}
Hiroyuki Kasai, Hiroyuki Sato, and Bamdev Mishra.
\newblock Riemannian stochastic recursive gradient algorithm with retraction
  and vector transport and its convergence analysis.
\newblock In \emph{International Conference on Machine Learning}, pages
  2521--2529, 2018.

\bibitem[LeCun(1998)]{lecun1998mnist}
Yann LeCun.
\newblock The mnist database of handwritten digits.
\newblock \emph{http://yann. lecun. com/exdb/mnist/}, 1998.

\bibitem[Li et~al.(2004)Li, Huang, Gu, and Tian]{li2004statistical}
Liyuan Li, Weimin Huang, Irene Yu-Hua Gu, and Qi~Tian.
\newblock Statistical modeling of complex backgrounds for foreground object
  detection.
\newblock \emph{IEEE Transactions on Image Processing}, 13\penalty0
  (11):\penalty0 1459--1472, 2004.

\bibitem[Li et~al.(2019)Li, Chen, Deng, Qu, Zhu, and So]{li2019nonsmooth}
Xiao Li, Shixiang Chen, Zengde Deng, Qing Qu, Zhihui Zhu, and Anthony Man~Cho
  So.
\newblock Weakly convex optimization over {S}tiefel manifold using {R}iemannian
  subgradient-type methods.
\newblock \emph{arXiv preprint arXiv:1911.05047}, 2019.

\bibitem[Ma and Aybat(2018)]{Ma-Aybat-PIEEE}
S.~Ma and N.~S. Aybat.
\newblock Efficient optimization algorithms for robust principal component
  analysis and its variants.
\newblock \emph{Proceedings of the IEEE}, 106\penalty0 (8):\penalty0
  1411--1426, 2018.

\bibitem[Mairal et~al.(2010)Mairal, Bach, Ponce, and Sapiro]{mairal2010online}
Julien Mairal, Francis Bach, Jean Ponce, and Guillermo Sapiro.
\newblock Online learning for matrix factorization and sparse coding.
\newblock \emph{Journal of Machine Learning Research}, 11\penalty0 (1), 2010.

\bibitem[Nemirovski et~al.(2009)Nemirovski, Juditsky, Lan, and
  Shapiro]{nemirovski2009robust}
Arkadi Nemirovski, Anatoli Juditsky, Guanghui Lan, and Alexander Shapiro.
\newblock Robust stochastic approximation approach to stochastic programming.
\newblock \emph{SIAM Journal on optimization}, 19\penalty0 (4):\penalty0
  1574--1609, 2009.

\bibitem[Nene et~al.(1996)Nene, Nayar, Murase, et~al.]{nene1996columbia}
Sameer~A Nene, Shree~K Nayar, Hiroshi Murase, et~al.
\newblock Columbia object image library (coil-20).
\newblock \emph{Technical report}, 1996.

\bibitem[Nguyen et~al.(2017)Nguyen, Liu, Scheinberg, and
  Takac]{nguyen2017sarah}
Lam~M Nguyen, Jie Liu, Katya Scheinberg, and Martin Takac.
\newblock {SARAH}: A novel method for machine learning problems using
  stochastic recursive gradient.
\newblock In \emph{Proceedings of the 34th International Conference on Machine
  Learning-Volume 70}, pages 2613--2621, 2017.

\bibitem[Ozoli{\c{n}}{\v{s}} et~al.(2013)Ozoli{\c{n}}{\v{s}}, Lai, Caflisch,
  and Osher]{Ozolins2013}
V.~Ozoli{\c{n}}{\v{s}}, R.~Lai, R.~Caflisch, and S.~Osher.
\newblock Compressed modes for variational problems in mathematics and physics.
\newblock \emph{Proceedings of the National Academy of Sciences}, 110\penalty0
  (46):\penalty0 18368--18373, 2013.

\bibitem[Pham et~al.(2019)Pham, Nguyen, Phan, and Tran-Dinh]{pham2019proxsarah}
Nhan~H Pham, Lam~M Nguyen, Dzung~T Phan, and Quoc Tran-Dinh.
\newblock {ProxSARAH}: An efficient algorithmic framework for stochastic
  composite nonconvex optimization.
\newblock \emph{arXiv:1902.05679}, 2019.

\bibitem[Rosasco et~al.(2014)Rosasco, Villa, and
  V{\~u}]{rosasco2014convergence}
Lorenzo Rosasco, Silvia Villa, and Bang~C{\^o}ng V{\~u}.
\newblock Convergence of stochastic proximal gradient algorithm.
\newblock \emph{arXiv preprint arXiv:1403.5074}, 2014.

\bibitem[Sjostrand et~al.(2007)Sjostrand, Rostrup, Ryberg, Larsen, Studholme,
  Baezner, Ferro, Fazekas, Pantoni, Inzitari, et~al.]{sjostrand2007sparse}
Karl Sjostrand, Egill Rostrup, Charlotte Ryberg, Rasmus Larsen, Colin
  Studholme, Hansjoerg Baezner, Jose Ferro, Franz Fazekas, Leonardo Pantoni,
  Domenico Inzitari, et~al.
\newblock Sparse decomposition and modeling of anatomical shape variation.
\newblock \emph{IEEE Transactions on Medical Imaging}, 26\penalty0
  (12):\penalty0 1625--1635, 2007.

\bibitem[Tang and Liu(2012)]{tang2012unsupervised}
Jiliang Tang and Huan Liu.
\newblock Unsupervised feature selection for linked social media data.
\newblock In \emph{Proceedings of the 18th ACM SIGKDD international conference
  on Knowledge discovery and data mining}, pages 904--912. ACM, 2012.

\bibitem[Wang and Lu(2016)]{online-spca-wang-lu-2016}
C.~Wang and Y.~M. Lu.
\newblock Online learning for sparse {PCA} in high dimensions: Exact dynamics
  and phase transitions.
\newblock \emph{https://arxiv.org/pdf/1609.02191.pdf}, 2016.

\bibitem[Wang et~al.(2021)Wang, Liu, Chen, Ma, Xue, and Zhao]{Wang-ManPL-2021}
Z.~Wang, B.~Liu, S.~Chen, S.~Ma, L.~Xue, and H.~Zhao.
\newblock A manifold proximal linear method for sparse spectral clustering with
  application to single-cell rna sequencing data analysis.
\newblock \emph{INFORMS J. Optimization}, 2021.

\bibitem[Wang et~al.(2019)Wang, Ji, Zhou, Liang, and
  Tarokh]{wang2018spiderboost}
Zhe Wang, Kaiyi Ji, Yi~Zhou, Yingbin Liang, and Vahid Tarokh.
\newblock Spider{B}oost and momentum: Faster stochastic variance reduction
  algorithms.
\newblock In \emph{NeurIPS}, 2019.

\bibitem[Wen and Yin(2013)]{Wen-Yin-2013}
Z.~Wen and W.~Yin.
\newblock A feasible method for optimization with orthogonality constraints.
\newblock \emph{Mathematical Programming}, 142\penalty0 (1-2):\penalty0
  397--434, 2013.

\bibitem[Xiao and Zhang(2014)]{xiao2014proximal}
Lin Xiao and Tong Zhang.
\newblock A proximal stochastic gradient method with progressive variance
  reduction.
\newblock \emph{SIAM Journal on Optimization}, 24\penalty0 (4):\penalty0
  2057--2075, 2014.

\bibitem[Xiao et~al.(2018)Xiao, Li, Wen, and Zhang]{xiao2018regularized}
Xiantao Xiao, Yongfeng Li, Zaiwen Wen, and Liwei Zhang.
\newblock A regularized semi-smooth {N}ewton method with projection steps for
  composite convex programs.
\newblock \emph{Journal of Scientific Computing}, 76\penalty0 (1):\penalty0
  364--389, 2018.

\bibitem[Yang et~al.(2014)Yang, Zhang, and Song]{yang2014optimality}
Wei~Hong Yang, Lei-Hong Zhang, and Ruyi Song.
\newblock Optimality conditions for the nonlinear programming problems on
  {R}iemannian manifolds.
\newblock \emph{Pacific Journal of Optimization}, 10\penalty0 (2):\penalty0
  415--434, 2014.

\bibitem[Yang and Xu(2015)]{yang2015streaming}
Wenzhuo Yang and Huan Xu.
\newblock Streaming sparse principal component analysis.
\newblock In \emph{International Conference on Machine Learning}, pages
  494--503. PMLR, 2015.

\bibitem[Yang et~al.(2011)Yang, Shen, Ma, Huang, and Zhou]{Yang2011}
Y.~Yang, H.~Shen, Z.~Ma, Z.~Huang, and X.~Zhou.
\newblock $\ell_{2,1}$-norm regularized discriminative feature selection for
  unsupervised learning.
\newblock In \emph{IJCAI}, volume~22, page 1589, 2011.

\bibitem[Zhang and Sra(2016)]{zhang2016first}
Hongyi Zhang and Suvrit Sra.
\newblock First-order methods for geodesically convex optimization.
\newblock In \emph{Conference on Learning Theory}, pages 1617--1638, 2016.

\bibitem[Zhang et~al.(2018)Zhang, Zhang, and Sra]{zhang2018r}
Jingzhao Zhang, Hongyi Zhang, and Suvrit Sra.
\newblock {R-SPIDER}: A fast {R}iemannian stochastic optimization algorithm
  with curvature independent rate.
\newblock \emph{arXiv preprint arXiv:1811.04194}, 2018.

\bibitem[Zhang et~al.(2017)Zhang, Lau, Kuo, Cheung, Pasupathy, and
  Wright]{Zhang-cvpr-2017}
Y.~Zhang, Y.~Lau, H.-W. Kuo, S.~Cheung, A.~Pasupathy, and J.~Wright.
\newblock On the global geometry of sphere-constrained sparse blind
  deconvolution.
\newblock In \emph{CVPR}, 2017.

\bibitem[Zhou et~al.(2019)Zhou, Yuan, and Feng]{zhou2018faster}
Pan Zhou, Xiao-Tong Yuan, and Jiashi Feng.
\newblock Faster first-order methods for stochastic non-convex optimization on
  {R}iemannian manifolds.
\newblock \emph{TPAMI}, 2019.

\bibitem[Zou et~al.(2006)Zou, Hastie, and Tibshirani]{Zou-spca-2006}
H.~Zou, T.~Hastie, and R.~Tibshirani.
\newblock Sparse principal component analysis.
\newblock \emph{J. Comput. Graph. Stat.}, 15\penalty0 (2):\penalty0 265--286,
  2006.

\end{thebibliography}


\appendix
\renewcommand{\thesection}{\Alph{section}}
\setcounter{section}{0}

\section{Auxiliary Definitions and Lemmas}\label{defs_and_lemmas}

In this section we give a few lemmas and definitions that are necessary to our analysis. These lemmas are proved in existing works, so we do not include the proof here.

\begin{definition}[Generalized Clarke subdifferential, see \cite{hosseini2011generalized}]\label{def-clarke}
For a locally Lipschitz function $F$ on the manifold $\M$, the Riemannian generalized directional derivative $F^{\circ}(X,\zeta)$ at $X\in\M$ in the direction $\zeta$ is defined by
\[
  \limsup\limits_{Y\rightarrow X,t\downarrow 0}\frac{F\circ \phi^{-1}(\phi(Y)+tD\phi(X)[\zeta])-f\circ \phi^{-1}(\phi(Y))}{t}.
\]
Here $(\phi,U)$ is a coordinate chart at $X$.
The Clarke subdifferential $\hat{\partial} F(X)$ at $X\in\M$ is:
\[
   \hat{\partial} F(X)=\{\xi\in \T_X\M :\langle \xi, V\rangle\leq F^{\circ}(X,\zeta),  \forall \zeta\in \T_X\M \}.
\]
\end{definition}

\begin{lemma}\label{batch_lemma}
Suppose $g_i$ is the unbiased and variance-bounded stochastic estimator of $g$ on randomly sampled instance $i$, \textit{i.e.} $\mathbb{E}_i[g_i] = g$ and $\mathbb{E}_i[\|g_i - g\|^2] \leq \sigma^2$. Then we can conclude that the estimator $g_\mathcal{S}:= \frac{1}{|\mathcal{S}|}\sum_{i\in\mathcal{S}} g_i$ based on randomly sampled mini-batch $\mathcal{S}$ is also unbiased and variance-bounded:
\begin{equation}\label{batch_lemma-eq}
\mathbb{E}_\mathcal{S}[g_\mathcal{S}] = g, \quad \mathbb{E}_\mathcal{S}[\|g_\mathcal{S}-g\|^2] \leq \frac{\sigma^2}{|\mathcal{S}|}.
\end{equation}
\end{lemma}


The following lemmas from previous works \citep{absil2009optimization, kasai2018riemannian, zhou2018faster} under Assumptions \ref{ass_retr_all}-\ref{ass_h} regarding retraction and vector transport are very useful.
\begin{lemma}[Retraction $L_R$ smoothness, Lemma 3.5 in \cite{kasai2018riemannian}]\label{retr_smooth_lemma}
If $f(\x)$ has an upper-bounded Hessian, there exists a neighborhood $\mathcal{U}$ of any $X\in \M$ and a constant $L_R>0$ such that $\forall \x, Y\in \mathcal{U}, \Retr_X(\xi) = Y, \xi\in \T_X\M$:
\begin{equation}\label{retr_l_smooth}
f(Y)\leq f(\x) + \langle \nabla f(\x), \xi\rangle + \frac{L_R}{2}\|\xi\|^2.
\end{equation}
\end{lemma}

\begin{lemma}[Lemma 3.7 in \cite{kasai2018riemannian}]\label{lemma-theta}
Under Assumption \ref{ass_retr_all}(ii), there exists a constant $\theta>0$, such that the following inequalities hold for any $\x, Y \in \mathcal{U}$:
\[
\|\Gamma_\eta\xi - P_\eta\xi\|\leq \theta \|\xi\|_\x\|\eta\|_\x, \quad \|\Gamma^{-1}_\eta\xi - P^{-1}_\eta\xi\|\leq \theta \|\chi\|_\x\|\eta\|_\x,
\]
where $\xi,\eta\in \T_\x\M$, $\chi\in \T_Y\M$, $\Retr_\x(\eta)=Y$.
\end{lemma}

\begin{lemma}[Lemma 4 in \cite{zhou2018faster}]\label{zhou_lemma4}
Given $\hat{\x}\in\M$ that does not depend on the update sequence $\{\x_t\}$, the following inequality about the retraction and vector transport holds:
\begin{equation}
    \mathbb{E}_i[\|\Gamma_{\x_t}^{\hat{\x}}(\nabla f_i(\x_t)) - \Gamma_{\x_{t-1}}^{\hat{\x}}\nabla f_i(\x_{t-1})\|^2] \leq 2\Theta^2\|\Retr_{\x_{t-1}}^{-1}(\x_t)\|^2,
\end{equation}
where $\Theta^2 = \theta^2G^2+2(1+c_R)L_H^2$ and $\theta$ is defined in Lemma \ref{lemma-theta}.
\end{lemma}

\begin{lemma} [Lemma 1 in \cite{zhou2018faster}] \label{lemma-zhou}
Let $n_t = \lceil t/q\rceil$, $(n_t-1)q\leq t\leq n_tq$, $t_0 = (n_t-1)q$, where $\lceil a\rceil$ denotes the smallest integer that is larger than $a$. Mini-batches $\mathcal{S}^1_t$,  $\mathcal{S}_t^2$ are selected as described in Algorithm \ref{prox-r-spiderboost-algo-chart}. Under the Assumptions \ref{ass_retr_all}-\ref{ass_h}, the estimation error between the R-SARAH estimator $V_t$ generated by Algorithm \ref{prox-r-spiderboost-algo-chart} and full gradient $\nabla f(X_t)$ is bounded by:
\[
    \mathbb{E}[\|\v_t - \nabla f(X_t) \|^2] \leq I\{|\mathcal{S}_t^1|<n\}\cdot\frac{\sigma^2}{|\mathcal{S}_t^1|} + \sum_{i=t_0}^{t-1}\frac{\Theta^2}{|\mathcal{S}_t^2|}\mathbb{E}[\|\Retr_{\x_{i}}^{-1}(\x_{i+1})\|^2],
\]
where $I\{\cdot\}$ denotes an indicator function.
\end{lemma}

For the ease of presentation, we adopt the following notation, which is consistent with the ones used in \eqref{subproblem-ManPG} and \eqref{subproblem}.
\begin{align}
\label{def-zeta} \zeta_t & := \argmin_{\zeta\in \T_{\x_t}\mathcal{M}} \{\phi_t(\zeta):=\langle \v_t,\zeta\rangle + \frac{1}{2\gamma}\|\zeta\|^2 + h(\x_t+\zeta)\}, \\
\label{def-xi}   \xi_t   & := \argmin_{\xi  \in \T_{\x_t}\mathcal{M}} \{\langle \nabla f(\x_t),\xi\rangle + \frac{1}{2\gamma}\|\xi\|^2 + h(\x_t+\xi)\}.
\end{align}
Moreover, note that according to the definition of $\mathcal{F}_t$, when we take conditional expectation $\mathbb{E}[\cdot\mid\mathcal{F}_{t-1}]$, $X_t$ in both R-ProxSGD and R-ProxSPB has been realized.

\section{Necessary Lemmas for Proving Theorem \ref{theorem_prox_r_sgd}}

\begin{lemma}\label{phi-t-descent}
The solution $\mathbf{\zeta}_t$ defined in \eqref{def-zeta} satisfies:
\begin{equation}\label{ssn_prop}
    \mathbb{E}[\phi_t(\eta_t\mathbf{\zeta}_t)|\mathcal{F}_{t-1}] -\phi_t(0) \leq \frac{(\eta_t-2)\eta_t}{2\gamma}\mathbb{E}[\|\mathbf{\zeta}_t\|^2|\mathcal{F}_{t-1}].
\end{equation}
\end{lemma}

\proof 
Note that $\phi_t(\mathbf{\zeta})$ is $(1/\gamma)$-strongly convex with respect to $\mathbf{\zeta}$. For $\mathbf{\zeta}_1, \mathbf{\zeta}_2 \in \T_{X_t}\mathcal{M}$, we have:
\begin{equation}\label{sc}
    \phi_t(\mathbf{\zeta}_2) \geq \phi_t(\mathbf{\zeta}_1) +\langle \hat{\partial} \phi_t(\mathbf{\zeta}_1), \mathbf{\zeta}_2-\mathbf{\zeta}_1\rangle + \frac{1}{2\gamma}\| \mathbf{\zeta}_2-\mathbf{\zeta}_1\|^2.
\end{equation}
Note that the optimality conditions of \eqref{def-zeta} are given by $0\in\mathrm{Proj}_{\T_{X_t}\mathcal{M}}\partial \phi_t(\mathbf{\zeta}_t)$. Therefore,
\begin{equation}\label{sc-1}
\langle \hat{\partial} \phi_t(\mathbf{\zeta}_1), \mathbf{\zeta}_2-\mathbf{\zeta}_1\rangle = \langle \mathrm{Proj}_{\T_{X_t}\mathcal{M}}  \partial \phi_t(\mathbf{\zeta}_1), \mathbf{\zeta}_2-\mathbf{\zeta}_1\rangle=0, \forall \mathbf{\zeta}_1, \mathbf{\zeta}_2 \in \T_{X_t}\mathcal{M}.
\end{equation}
Letting $\mathbf{\zeta}_1=\mathbf{\zeta}_t$ and $\mathbf{\zeta}_2 = 0$ in \eqref{sc}, and combining with \eqref{sc-1}, we have
\[
    \phi_t(0)\geq \phi_t(\mathbf{\zeta}_t)  + \frac{1}{2\gamma}\|\mathbf{\zeta}_t\|^2,
\]
which is equivalent to:
\begin{equation}\label{sc2}
    h(X_t+\mathbf{\zeta}_t) - h(X_t) \leq \langle -\v_t, \mathbf{\zeta}_t\rangle -\frac{1}{\gamma}\|\mathbf{\zeta}_t\|^2.
\end{equation}
According to the definition of $\phi_t$, $\phi_t(\eta_t\mathbf{\zeta}_t) -\phi_t(0)$ can be written as:
\begin{equation}\label{phi_t_def}
    \phi_t(\eta_t\mathbf{\zeta}_t) -\phi_t(0) = \eta_t\langle \v_t, \mathbf{\zeta}_t\rangle +\frac{\eta_t^2}{2\gamma}\|\mathbf{\zeta}_t\|^2 + h(X_t+\eta_t\mathbf{\zeta}_t)-h(X_t).
\end{equation}
From \eqref{sc2} and the convexity of $h$: $ h(X_t+\eta_t\mathbf{\zeta}_t) \leq \eta_t h(X_t+\mathbf{\zeta}_t) + (1-\eta_t)h(X_t)$, $\eta_t\in (0,1]$, we have
\begin{equation}\label{convexity}
    h(X_t+\eta_t\mathbf{\zeta}_t)-h(X_t) \leq -\eta_t\langle \v_t, \mathbf{\zeta}_t\rangle - \frac{\eta_t}{\gamma}\|\mathbf{\zeta}_t\|^2.
\end{equation}
Combine \eqref{phi_t_def} and \eqref{convexity} and take expectation conditioned on $\mathcal{F}_{t-1}$ on both sides, we get the desired result. 
\endproof

The following lemma justifies why $G(\x, \nabla f(\x), \gamma)$ is valid for defining the $\epsilon$-stationary solution.
\begin{lemma}\label{G-stationary-valid}
If $G(\x, \nabla f(\x), \gamma)  = 0$, and the retraction is given by the Polar decomposition: $\Retr_\x(\xi) = (\x+\xi)(I+\xi^\top\xi)^{-\frac{1}{2}}$, then $\x$ is a stationary point of problems \eqref{prob_1}, i.e., $0\in \nabla f(\x) + \mathrm{Proj}_{\T_{\x}\mathcal{M}}\partial h(\x)$.
\end{lemma}

To prove Lemma \ref{G-stationary-valid}, we first need to show the following Lemma.
\begin{lemma}\label{G-stationary-valid-affiliate}
Consider $\x\in\mathcal{M}$, $\mathcal{M}$ is the Stiefel manifold and $\xi\in T_\x\mathcal{M}$. If $\x = \Retr_\x(\xi)$, where the retraction is given by the Polar decomposition: $\Retr_\x(\xi) = (\x+\xi)(I+\xi^\top\xi)^{-\frac{1}{2}}$, then $\xi = \mathbf{0}_\x$.
\end{lemma}

\proof 
If $\x =\Retr_\x(\xi) = (\x+\xi)(I+\xi^\top\xi)^{-\frac{1}{2}}$, then we have
\begin{equation}\label{solve_eq-1}
\x+\xi = \x(I+\xi^\top\xi)^{\frac{1}{2}}.
\end{equation}
Since $\x^\top\x =I$, \eqref{solve_eq-1} leads to
\begin{equation}\label{solve_eq-2}
\x^\top\x + \xi^\top\x = (I+\xi^\top\xi)^{\frac{1}{2}}
\end{equation}
and
\begin{equation}\label{solve_eq-3}
\x^\top\x + \x^\top\xi = (I+\xi^\top\xi)^{\frac{1}{2}}.
\end{equation}
Since $\xi\in \T_\x\mathcal{M}$, we have $\xi^\top\x + \x^\top\xi = 0$. Adding \eqref{solve_eq-2} and \eqref{solve_eq-3} gives $2I = 2(I+\xi^\top\xi)^{\frac{1}{2}}$, which implies $\xi = \mathbf{0}_\x$. 
\endproof

Now we are ready to give the proof of Lemma \ref{G-stationary-valid}.
\proof 
If $G(\x_t, \nabla f(\x_t), \gamma)  = 0$, we have $\xi_t = \mathbf{0}_{\x_t}$ because of Lemma \ref{G-stationary-valid-affiliate}. According to \cite{yang2014optimality}, the optimality conditions of \eqref{def-xi} are given by
\[
    0\in \nabla f(\x_t) + \frac{1}{\gamma}\xi_t + \mathrm{Proj}_{\T_{\x_t}\mathcal{M}}\partial h(\x_t+\xi_t).
\]
Thus, $G(\x_t, \nabla f(\x_t), \gamma)  = 0$ leads to that $0\in \nabla f(\x_t) + \mathrm{Proj}_{\T_{\x_t}\mathcal{M}}\partial h(\x_t)$, which means $\x_t$ is a stationary point of problem \eqref{prob_1}. 
\endproof

The following lemma shows the progress of the algorithm in one iteration in terms of objective function value.

\begin{lemma}\label{F-descent}
Denote $X_t^+:=X_t+\eta_t\mathbf{\zeta}_t$. The following inequality holds:
\[
 F(X_{t+1})\!-\!F(X_t) \!\leq\! \frac{(L_R\gamma - 1)\eta_t^2}{2\gamma}\|\mathbf{\zeta}_t\|^2  \!+\! h(X_{t+1})\!-\! h(X_t^+) + \phi_t(\eta_t\mathbf{\zeta}_t)-\phi_t(0) + \eta_t\langle \nabla f(X_t)-\v_t, \mathbf{\zeta}_t\rangle.
\]
\end{lemma}

\proof 
Consider the update $X_{t+1} = \Retr_{X_t} (\eta_t\zeta_t)$. By applying Lemma \ref{retr_smooth_lemma} with $\x = X_t, Y = X_{t+1}$ and $\xi = \eta_t\zeta_t$, we get
\[
f(X_{t+1})-f(X_t) \leq \eta_t\langle \nabla f(X_t), \zeta_t\rangle + \frac{L_R\eta_t^2}{2}\|\zeta_t\|^2,
\]
which leads to:
\begin{eqnarray}\label{difference_of_F}
    F(X_{t+1})-F(X_t) \leq \frac{L_R\eta_t^2}{2}\|\mathbf{\zeta}_t\|^2 + \eta_t\langle \nabla f(X_t), \mathbf{\zeta}_t\rangle + h(X_{t+1})-h(X_t).
\end{eqnarray}
Denote $X_t^+:=X_t+\eta_t\mathbf{\zeta}_t$. The definition of $\phi_t$ indicates:
\begin{equation}\label{previous}
    \eta_t\langle \v_t, \mathbf{\zeta}_t\rangle = \phi_t(\eta_t\mathbf{\zeta}_t)-\phi_t(0) - \frac{\eta_t^2}{2\gamma}\|\mathbf{\zeta}_t\|^2 - h(X_t^+) + h(X_t).
\end{equation}
Combining \eqref{difference_of_F} and \eqref{previous} gives the desired result. 
\endproof

The following lemma gives an upper bound to the size of $G(\x_t, \nabla f(\x_t), \gamma)$.

\begin{lemma}\label{G-bounded}
With $\zeta_t$ and $\xi_t$ defined in \eqref{def-zeta} and \eqref{def-xi}, for $G(\x_t, \nabla f(\x_t), \gamma) = \frac{1}{\gamma}(\x_t - \mathrm{Retr}_{\x_t}(\xi_t))$, it holds that
\begin{equation}\label{G-bounded-eq}
\|G(\x_t, \nabla f(\x_t), \gamma)\|^2 \!\leq\! 2M_1^2(7\|\zeta_t\|^2 \!+\! 4\gamma \|V_t\!-\! \nabla f(\x_t)\|^2).
\end{equation}
\end{lemma}

\proof 
Let $G(\x_t, V_t, \gamma) = \frac{1}{\gamma}(\x_t - \mathrm{Retr}_{\x_t}(\gamma\zeta_t))$. We first have the following trivial inequality:
\begin{eqnarray}\label{start_eq}
    \|G(\x_t, \nabla f(\x_t), \gamma)\|^2 \leq 2\|G(\x_t, V_t, \gamma)\|^2+ 2\|G(\x_t, V_t, \gamma) - G(\x_t, \nabla f(\x_t), \gamma)\|^2.
\end{eqnarray}
The first term on the right hand side of \eqref{start_eq} can be bounded based on the property of retraction in Assumption \ref{ass_retr_all} (iii):
\begin{equation}\label{start_eq-1}
\|G(\x_t, V_t, \gamma)\|^2 = \frac{1}{\gamma^2}\|\x_t - \mathrm{Retr}_{X_t}(\gamma \zeta_t)\|^2 \leq M_1^2\|\zeta_t\|^2.
\end{equation}
The second term on the right hand side of \eqref{start_eq} can be bounded as:
\begin{eqnarray}
    \|G(\x_t, V_t, \gamma) - G(\x_t, \nabla f(\x_t), \gamma)\|^2 \leq \frac{2\|\x_t - \mathrm{Retr}_{\x_t}(\gamma\zeta_t)\|^2}{\gamma^2}+ \frac{2}{\gamma^2}\|\x_t - \mathrm{Retr}_{\x_t}(\gamma\xi_t)\|^2,
\end{eqnarray}
which further implies
\begin{eqnarray}
    \|G(\x_t, V_t, \gamma) - G(\x_t, \nabla f(\x_t), \gamma)\|^2 \leq 2M_1^2(\|\zeta_t\|^2 + \|\xi_t\|^2)\leq 2M_1^2(3\|\zeta_t\|^2 + 2\|\xi_t-\zeta_t\|^2).
\end{eqnarray}
The optimality conditions of \eqref{def-zeta} and \eqref{def-xi} are given by (see \cite{yang2014optimality}):
\begin{align}
\label{opt-zeta} &  0 \in V_t + \frac{1}{\gamma}\zeta_t + \mathrm{Proj}_{T_{\x_t}\mathcal{M}}\partial h(\x_t+\zeta_t), \\
\label{opt-xi}   &  0 \in \nabla f(X_t) + \frac{1}{\gamma}\xi_t + \mathrm{Proj}_{T_{\x_t}\mathcal{M}}\partial h(\x_t+\xi_t).
\end{align}
Let $\x_t^\dagger = \x_t +\xi_t$ and $\x_t^+ = \x_t + \zeta_t$. \eqref{opt-zeta} and \eqref{opt-xi} indicate that for any $\mathbf{u} \in \T_{\x_t}\mathcal{M}$, there exist $p^+\in \partial h(\x_t^+)$ and $p^\dagger\in \partial h(\x_t^\dagger)$ such that
\begin{eqnarray}
& \langle \frac{1}{\gamma} \zeta_t + V_t + \mathrm{Proj}_{T_{\x_t}\mathcal{M}} p^+, \mathbf{u} - \x_t^+ \rangle \geq 0, \label{proof-lem3-1}\\
& \langle \frac{1}{\gamma} \xi_t + \nabla f(\x_t) + \mathrm{Proj}_{T_{\x_t}\mathcal{M}} p^\dagger, \mathbf{u} - \x_t^\dagger \rangle \geq 0. \label{proof-lem3-2}
\end{eqnarray}
Let $\mathbf{u} = \x_t^\dagger$ in \eqref{proof-lem3-1} and $\mathbf{u} = \x_t^+$ in \eqref{proof-lem3-2}. Since $\x_t^\dagger - \x_t^+$ and $\x_t^+ - \x_t^\dagger$ both lie in $\T_{\x_t}\mathcal{M}$, we have $\langle \mathrm{Proj}_{T_{\x_t}\mathcal{M}}p^+, \x_t^+ - \x_t^\dagger\rangle =  \langle p^+, \x_t^+ - \x_t^\dagger\rangle$ and $\langle \mathrm{Proj}_{T_{\x_t}\mathcal{M}}p^\dagger, \x_t^\dagger - \x_t^+\rangle =  \langle p^\dagger, \x_t^\dagger - \x_t^+\rangle$. Therefore, \eqref{proof-lem3-1} and \eqref{proof-lem3-2} reduce to:
\begin{eqnarray}
& \langle \frac{1}{\gamma} \zeta_t + V_t +  p^+, \x_t^\dagger - \x_t^+ \rangle \geq 0, \label{proof-lem3-3} \\
& \langle \frac{1}{\gamma} \xi_t + \nabla f(\x_t) + p^\dagger, \x_t^+ - \x_t^\dagger \rangle \geq 0. \label{proof-lem3-4}
\end{eqnarray}
By using the convexity of $h(\x)$, we have $\langle p^+, \x_t^+ - \x_t^\dagger\rangle \geq h(\x_t^+) - h(\x_t^\dagger)$, and $\langle p^\dagger, \x_t^\dagger - \x_t^+\rangle\geq h(\x_t^\dagger) - h(\x_t^+)$. Therefore, \eqref{proof-lem3-3} and \eqref{proof-lem3-4} reduce to:
\begin{eqnarray}
& \langle V_t, \x_t^\dagger - \x_t^+ \rangle \geq \frac{1}{\gamma}\langle\zeta_t, \x_t^+ - \x_t^\dagger \rangle + h(\x_t^+) - h(\x_t^\dagger),\label{proof-lem3-5}\\
& \langle \nabla f(\x_t), \x_t^+ - \x_t^\dagger \rangle \geq \frac{1}{\gamma}\langle\xi_t, \x_t^\dagger - \x_t^+ \rangle +  h(\x_t^\dagger) - h(\x_t^+).\label{proof-lem3-6}
\end{eqnarray}
Summing up \eqref{proof-lem3-5} and \eqref{proof-lem3-6} gives: (note that $\x_t^\dagger -\x_t^+ = \xi_t - \zeta_t$):
\begin{eqnarray}
\|V_t - \nabla f(\x_t)\|\|\x_t^+- \x_t^\dagger\| \geq \langle V_t - \nabla f(\x_t), \x_t^\dagger - \x_t^+\rangle  \geq \frac{1}{\gamma}\langle \xi_t - \zeta_t, \x_t^\dagger -\x_t^+\rangle = \frac{1}{\gamma}\|\x_t^\dagger -\x_t^+\|^2,
\end{eqnarray}
which further implies $\|\xi_t - \zeta_t\| = \|\x_t^\dagger -\x_t^+\| \leq \gamma \|V_t - \nabla f(\x_t)\|$. We hence have:
\[
\|G(\x_t, V_t, \gamma)  -  G(\x_t, \nabla f(\x_t), \gamma)\|^2 \leq 2M_1^2(3\|\zeta_t\|^2 + 2\|\xi_t-\zeta_t\|^2) \leq 6M_1^2\|\zeta_t\|^2 + 4M_1^2 \gamma \|V_t - \nabla f(\x_t)\|^2,
\]
which combining with \eqref{start_eq} and \eqref{start_eq-1} completes the proof. 
\endproof

The following lemma shows the progress of R-ProxSGD in one iteration in terms of objective function value.
\begin{lemma}\label{F-descent-SGD}
The sequences $\{X_t\}_{t=1}^{T+1}$ and $\{\mathbf{\zeta}_t\}_{t=1}^{T}$ generated by R-ProxSGD (Algorithm \ref{prox-r-sgd-algo-chart}) satisfy the following inequality: 
\begin{equation}\label{lem:F-descent-SGD-eq}
\mathbb{E}[F(X_{t+1}) \!-\! F(X_t)]\!\leq\! \big(\tilde{L}\eta_t^2-\frac{1}{\gamma}\eta_t+\frac{1}{2}\big)\mathbb{E}[\|\mathbf{\zeta}_t\|^2]  \!+\!\frac{\eta_t^2\sigma^2}{2|\mathcal{S}_t|},
\end{equation}
where $\tilde{L} = (L_R/2+L_hM_2)$.
\end{lemma}

\proof 
Denote $\x_t^+=\x_t+\eta_t\zeta_t$. Assumptions \ref{ass_retr_all}(iii) and \ref{ass_h} yield the following inequalities:
\[h(X_{t+1})-h(X_t^+) \leq L_h\|X_{t+1}-X_t^+\|\leq L_hM_2\eta_t^2\|\mathbf{\zeta}_t\|^2,\]
which together with Lemma \ref{F-descent} and Young's inequality gives
\begin{eqnarray}\label{re-written_difference_of_F}
F(X_{t+1})-F(X_t) \leq \left(\frac{L_R\eta_t^2}{2}-\frac{\eta_t^2}{2\gamma}+L_hM_2\eta_t^2+\frac{1}{2}\right)\|\mathbf{\zeta}_t\|^2 +\frac{\eta_t^2}{2}\|\nabla f(X_t)-\v_t\|^2 + \phi_t(\eta_t\mathbf{\zeta}_t)-\phi_t(0).
\end{eqnarray}
Taking expectation conditioned on $\mathcal{F}_{t-1}$ to both side of \eqref{re-written_difference_of_F}, we get:
\begin{align}\label{cond_expected_difference_of_F}
    \mathbb{E}[F(X_{t+1})\mid \mathcal{F}_{t-1}] -F(X_t)\leq & \left(\bar{L}\eta_t^2+\frac{1}{2}\right)\mathbb{E}[\|\mathbf{\zeta}_t\|^2\mid \mathcal{F}_{t-1}]  +\frac{\eta_t^2}{2}\mathbb{E}[\|\nabla f(X_t)-\v_t\|^2 \mid \mathcal{F}_{t-1}] \\ & +\mathbb{E}[\phi_t(\eta_t\mathbf{\zeta}_t)\mid\mathcal{F}_{t-1}]-\phi_t(0),\nonumber
\end{align}
where $\bar{L} := \frac{L_R}{2}-\frac{1}{2\gamma}+L_hM_2$.
Using Lemma \ref{phi-t-descent} and taking the whole expectation on both sides of \eqref{cond_expected_difference_of_F} completes the proof. 
\endproof

{
\section{Proof of Theorem \ref{theorem_prox_r_sgd}}

We can re-arrange terms in \eqref{lem:F-descent-SGD-eq} as follows for $0<\eta_t\leq 1$ (note $|\mathcal{S}_t|=s$ for all $t$):
\begin{equation}\label{final}
    \left(\frac{1}{\gamma}\eta_t-\tilde{L}\eta_t^2 -\frac{1}{2}\right)\mathbb{E}[\|\mathbf{\zeta}_t\|^2] \leq\mathbb{E}[F(X_{t})] -\mathbb{E}[F(X_{t+1})]  +\frac{\eta_t^2\sigma^2}{2s}.
\end{equation}
If we choose $\gamma$ small enough such that $\gamma\leq\frac{2\eta_t}{2\tilde{L}\eta_t^2+\eta_t+1}$ for all $t=0,\ldots,T-1$, then
\begin{equation}\label{gamma}
\frac{1}{\gamma}\eta_t-\tilde{L}\eta_t^2 -\frac{1}{2}\geq \frac{\eta_t}{2}, t=0,\ldots,T-1.
\end{equation}
Combining \eqref{final} and \eqref{gamma} yields:
\begin{equation}\label{final-1}
    \frac{1}{2}\mathbb{E}[\|\mathbf{\zeta}_t\|^2] \leq\frac{\mathbb{E}[F(X_{t})] -\mathbb{E}[F(X_{t+1})]}{\eta_t}  +\frac{\eta_t\sigma^2}{2s}. 
\end{equation}
We choose $\eta_t$ as a constant $\eta_t = \eta \in (0,1)$. 
Denote $\Delta_0:= F(X_{0}) - F(X^*)$, where $X^*$ is a global optimal solution to the problem \eqref{prob_1}. Summing up \eqref{final-1} for $t=0,\ldots,T-1$ and dividing both sides by $T$, we get: 
\begin{equation}\label{proof-theorem5-1}
        \frac{1}{T}\sum_{t=0}^{T-1}\mathbb{E}[\|\mathbf{\zeta}_t\|^2]\leq \frac{2\Delta_0}{T\eta} + \frac{\sigma^2\eta}{s}, 
\end{equation}
where we used the fact that $\mathbb{E}[F(X_{0})] -\mathbb{E}[F(X_T)] \leq \Delta_0$.
Moreover, \eqref{batch_lemma-eq} indicates that
\begin{equation}\label{proof-theorem5-3}
\frac{1}{T}\sum_{t=0}^{T-1}\mathbb{E}[\|V_t - \nabla f(\x_t)\|^2 \leq \frac{\sigma^2}{s}.
\end{equation}
Combining \eqref{G-bounded-eq}, \eqref{proof-theorem5-1} and \eqref{proof-theorem5-3} yields:
\[
\frac{1}{T}\sum_{t=0}^{T-1}\mathbb{E}[\|G(\x_t, \nabla f(\x_t), \gamma)\|^2] \leq 14M_1^2\Big(\frac{2\Delta_0}{T\eta} + \frac{\sigma^2\eta}{s}\Big) + 8M_1^2\gamma\sigma^2/s,
\]
which together with Jensen's inequality and the convexity of $\|\cdot\|^2$ implies that:
\begin{align}
\label{proof-theorem5-2} & \left(\mathbb{E}\left[\frac{1}{T}\sum_{t=0}^{T-1}\|G(\x_t, \nabla f(\x_t), \gamma)\|\right]\right)^2 \\
\nonumber           \leq & \mathbb{E}\left[\left(\frac{1}{T}\sum_{t=0}^{T-1}\|G(\x_t, \nabla f(\x_t), \gamma)\|\right)^2\right] \\
\nonumber           \leq & \frac{1}{T}\sum_{t=0}^{T-1}\mathbb{E}[\|G(\x_t, \nabla f(\x_t), \gamma)\|^2] \\
\nonumber           \leq & 14M_1^2\Big(\frac{2\Delta_0}{T\eta} + \frac{\sigma^2\eta}{s}\Big) + 8M_1^2\gamma\sigma^2/s.
\end{align}
By setting $s = (M_1^2\sigma^2(28\eta+16\gamma))\epsilon^{-2}$, we know that as long as
\begin{equation}\label{inequality-T}
T \geq \frac{56M_1^2\Delta_0}{\eta\epsilon^{2}},
\end{equation}
the right hand side of \eqref{proof-theorem5-2} is upper bounded by $\epsilon^2$, that is:
\begin{equation}\label{proof-theorem5-4}
\left(\mathbb{E}\left[\frac{1}{T}\sum_{t=0}^{T-1}\|G(\x_t, \nabla f(\x_t), \gamma)\|\right]\right) \leq \epsilon.
\end{equation}
Therefore, for an index $\nu$ that is uniformly sampled from $\{0,\ldots,T-1\}$, we have 
\[
\mathbb{E}[\|G(\x_\nu, \nabla f(\x_\nu), \gamma)\|] \leq \epsilon,
\] 
i.e., $\x_\nu$ is an $\epsilon$-stochastic stationary point of problem \eqref{prob_1}. Condition \eqref{inequality-T} shows that the number of iterations needed by R-ProxSGD for obtaining an $\epsilon$-stochastic stationary point is $T=O(\epsilon^{-2})$, which immediately implies that the total IFO complexity is $O(\epsilon^{-4})$. This completes the proof of Theorem \ref{theorem_prox_r_sgd}.


}

\section{Necessary Lemma for Proving Theorem \ref{theorem_spiderboost}}

Similar to Lemma \ref{F-descent-SGD}, the following lemma gives the progress of R-ProxSPB in one iteration in terms of the objective function value.

\begin{lemma}\label{F-descent-SPB}
The sequences $\{X_t\}_{t=1}^{T+1}$ and $\{\mathbf{\zeta}_t\}_{t=1}^T$ generated by R-ProxSPB (Algorithm \ref{prox-r-spiderboost-algo-chart}) satisfy the following inequality:
\begin{equation}\label{F-descent-SPB-eq}
\mathbb{E}[F(X_{t+1})-F(X_t)]\leq \eta\big(\tilde{L}\eta-\frac{1}{\tilde{\gamma}}\big)\mathbb{E}[\|\mathbf{\zeta}_t\|^2] + I\{|\mathcal{S}_t^1|<n\}\frac{\eta\sigma^2}{2|\mathcal{S}_t^1|}+ \sum_{i=(n_t-1)q}^t\frac{\Theta^2\eta^3c_E}{2|\mathcal{S}_t^2|}\mathbb{E}[\|\mathbf{\zeta}_i\|^2],
\end{equation}
where $\tilde{L} = L_R/2+L_hM_2$ and $\tilde{\gamma} = \frac{2\gamma}{2-\gamma}$.
\end{lemma}

\proof 
Similar to the proof of Lemma \ref{F-descent-SGD}, by using Lemma \ref{F-descent}, Assumptions \ref{ass_retr_all}(iii) and \ref{ass_h}, and Young's inequality, we have:
\begin{align}\label{re-written_difference_of_F_spider}
    & F(X_{t+1})-F(X_t) \\
    \nonumber \leq & (\frac{L_R\eta^2}{2}-\frac{\eta^2}{2\gamma}+L_hM_2\eta^2+\frac{\eta}{2})\|\mathbf{\zeta}_t\|^2 +\frac{\eta}{2}\|\v_t - \nabla f(X_t)\|^2 + \phi_t(\eta\mathbf{\zeta}_t)-\phi_t(0).
\end{align}
Taking conditional expectation on both sides of \eqref{re-written_difference_of_F_spider} conditioned on $\mathcal{F}_{t-1}$, we have:
\begin{align}\label{cond_expected_difference_of_F_spider}
    & \mathbb{E}[F(X_{t+1})\mid\mathcal{F}_{t-1}] -F(X_t) \\
    \nonumber \leq & \eta(\bar{L}\eta+\frac{1}{2})\mathbb{E}[\|\mathbf{\zeta}_t\|^2\mid \mathcal{F}_{t-1}] +\frac{\eta}{2}\mathbb{E}[\|\v_t - \nabla f(X_t)\|^2 \mid \mathcal{F}_{t-1}] +\mathbb{E}[ \phi_t(\eta\mathbf{\zeta}_t)\mid \mathcal{F}_{t-1}]-\phi_t(0),
\end{align}
where $\bar{L} := \frac{L_R}{2}-\frac{1}{2\gamma}+L_hM_2$. Taking the whole expectation on both sides of \eqref{cond_expected_difference_of_F_spider} yields:
\begin{align*}
 & \mathbb{E}[F(X_{t+1})-F(X_t)] \\
 \stackrel{(i)}{\leq}   & \eta\big(\tilde{L}\eta-\frac{1}{\tilde{\gamma}}\big)\mathbb{E}[\|\mathbf{\zeta}_t\|^2]  +\frac{\eta}{2}\mathbb{E}[\|\v_t - \nabla f(X_t)\|^2]\\
 \stackrel{(ii)}{\leq}  & \eta\big(\tilde{L}\eta - \frac{1}{\tilde{\gamma}}\big)\mathbb{E}[\|\mathbf{\zeta}_t\|^2]  +  I\{|\mathcal{S}_t^1|<n\}\frac{\eta\sigma^2}{2|\mathcal{S}_t^1|}  +  \sum_{i=t_0}^t\frac{\Theta^2\eta}{2|\mathcal{S}_t^2|}\mathbb{E}[\|\Retr^{-1}_{\x_i}(\x_{i+1)}\|^2]\\
 \stackrel{(iii)}{\leq} & \eta\big(\tilde{L}\eta-\frac{1}{\tilde{\gamma}}\big)\mathbb{E}[\|\mathbf{\zeta}_t\|^2] +I\{|\mathcal{S}_t^1|<n\}\frac{\eta\sigma^2}{2|\mathcal{S}_t^1|} + \sum_{i=t_0}^t\frac{\Theta^2\eta^3c_E}{2|\mathcal{S}_t^2|}\mathbb{E}[\|\mathbf{\zeta}_i\|^2],
\end{align*}
where (i) is from Lemma \ref{phi-t-descent}, (ii) is due to Lemma \ref{lemma-zhou}, and (iii) is due to the update $X_{t+1} = \Retr_{X_t}(\eta\mathbf{\zeta}_t)$ and the Assumption \ref{ass_retr_all}(ii). This completes the proof. 
\endproof

\section{Proof of Theorem \ref{theorem_spiderboost}}

Let $n_t=\lceil t/q\rceil, t_0 = (n_t-1)q$. Since the length of recursion of $\v_t$ is $q$ in R-ProxSPB, we calculate the telescoping sum of \eqref{F-descent-SPB-eq} from $t_0=(n_t-1)q$ to $t+1\leq n_tq$:
\begin{align}\label{telescopic_q}
 & \mathbb{E}[F(X_{t+1})-F(X_{t_0})] \\
 \nonumber \leq & \eta\left(\tilde{L}\eta-\frac{1}{\tilde{\gamma}}\right)\sum_{i=t_0}^t\mathbb{E}[\|\mathbf{\zeta}_i\|^2] + \sum_{i=t_0}^t I\{|\mathcal{S}_t^1|<n\}\frac{\eta\sigma^2}{2|\mathcal{S}_t^1|} + \frac{\Theta^2 c_E\eta^3}{2|\mathcal{S}_t^2|}\sum_{j=t_0}^t\sum_{i=t_0}^j\mathbb{E}[\|\mathbf{\zeta}_i\|^2].
\end{align}
By noting $\sum_{j=t_0}^t\sum_{i=t_0}^j\mathbb{E}[\|\mathbf{\zeta}_i\|^2]\leq q\sum_{i=t_0}^t\mathbb{E}[\|\mathbf{\zeta}_i\|^2]$, $\tilde{\gamma}=2\gamma/(2-\gamma)=1/2$ (since $\gamma=2/5$), and $|\mathcal{S}_t^2|=q$ for all $t$, \eqref{telescopic_q} can be reduced to:
\begin{align}\label{telescopic_q_rewritten}
\mathbb{E}[F(X_{t+1})-F(X_{t_0})] \leq \sum_{i=t_0}^t I\{|\mathcal{S}_t^1|<n\}\frac{\eta\sigma^2}{2|\mathcal{S}_t^1|}+ \eta\left(\frac{c_E\Theta^2\eta^2}{2} + \tilde{L}\eta-2\right)\sum_{i=t_0}^t\mathbb{E}[\|\mathbf{\zeta}_i\|^2].
\end{align}
Moreover, the choice of $\eta$: $0< \eta \leq (-\tilde{L}+\sqrt{\tilde{L}^2+2c_E\Theta^2})/(c_E\Theta^2)$ guarantees that
\[\frac{c_E\Theta^2\eta^2}{2} + \tilde{L}\eta-2\leq -1.\]
Therefore, \eqref{telescopic_q_rewritten} reduces to
\begin{align}\label{telescopic_q_rewritten-1}
\eta\sum_{i=t_0}^t\mathbb{E}[\|\mathbf{\zeta}_i\|^2] \leq -\mathbb{E}[F(X_{t+1})-F(X_{t_0})] + \sum_{i=t_0}^t I\{|\mathcal{S}_t^1|<n\}\frac{\eta\sigma^2}{2|\mathcal{S}_t^1|}.
\end{align}

\subsection{Finite-sum case}
In the finite-sum case, we have $|\mathcal{S}_t^1| = n$, which implies that $I\{|\mathcal{S}_t^1|<n\}=0$. Therefore, \eqref{telescopic_q_rewritten-1} reduces to:
\begin{equation}\label{finite-sum-1}
    \sum_{i=t_0}^t\mathbb{E}[\|\mathbf{\zeta}_i\|^2] \leq \frac{\mathbb{E}[F(X_{(n_t-1)q})- F(X_{t+1})]}{\eta}.
\end{equation}
We now calculate the telescoping sum for \eqref{finite-sum-1} for all length-$q$ epochs that $t+1 = q, 2q, ..., Kq$ ($K= \lfloor \frac{T}{q}\rfloor$) and the telescoping sum from $t=Kq+1$ to $T-1$. This results in:
\begin{align}\label{proof-spb-0}
\frac{1}{T}\sum_{t=0}^{T-1}\mathbb{E}[\|\mathbf{\zeta}_t\|^2] =  \frac{1}{T}\left(\sum_{t=0}^{Kq-1}\mathbb{E}[\|\mathbf{\zeta}_t\|^2]+\sum_{t=Kq}^{T-1}\mathbb{E}[\|\mathbf{\zeta}_t\|^2]\right)\leq \frac{ \mathbb{E}[F(X_0)-F(X_T)]}{T\eta}\leq \frac{\Delta_0}{\eta T}.
\end{align}
Moreover, Lemma \ref{lemma-zhou} yields that
\begin{equation}\label{proof-spb-1}
     \mathbb{E}[\|V_t - \nabla f(X_t) \|^2] \leq \sum_{i=t_0}^{t-1}\frac{\Theta^2 c_E^2\eta^2}{q}\mathbb{E}[\|\zeta_i\|^2].
\end{equation}
Summing up \eqref{proof-spb-1} over $t= 0,\ldots,T-1$, we get
\begin{align}\label{proof-spb-2}
\frac{1}{T}\sum_{t=0}^{T-1}\mathbb{E}[\|V_t - \nabla f(X_t) \|^2] \leq\frac{1}{T} \sum_{t=1}^{T-1} \sum_{i=t_0}^{t-1}\frac{\Theta^2 c_E^2\eta^2}{q}\mathbb{E}[\|\zeta_i\|^2]\leq \frac{1}{T} \sum_{t=0}^{T-1} \sum_{i=t_0}^{t}\frac{\Theta^2 c_E^2\eta^2}{q}\mathbb{E}[\|\zeta_i\|^2].
\end{align}
Note that $\sum_{j=t_0}^t\sum_{i=t_0}^j\mathbb{E}[\|\mathbf{\zeta}_i\|^2]\leq q\sum_{i=t_0}^t\mathbb{E}[\|\mathbf{\zeta}_i\|^2]$. This together with \eqref{proof-spb-2} yields
\begin{align}\label{proof-spb-3}
\frac{1}{T}\sum_{t=0}^{T-1}\mathbb{E}[\|V_t - \nabla f(X_t) \|^2] \leq  \frac{\Theta^2 c_E^2\eta^2}{T}\sum_{t=0}^{T-1}\mathbb{E}[\|\zeta_t\|^2].
\end{align}
Now combining Lemma \ref{G-bounded}, \eqref{proof-spb-0} and \eqref{proof-spb-3}, we have that:
\begin{equation}\label{proof-spb-4}
\frac{1}{T}\sum_{t=0}^{T-1}\mathbb{E}[\|G(\x_t, \nabla f(\x_t), \gamma)\|^2] \leq 14M_1^2\frac{\Delta_0}{\eta T} + 8M_1^2\gamma\Theta^2 c_E^2 \eta\frac{\Delta_0}{T}.
\end{equation}
Again, using Jensen's inequality and the convexity of $\|\cdot\|^2$, \eqref{proof-spb-4} gives:
\begin{align}
\label{proof-theorem2-0} & \left(\mathbb{E}\left[\frac{1}{T}\sum_{t=0}^{T-1}\|G(\x_t, \nabla f(\x_t), \gamma)\|\right]\right)^2 \\
\nonumber           \leq & \mathbb{E}\left[\left(\frac{1}{T}\sum_{t=0}^{T-1}\|G(\x_t, \nabla f(\x_t), \gamma)\|\right)^2\right] \\
\nonumber           \leq & \frac{1}{T}\sum_{t=0}^{T-1}\mathbb{E}[\|G(\x_t, \nabla f(\x_t), \gamma)\|^2] \\
\nonumber           \leq & 14M_1^2\frac{\Delta_0}{\eta T} + 8M_1^2\gamma\Theta^2 c_E^2 \eta\frac{\Delta_0}{T}.
\end{align}
Hence, we know that as long as
\begin{equation}\label{proof-spb-T-bound}
T \geq \left(14M_1^2\frac{\Delta_0}{\eta} + 8M_1^2\gamma\Theta^2 c_E^2 \eta\Delta_0\right)\epsilon^{-2},
\end{equation}
the right hand side of \eqref{proof-theorem2-0} is upper bounded by $\epsilon^2$, which implies that if index $\nu$ is uniformly sampled from $\{0,\ldots,T-1\}$, then
\[\mathbb{E}\left[\|G(\x_\nu, \nabla f(\x_\nu), \gamma)\|\right] \leq\epsilon.\]
That is, $X_\nu$ is an $\epsilon$-stochastic stationary point of problem \eqref{prob_1}. Equation \eqref{proof-spb-T-bound} then implies that the number of iterations needed by R-ProxSPB for obtaining an $\epsilon$-stochastic stationary point of problem \eqref{prob_1} in the finite-sum case is $T=\O(\epsilon^{-2})$.
Furthermore, the IFO complexity of R-ProxSPB under the finite-sum setting is:
\begin{equation}
    \lceil T/q \rceil\cdot |\mathcal{S}_t^1| + T\cdot |\mathcal{S}_t^2| \leq \frac{T+q}{q}n+T\sqrt{n} = \O(\sqrt{n}\epsilon^{-2}+ n),
\end{equation}
where the equality is due to $q=\sqrt{n}$.

\subsection{Online setting}

In the online case, $I\{|\mathcal{S}_t^1|<n\}=1$. Since $|\mathcal{S}_t^1|$ is the same for all $t$, we denote $s:= |\mathcal{S}_t^1|$. In this case, \eqref{telescopic_q_rewritten-1} reduces to
\begin{equation}\label{proof-spb-online-1}
    \sum_{i=t_0}^t\mathbb{E}[\|\mathbf{\zeta}_i\|^2] \leq \frac{\mathbb{E}[F(X_{(n_t-1)q})-F(X_{t+1})]}{\eta} + \frac{1}{2}\sum_{i=t_0}^t\frac{\sigma^2}{|\mathcal{S}_t^1|}.
\end{equation}
We calculate the telescoping sum for \eqref{proof-spb-online-1} for all length-$q$ epochs that $t+1 = q, 2q, \ldots, Kq$ ($K= \lfloor \frac{T}{q}\rfloor$) and the telescoping sum from $t=Kq$ to $T-1$. This gives:
\begin{align}\label{proof-spb-online-2}
\frac{1}{T}\sum_{t=0}^{T-1}\mathbb{E}[\|\mathbf{\zeta}_t\|^2] =  \frac{1}{T}\left(\sum_{t=0}^{Kq-1}\mathbb{E}[\|\mathbf{\zeta}_t\|^2]+\sum_{t=Kq}^{T-1}\mathbb{E}[\|\mathbf{\zeta}_t\|^2]\right)
\leq \frac{\Delta_0}{\eta T} + \frac{\sigma^2}{2s}.
\end{align}
Note that Lemma \ref{lemma-zhou} gives:
\[
     \mathbb{E}[\|V_t - \nabla f(X_t) \|^2] \leq \frac{\sigma^2}{s} +  \sum_{i=t_0}^{t-1}\frac{\Theta^2 c_E^2\eta^2}{|\mathcal{S}_t^2|}\mathbb{E}[\|\zeta_i\|^2],
\]
which further implies:
\begin{equation}\label{proof-spb-online-3}
    \frac{1}{T}\sum_{t=0}^{T-1}\mathbb{E}[\|V_t - \nabla f(X_t) \|^2] \leq \frac{\sigma^2}{s} + \frac{\Theta^2 c_E^2\eta^2}{T}\sum_{t=0}^{T-1}\mathbb{E}[\|\zeta_t\|^2].
\end{equation}
Now combining Lemma \ref{G-bounded}, \eqref{proof-spb-online-2} and \eqref{proof-spb-online-3}, we have that:
\begin{equation}\label{proof-spb-online-4}
\frac{1}{T}\sum_{t=0}^{T-1}\mathbb{E}[\|G(\x_t, \nabla f(\x_t), \gamma)\|^2] \leq \left(14M_1^2+8M_1^2\gamma\Theta^2 c_E^2\eta^2\right)\frac{\Delta_0}{\eta T} + M_1^2\frac{59\sigma^2}{5s}.
\end{equation}
Again, using Jensen's inequality and the convexity of $\|\cdot\|^2$, \eqref{proof-spb-online-4} gives:
\begin{align}
\label{proof-theorem2-1} & \left(\mathbb{E}\left[\frac{1}{T}\sum_{t=0}^{T-1}\|G(\x_t, \nabla f(\x_t), \gamma)\|\right]\right)^2 \\
\nonumber           \leq & \mathbb{E}\left[\left(\frac{1}{T}\sum_{t=0}^{T-1}\|G(\x_t, \nabla f(\x_t), \gamma)\|\right)^2\right] \\
\nonumber           \leq & \frac{1}{T}\sum_{t=0}^{T-1}\mathbb{E}[\|G(\x_t, \nabla f(\x_t), \gamma)\|^2] \\
\nonumber           \leq & \left(14M_1^2+8M_1^2\gamma\Theta^2 c_E^2\eta^2\right)\frac{\Delta_0}{\eta T} + M_1^2\frac{59\sigma^2}{5s}.
\end{align}
Now, by choosing
\begin{equation}\label{proof-spb-online-T-bound}
T = \left(\frac{2(14M_1^2+8M_1^2\gamma\Theta^2 c_E^2\eta^2)\Delta_0}{\eta}\right)\epsilon^{-2}, \quad \mbox{and} \quad s = \frac{118M_1^2\sigma^2}{5}\epsilon^{-2},
\end{equation}
we know that the right hand side of \eqref{proof-theorem2-1} is equal to $\epsilon^2$, which implies that if index $\nu$ is uniformly sampled from $\{0,\ldots,T-1\}$, then
\[\mathbb{E}\left[\|G(\x_\nu, \nabla f(\x_\nu), \gamma)\|\right] \leq\epsilon.\]
That is, $X_\nu$ is an $\epsilon$-stochastic stationary point of problem \eqref{prob_1}. Equation \eqref{proof-spb-online-T-bound} then implies that the number of iterations needed by R-ProxSPB for obtaining an $\epsilon$-stochastic stationary point of problem \eqref{prob_1} in the finite-sum case is $T=\O(\epsilon^{-2})$, and moreover, this needs to require the batch size $|\mathcal{S}_t^1|=s=\O(\epsilon^{-2})$ for all $t$. Furthermore, the IFO complexity of R-ProxSPB under the online setting is given by:
\[
    \lceil T/q\rceil\cdot |\mathcal{S}_t^1| + T\cdot |\mathcal{S}_t^2| \leq \frac{T+q}{q}\O(\epsilon^{-2})+T q = \O(\epsilon^{-3}),
\]
where the equality is due to $q = \epsilon^{-1}$.

\end{document}